\documentclass[11pt]{article}
\usepackage{amsmath}
\usepackage{algorithm}
 \usepackage{subcaption}
\usepackage[noend]{algpseudocode}
\usepackage{color}
\usepackage{bm}
\usepackage[authoryear]{natbib}
\usepackage[colorlinks,citecolor=blue,linkcolor=blue,urlcolor=blue]{hyperref}
\usepackage{latexsym, epsfig, amssymb, amsmath, amsthm, graphicx, mathrsfs, booktabs}
\usepackage{rotating, tabularx, setspace,paralist}
\usepackage{color}
\usepackage{graphicx}
\usepackage[OT1]{fontenc}
\usepackage{lscape}
\usepackage{multirow,multicol}

\usepackage{anysize}
\renewcommand{\baselinestretch}{1.2}
\marginsize{1.1in}{1.1in}{0.55in}{0.8in} %

\makeatletter
\def\singlespace{\def\baselinestretch{1}\@normalsize}

\renewcommand{\theequation}{
\arabic{equation}%
}

\makeatother

\newcommand{\ba}{\mbox{\bf a}}

\newcommand{\bu}{\mbox{\bf u}}
\newcommand{\bv}{\mbox{\bf v}}

\newcommand{\bx}{\mbox{\bf x}}
\newcommand{\by}{\mbox{\bf y}}

\newcommand{\bA}{\mbox{\bf A}}
\newcommand{\bB}{\mbox{\bf B}}

\newcommand{\bD}{\mbox{\bf D}}

\newcommand{\bH}{\mbox{\bf H}}
\newcommand{\bI}{\mbox{\bf I}}

\newcommand{\bQ}{\mbox{\bf Q}}
\newcommand{\bR}{\mbox{\bf R}}

\newcommand{\bM}{\mbox{\bf M}}
\newcommand{\bV}{\mbox{\bf V}}
\newcommand{\bW}{\mbox{\bf W}}
\newcommand{\bX}{\mbox{\bf X}}

\newcommand{\bone}{\mbox{\bf 1}}

\newcommand{\bPi}{\mbox{\boldmath $\Pi$}}

\newcommand{\bTheta}{\mbox{\boldmath $\Theta$}}

\newcommand{\bLam}{\mbox{\boldmath $\Lambda$}}

\newcommand{\bpi}{\mbox{\boldmath $\pi$}}

\newcommand{\bbell}{\mbox{\boldmath $\ell$ }}

\newcommand{\sbx}{\scriptsize{\bx}}
\newcommand{\sby}{\scriptsize{\by}}
\newcommand{\sbu}{\scriptsize{\bu}}
\newcommand{\sbv}{\scriptsize{\bv}}

\newcommand{\var}{\mathrm{var}}

\newcommand{\tr}{\mathrm{tr}}
\newcommand{\diag}{\mathrm{diag}}

\def\toD{\overset{\mathscr{D}}{\longrightarrow}}
\def\toP{\overset{\mathrm{P}}{\longrightarrow}}

\def\independenT#1#2{\mathrel{\setbox0\hbox{$#1#2$}%
\copy0\kern-\wd0\mkern4mu\box0}}

\def\beginn{\begin{eqnarray*}}
\def\endn{\end{eqnarray*}}
\def\beginy{\begin{eqnarray}}
\def\endy{\end{eqnarray}}
\def\begine{\begin{enumerate}}
\def\ende{\end{enumerate}}

\theoremstyle{plain}
\newtheorem{thm}{Theorem}

\newtheorem{lem}{Lemma}
\newtheorem{coro}{Corollary}
\newtheorem{rmk}{Remark}

\newtheorem{deff}{Definition}

\newtheorem{cond}{Condition}
\newcommand{\non}{\nonumber \\}
\newcommand{\bbA}{\mbox{\bf A}}

\newcommand{\bba}{\mbox{\bf a}}
\newcommand{\bbB}{\mbox{\bf B}}
\newcommand{\bbC}{\mbox{\bf C}}
\newcommand{\bbc}{\mbox{\bf c}}
\newcommand{\bbD}{\mbox{\bf D}}

\newcommand{\bbe}{\mbox{\bf e}}

\newcommand{\bbF}{\mbox{\bf F}}
\newcommand{\bbP}{\mbox{\bf P}}

\newcommand{\bbG}{\mbox{\bf G}}

\newcommand{\bbH}{\mbox{\bf H}}
\newcommand{\bbh}{\mbox{\bf h}}

\newcommand{\bbI}{\mbox{\bf I}}
\newcommand{\bbi}{\mbox{\bf i}}
\newcommand{\bbj}{\mbox{\bf j}}
\newcommand{\bbJ}{\mbox{\bf J}}

\newcommand{\bbM}{\mbox{\bf M}}

\newcommand{\bbQ}{\mbox{\bf Q}}

\newcommand{\bbR}{\mbox{\bf R}}

\newcommand{\bbs}{\mbox{\bf s}}

\newcommand{\bbS}{\mbox{\bf S}}

\newcommand{\bbu}{\mbox{\bf u}}
\newcommand{\bbV}{\mbox{\bf V}}
\newcommand{\bbv}{\mbox{\bf v}}

\newcommand{\bbW}{{\bf W}}

\newcommand{\bbX}{\mbox{\bf X}}

\newcommand{\bbx}{\mbox{\bf x}}

\newcommand{\bby}{\mbox{\bf y}}

\newcommand{\bbb}{\mbox{\bf b}}

\newcommand{\bbL}{\mbox{\bf L}}

\newcommand{\ep}{\ensuremath{\epsilon}}

\newcommand{\mb}[1]{\ensuremath{\mathbb #1}}

\renewcommand{\baselinestretch}{1.3}

\begin{document}

\title{Asymptotic Theory of Eigenvectors for Random Matrices with Diverging Spikes %
\thanks{Jianqing Fan is Frederick L. Moore '18 Professor of Finance, Department of Operations Research and Financial Engineering, Princeton University, Princeton, NJ 08544, USA (E-mail: jqfan@princeton.edu). %
Yingying Fan is Professor and Dean's Associate Professor in Business Administration, Data Sciences and Operations Department, Marshall School of Business, University of Southern California, Los Angeles, CA 90089 (E-mail: \textit{fanyingy@marshall.usc.edu}). %
Xiao Han is Postdoctoral Scholar, Data Sciences and Operations Department, Marshall School of Business, University of Southern California, Los Angeles, CA 90089 (E-mail: \textit{xhan011@e.ntu.edu.sg}). %
Jinchi Lv is Kenneth King Stonier Chair in Business Administration and Professor, Data Sciences and Operations Department, Marshall School of Business, University of Southern California, Los Angeles, CA 90089 (E-mail: \textit{jinchilv@marshall.usc.edu}). %
This work was supported by NIH grants R01-GM072611-14 and  1R01GM131407-01, NSF grants DMS-1662139, DMS-1712591, and DMS-1953356, NSF CAREER Award DMS-1150318, a grant from the Simons Foundation, and Adobe Data Science Research Award. %
The authors sincerely thank the Joint Editor, Associate Editor, and referees for their valuable comments that helped improve the paper substantially.
}
\date{October 13, 2020}
\author{Jianqing Fan$^1$, Yingying Fan$^2$, Xiao Han$^2$ and Jinchi Lv$^2$
\medskip\\
Princeton University$^1$ and University of Southern California$^2$
\\
} %
}

\maketitle

\begin{abstract}
Characterizing the asymptotic distributions of eigenvectors for large random matrices poses important challenges yet can provide useful insights into a range of statistical applications. To this end, in this paper we introduce a general framework of asymptotic theory of eigenvectors (ATE) for large spiked random matrices with diverging spikes and heterogeneous variances, and establish the asymptotic properties of the spiked eigenvectors and eigenvalues for the scenario of the generalized Wigner matrix noise. Under some mild regularity conditions, we provide the asymptotic expansions for the spiked eigenvalues and show that they are asymptotically normal after some normalization. For the spiked eigenvectors, we establish asymptotic expansions for the general linear combination and further show that it is asymptotically normal after some normalization, where the weight vector can be arbitrary. We also provide a more general asymptotic theory for the spiked eigenvectors using the bilinear form. Simulation studies verify the validity of our new theoretical results. Our family of models encompasses many popularly used ones such as the stochastic block models with or without overlapping communities for network analysis and the topic models for text analysis, and our general theory can be exploited for statistical inference in these large-scale applications.
\end{abstract}

\textit{Running title}: Asymptotic Theory for Eigenvectors


\section{Introduction} \label{Sec1}

The big data era has brought us a tremendous amount of both structured and unstructured data including networks and texts in many modern applications. For network and text data, we are often interested in learning the cluster and other structural information for the underlying network communities and text topics. In these large-scale applications, we are given a network data matrix or can create such a matrix by calculating some similarity measure between text documents, where each entry of the data matrix is binary indicating the absence or presence of a link, or continuous indicating the strength of similarity between each pair of nodes or documents. Such applications naturally give rise to random matrices that can be used to reveal interesting latent structures of networks and texts for effective predictions and recommendations.

Random matrix has been widely exploited to model the interactions among the nodes of a network for applications ranging from physics and social sciences to genomics and neuroscience. Random matrix theory (RMT) has a long history and was originated by Wigner in \cite{W55} for modeling the nucleon-nucleus interactions to understand the behavior of atomic nuclei and link the spacings of the levels of atomic nuclei to those of the eigenvalues of a random matrix. 
See, for example, \cite{B99} for a review of some classical technical tools such as the moment method and Stieltjes transform as well as some more recent developments on the RMT, and \cite{M04, T04, BS06} for detailed book-length accounts of the topic of random matrices.

There is a rich recent literature in mathematics on the asymptotic behaviors of eigenvalues and eigenvectors of random matrices \citep{erdHos2013,bourgade2018,bourgade2017,rudelson2016no,dekel2007}. The main challenge in many  RMT problems is caused by the strong dependence of eigenvalues if they are close to each other.  Using the terminologies in RMT,  four regimes are often of interests: bulk, subcritical edge, critical edge, and supercritical regimes. The first three regimes all have eigenvalues that are highly correlated with each other, and the last regime has weaker interactions among the  eigenvalues.  The last regime can be further divided into two categories according to the relative strength of spiked eigenvalues compared to noise, which can be roughly understood as the signal-to-noise ratio.   There have been exciting mathematical developments in the recent mathematical literature when the smallest spiked eigenvalue has the same order as the noise \citep{capitaine2018non, AJ13, bao2018singular}. Due to the challenge caused by constant signal-to-noise ratio, these existing results often take complicated forms and the asymptotic distributions depend generally on the noise matrix distribution in a complex way, limiting their practical usage to statisticians. In this paper, we consider the setting of \textit{diverging spikes} where the spiked eigenvalues are an order of magnitude larger than the noise level asymptotically.  Although mathematically easier, such random matrices are of great interests to statisticians, because many statistical applications such as network analysis and text analysis often fall into this regime. Yet there lack any formal results on the asymptotic expansions and asymptotic  distributions of spiked eigenvectors even in this setting. This motivates our study in this paper.

There is a larger literature on the limiting distributions of eigenvalues than eigenvectors in RMT. For instance, the limiting spectral distribution of the Wigner matrix was generalized by \cite{A67} and \cite{A71}. \cite{MP67} established the well-known Marchenko--Pastur law for the limiting spectral distribution of the sample covariance matrix including the Wishart matrix which plays an important role in statistical applications. 
In contrast, the asymptotic distribution of the largest nonspiked eigenvalue of Wigner matrix with Gaussian ensemble was revealed to be the Tracy--Widom law in \cite{TW94} and \cite{TW96}. More recent developments on the asymptotic distribution of the largest nonspiked eigenvalue include \cite{J01}, \cite{K07}, \cite{J08}, \cite{EY11}, and \cite{KY17}. See also \cite{F81}, \cite{BS05}, \cite{BY08}, \cite{AJ13}, \cite{PS13}, \cite{RJ13}, \cite{AJ14}, and \cite{WF17} for the asymptotic distributions of the spiked eigenvalues of various random matrices and sample covariance matrices. 
For the eigenvectors, \cite{capitaine2018non} and \cite{bao2018singular} established their asymptotic distributions, which depend on the specific distribution of the Wigner matrix in a complicated way, in the challenging setting of constant signal-to-noise ratio.
There is also a growing literature on the specific scenario and applications of large network matrices. To ensure consistency, \cite{JL09} proposed the sparse principal component analysis to reduce the noise accumulation in high-dimensional random matrices. See, for example, \cite{M01}, \cite{ST07}, \cite{BC09}, \cite{DKM11}, \cite{RC11}, \cite{L16}, \cite{A17}, \cite{JKL17}, \cite{CL18}, and \cite{V18}.

Matrix perturbation theory has been commonly used to characterize the deviations of empirical eigenvectors from the population ones, often under the average errors \citep{HJ12}. In contrast, recently \cite{FWZ18} and \cite{AFWZ17} investigated random matrices with low expected rank and provided a tight bound for the difference between the empirical eigenvector and some linear transformation of the population eigenvector through a delicate entrywise eigenvector analysis for the first-order approximation under the maximum norm. 
See also \cite{P07}, \cite{KK16}, \cite{KX16}, and \cite{WF17} for the asymptotics of empirical eigenstructure for large random matrices. Yet despite these endeavors, the precise asymptotic distributions of the eigenvectors for large spiked random matrices still remain largely unknown even for the case of Wigner matrix noise. Indeed characterizing the exact asymptotic distributions of  eigenvectors in such setting can provide useful insights into a range of statistical applications that involve the eigenspaces. In this sense, the asymptotic expansions and asymptotic distributions of eigenvectors established in this paper complement the existing work in the statistics literature.

The major contribution of this paper is introducing a general framework of asymptotic theory of eigenvectors (ATE) for large spiked random matrices with diverging spikes, where the mean matrix is low-rank and the noise matrix is the generalized Wigner matrix. 
The generalized Wigner matrix refers to a symmetric random matrix whose diagonal and upper diagonal entries are independent with zero mean, allowing for heterogeneous variances.  
Our family of models includes a variety of popularly used ones such as the stochastic block models with or without overlapping communities for network analysis and the topic models for text analysis. 
Under some mild regularity conditions, we establish the asymptotic expansions for the spiked eigenvalues and prove that they are asymptotically normal after some normalization. For the spiked eigenvectors, we provide asymptotic expansions for the general linear combination and further establish that it is asymptotically normal after some normalization for arbitrary weight vector. We also present a more general asymptotic theory for the spiked eigenvectors based on the bilinear form. To the best of our knowledge, these theoretical results are new to the literature. Our general theory can be exploited for statistical inference in a range of large-scale applications including network analysis and text analysis. For detailed comparisons with the literature, see Section \ref{Sec3.6}.

The rest of the paper is organized as follows. Section \ref{Sec2} presents the model setting and theoretical setup for ATE. We establish the asymptotic expansions and asymptotic distributions for the spiked eigenvectors as well as the asymptotic distributions for the spiked eigenvalues in Section \ref{Sec3}. Several specific statistical applications of our new asymptotic theory are discussed in Section \ref{new.Sec4}. Section \ref{Sec4} presents some numerical examples to demonstrate our theoretical results. We further provide a more general asymptotic theory extending the results from Section \ref{Sec3} using the bilinear form in Section \ref{Sec5}. Section \ref{Sec6} discusses some implications and extensions of our work. The proofs of main results are relegated to the Appendix. Additional technical details are provided in the Supplementary Material.

\section{Model setting and theoretical setup} \label{Sec2}

\subsection{Model setting} \label{Sec2.1}

As mentioned in the introduction, we focus on the class of large spiked symmetric random matrices with low-rank mean matrices and generalized Wigner matrices of noises. It is worth mentioning that our definition of the generalized Wigner matrix specified in Section \ref{Sec1} is broader than the conventional one in the classical RMT literature; see, for example,  \cite{Yau2012} for the formal mathematical definition with additional assumptions. To simplify the technical presentation, consider an $n \times n$ symmetric random matrix with the following structure
\begin{align} \label{model}
\bbX = \bbH + \bbW,
\end{align}
where $\bH=\bbV\bbD\bbV^T$ is a deterministic latent mean matrix of low rank structure, $\bbV=(\bbv_1,\cdots,\bbv_K)$ is an $n\times K$ orthonormal matrix of population eigenvectors $\bbv_k$'s with $\bbV^T\bbV=\bbI_K$, $\bbD=\diag(d_1,\cdots,d_K)$ is a diagonal matrix of population eigenvalues $d_k$'s with $|d_1| \geq \cdots \geq |d_K|>0$, 
 and $\bbW=(w_{ij})_{1\leq i,j \leq n}$ is a symmetric random matrix of independent noises on and above the diagonal with zero mean $\mathbb E w_{ij} = 0$, variances $\sigma_{ij}^2=\mathbb{E}w_{ij}^2$, and  $ \max_{1\leq i, j\leq n}|w_{ij}|\leq 1$. The rank $K$ of the mean part is assumed typically to be a smaller order of the random matrix size $n$, which is referred to as matrix dimensionality hereafter for convenience. The bounded assumption on $w_{ij}$ is made frequently for technical simplification and satisfied in many real applications such as network analysis and text analysis. 
 {It can be relaxed to
 $ \mathbb{E}|w_{ij}|^l\leq C^{l-2}\mathbb{E}|w_{ij}|^2,\  l\ge 2,$ with $C$ some positive constant, and all the proofs and results can carry through.}   

In practice, it is either matrix $\bbX$ or matrix $\bbX - \diag(\bbX)$ that is readily available to us, where $\diag(\cdot)$ denotes the diagonal part of a matrix. In the context of graphs, random matrix $\bbX$ characterizes the connectivity structure of a graph with self loops, while random matrix $\bbX - \diag(\bbX)$ corresponds to a graph without self loops. In the latter case, the observed data matrix can be decomposed as
\begin{align}\label{model-noloop}
\bbX - \diag(\bbX) = \bbH + \left[\bbW - \diag(\bbX)\right].
\end{align}
Observe that $\bbW - \diag(\bbX)$ has the similar structure as $\bbW$ in the sense of being symmetric and having bounded independent entries on and above the diagonal, by assuming that $\diag(\bX)$ has bounded entries for such a case. Thus models \eqref{model} and \eqref{model-noloop} share the same decomposition of a deterministic low rank matrix plus some symmetric noise matrix of bounded entries, which is roughly all we need for the theoretical framework and technical analysis.  For these reasons, to simplify the technical presentation we abuse slightly the notation by using $\bX$ and $\bbW$ to represent the observed data matrix and the latent noise matrix, respectively, in either model \eqref{model} or model \eqref{model-noloop}.  Therefore, throughout the paper the data matrix $\bbX$ may have diagonal entries all equal to zero and correspondingly the  noise matrix $\bbW$ may have a nonzero diagonal mean matrix, and our theory covers both cases.

In either of the two scenarios discussed above, we are interested in inferring the structural information in models \eqref{model} and \eqref{model-noloop}, which often boils down to the latent eigenstructure $(\bbD, \bbV)$.
Since both the eigenvector matrix $\bV$ and eigenvalue matrix $\bD$ are unavailable to us, we resort to the observable random data matrix $\bX$ for extracting the structural information. To this end, we conduct a spectral decomposition of $\bX$, and denote by $\lambda_1, \cdots , \lambda_n$ its eigenvalues and  $\widehat\bbv_1, \cdots, \widehat{\bbv}_n$ the corresponding eigenvectors. Without loss of generality, assume that $|\lambda_1| \geq \cdots \geq |\lambda_n|$ and denote by $\widehat\bbV=(\widehat \bbv_1,\cdots,\widehat \bbv_K)$ an $n \times K$ matrix of spiked eigenvectors. As mentioned before, we aim at investigating the precise asymptotic behavior of the spiked empirical eigenvalues $\lambda_1,\cdots, \lambda_K$ and spiked empirical eigenvectors $\widehat \bv_1, \cdots, \widehat \bv_K$ of data matrix $\bX$. It is worth mentioning that our definition of spikedness differs from the conventional one in that the underlying rank order depends on the magnitude of eigenvalues instead of the nonnegative eigenvalues that are usually assumed.

One concrete example is the stochastic block model (SBM), where the latent mean matrix $\bbH$ takes the form
$
\bbH = \bPi\bbP\bPi^T
$
with $\bPi = (\bpi_1,\cdots, \bpi_n)^T \in \mathbb{R}^{n\times K}$ a matrix of community membership vectors and $\bbP = (p_{kl})\in \mathbb{R}^{K\times K}$ a nonsingular matrix with $p_{kl}\in[0,1]$ for $1\le k,l\le K$. Here, for each $1 \leq i \leq n$, $\bpi_i \in \{\bbe_1,\cdots, \bbe_K\}$ with $\bbe_j\in\mathbb{R}^K$, $1\leq j \leq K$, a unit vector with the $k$th component being one and all other components being zero. It is well known that the community information of the SBM is encoded completely in the eigenstructure of the mean matrix $\bbH$, which  serves as one of our motivations for investigating the precise asymptotic distributions of the empirical eigenvectors and eigenvalues.

\subsection{Theoretical setup} \label{Sec2.2}

We first introduce some notation that will be used throughout the paper. We use $a\ll b$ to represent $a/b\rightarrow0$ as matrix size $n$ increases. We say that an event $\mathcal{E}_n$ holds with significant probability  if $\mathbb{P}(\mathcal{E}_n)=1-O(n^{-l})$ for some positive constant $l$ and sufficiently large $n$. For a matrix $\bA$, we use $\lambda_j(\bbA)$ to denote the $j$th largest eigenvalue in magnitude, and $\|\bbA\|_F$,  $\|\bA\|$, and $\|\bbA\|_{\infty}$ to denote the Frobenius norm, the spectral norm, and the matrix entrywise maximum norm, respectively. Denote by $\bA_{- k}$ the submatrix of $\bA$ formed by removing the $k$th column.  For any $n$-dimensional unit vector $\bx = (x_1,\cdots, x_n)^T$, let $d_{\sbx}=\|\bx\|_\infty$  represent the maximum norm of the vector.

We next introduce a definition that plays a key role in proving all asymptotic normality results in this paper.

\begin{deff}\label{def: 1}
A pair of unit vectors $(\bx, \by)$ of appropriate dimensions is said to satisfy the $\bW^{l}$-CLT condition for some positive integer $l$ if $\bx^T(\bW^l - \mb{E}\bW^l)\by$ is asymptotically standard normal after some normalization, where CLT refers to the central limit theorem.
\end{deff}

Lemmas \ref{0524-1} and \ref{0524-1h} below provide some sufficient conditions under which $(\bx, \by)$ can satisfy the
$\bW^{l}$-CLT condition defined in Definition \ref{def: 1} for $l=1$ and $2$, which is all we need for our technical analysis of asymptotic distributions.  In this paper, we apply these lemmas with either $\bx$ or $\by$ equal to $\bbv_k$. Therefore,  a sufficient condition for the results in our paper is that $\|\bbv_k\|_{\infty}$ is small enough. 	

\begin{lem}\label{0524-1}
Assume that $n$-dimensional unit vectors $\bbx$ and $\bby$ satisfy
\begin{equation}\label{eq: lemma1-assp}
		\|\bbx\|_{\infty}\|\bby\|_{\infty} \ll 
\left[\var(\bbx^T\bbW\bby)\right]^{1/2} = s_n.
		\end{equation}
Then $\bbx^T\bbW\bby$ satisfies the Lyapunov condition for CLT and we have $(\bbx^T\bbW\bby-\mathbb{E}\bbx^T\bbW\bby)/s_n \toD N(0,1)$ as $n \rightarrow \infty$, which entails that $(\bx, \by)$ satisfies the
$\bW^{l}$-CLT condition with $l = 1$.
\end{lem}


To introduce $\bW^{2}$-CLT, for any given 
unit vectors $\bbx = (x_1,\cdots, x_n)^T$ and $\bby = (y_1,\cdots, y_n)^T$,   we denote respectively $s^2_{\sbx,\sby}$ and $\kappa_{\sbx,\sby}$ the mean and variance of the random variable
\begin{align}\label{0524.3a}
&\sum_{1\leq k,i\leq n,\, k\le i}\sigma^2_{ki}\Big[\sum_{1\leq l<k\leq n}w_{il}(x_ky_l+y_kx_l)+\sum_{1\leq l<i\leq n}w_{kl}(x_iy_l+y_ix_l)\nonumber\\
&\quad+(1-\delta_{ki})\mathbb{E}w_{ii}(x_iy_k+x_ky_i)\Big]^2+2\sum_{1\leq k,i\leq n,\,k\le i}\gamma_{ki}(x_ky_k+x_iy_i) \nonumber\\
&\quad \times \Big[\sum_{1\leq l<k\leq n}w_{il}(x_ky_l+y_kx_l)+\sum_{1\leq l<i\leq n}w_{kl}(x_iy_l+y_ix_l)\nonumber\\
&\quad+(1-\delta_{ki})\mathbb{E}w_{ii}(x_iy_k+x_ky_i)\Big]+\sum_{1\leq k,i\leq n,\,k\le i}\kappa_{ki}(x_ky_k+x_iy_i)^2,
\end{align}
where $\gamma_{ki}=\mathbb{E}w_{ki}^3$  
and $\kappa_{ki}=\mathbb{E}(w_{ki}^2-\sigma_{ki}^2)^2$ for $k\neq i$,   $\gamma_{kk}=2(\mathbb{E}\omega_{kk}^3-\sigma_{kk}^2\mathbb{E}\omega_{kk})$, $\kappa_{kk}=4\mathbb{E}(\omega_{kk}^2-\sigma_{kk}^2)^2$ with $\omega_{kk}=2^{-1}w_{kk}$, $\sigma_{kk}^2=\mathbb{E}\omega_{kk}^2$, and $\delta_{ki} = 1$ when $k = i$ and 0 otherwise. It is worth mentioning that the random variable given in (\ref{0524.3a}) coincides with the one defined in (\ref{0524.3}) in Section \ref{SecB.2} of Supplementary Material, which is simply the conditional variance of random variable $\bbx^T(\bbW^2-\mathbb{E}\bbW^2)\bby$ given in (\ref{0524.2}) when expressed as a sum of martingale differences with respect to a suitably defined $\sigma$-algebra; see Section \ref{SecB.2} for more technical details and the precise expressions for $s^2_{\sbx,\sby}$ and $\kappa_{\sbx,\sby}$ given in (\ref{0525.1}) and (\ref{0929.2h}), respectively.

\begin{lem}\label{0524-1h}
Assume that $n$-dimensional unit vectors $\bbx$ and $\bby$ satisfy $\|\bbx\|_\infty \|\bby\|_{\infty}\rightarrow 0$,
	$\kappa_{\sbx,\sby}^{1/4}\ll s_{\sbx,\sby}$,  and $s_{\sbx,\sby}\rightarrow \infty$. Then we have
	$[\bbx^T(\bbW^2-\mathbb{E}\bbW^2)\bby]/s_{\sbx,\sby} \toD N(0,1)$ as $n\rightarrow \infty$, which entails that $(\bx, \by)$ satisfies the
	$\bW^{2}$-CLT condition.
\end{lem}

\begin{rmk}\label{2rvm1}
		To provide more insights into the conditions of Lemmas \ref{0524-1} and \ref{0524-1h}, we discuss the special case of standard Wigner matrix where $\sigma_{ij}^2=p(1-p)$ with $p$ the expected value of entries of $\bbX$. Then $s_n^2 := \var(\bbx^T\bbW\bby)\in [p(1-p),2p(1-p)]$ and condition \eqref{eq: lemma1-assp} in Lemma \ref{0524-1} reduces to
		$$\|\bbx\|_{\infty}\|\bby\|_{\infty} \ll \sqrt{\var(\bbx\bbW\bby)}\sim \sqrt{p(1-p)}.$$
Moreover, \eqref{neweq025} in the Supplementary Material ensures that Lemma \ref{0524-1h} holds under the following sufficient conditions
		\begin{equation}
		\|\bbx\|_\infty \|\bby\|_{\infty}\rightarrow 0, \  n^{3/2}p(1-p)\|\bbx\|_\infty^2\|\bby\|_{\infty}^2\rightarrow 0, \text{ and } p(1-p) n\rightarrow \infty.
		\end{equation}
Thus if either $\|\bbx\|_\infty $ or $\|\bby\|_{\infty}$ is small enough, both lemmas hold. Indeed in this scenario, direct calculations show that $\bbs^2_{\bbx,\bby}\sim np(1-p)$. 	
	\end{rmk}

We see from Lemmas \ref{0524-1} and \ref{0524-1h} that the
$\bW^{l}$-CLT condition defined in Definition \ref{def: 1} can indeed be satisfied under some mild regularity conditions. In particular, Definition \ref{def: 1} is important to our technical analysis since to establish the asymptotic normality of the spiked eigenvectors and spiked eigenvalues, we first need to expand the target to the form of $\bx^T(\bW^l - \mb{E}\bW^l)\by$ with $l$ some positive integer plus some small order term, and then the asymptotic normality follows naturally if $(\bx, \by)$ satisfies the $\bW^{l}$-CLT condition. To facilitate our technical presentation, let us introduce some further notation. For any $t\neq 0$ and given 
matrices $\bM_1$ and $\bM_2$ of appropriate dimensions, we define the function
\begin{equation}\label{0619.0}
\mathcal{R}(\bM_1,\bbM_2,t)=-\sum_{{l=0,\,l\neq 1}}^L t^{-(l+1)}\bbM_1^T\mathbb{E}\bbW^l\bbM_2,
\end{equation}
where $L$ is some sufficiently large positive integer that will be specified later in our technical analysis. For each $1\leq k \leq K$, any given 
matrices $\bM_1$ and $\bM_2$ of appropriate dimensions, and  $n$-dimensional vector $\bu$, we further define functions
\begin{align}\label{0619.1}
&\mathcal{P}(\bM_1,\bM_2,t)=t \mathcal{R}(\bM_1,\bM_2,t), \quad
\mathcal{\widetilde P}_{k,t}=\left[t^2 (A_{\bbv_k,k,t}/t)'\right]^{-1},\\
&
\bbb_{\sbu,k,t}=\bbu-\bbV_{-k}\left[(\bbD_{-k})^{-1}+\mathcal{R}(\bbV_{-k},\bbV_{-k},t)\right]^{-1}\mathcal{R}(\bbu,\bbV_{-k},t)^T, \label{0619.3}
\end{align}
where $\bbD_{-k}$ denotes the submatrix of the diagonal matrix $\bbD$ by removing the $k$th row and $k$th column,
\begin{equation}\label{0619.2}
A_{\sbu,k,t}=\mathcal{P}(\bbu,\bbv_k,t)-\mathcal{P}(\bbu,\bbV_{-k},t)\left[t(\bbD_{-k})^{-1}+\mathcal{P}(\bbV_{-k},\bbV_{-k},t)\right]^{-1}\mathcal{P}(\bbV_{-k},\bbv_k,t),
\end{equation}
$(\cdot)'$ denotes the derivative with respect to scalar $t$ or complex variable $z$ throughout the paper, and the rest of notation is the same as introduced before.

\section{Asymptotic distributions of spiked eigenvectors} \label{Sec3}

\subsection{Technical conditions} \label{Sec3.1}

To facilitate our technical analysis, we need some basic regularity conditions.

\begin{cond}\label{as1}
	Assume that $\alpha_n = \|\mathbb{E}(\bbW-\mathbb{E}\bbW)^2\|^{1/2} \rightarrow \infty$ as $n\rightarrow \infty$.
\end{cond}

\begin{cond}\label{as3}
	There exists a positive constant $c_0<1$ such that $\min\{|d_i|/|d_{j}|: 1\le i<j\le K+1,\, d_i\neq -d_j\}\ge 1+c_0$. In addition, either of the following two conditions holds:
	\begin{itemize}
	\item[i)] $|d_K|/(n^{\ep}\alpha_n)\rightarrow \infty$ with some small positive constant $\ep$,
	\item[ii)] $\max_{i,j}\var(w_{ij})\le (c^2_1\alpha_n^2)/n$ and $|d_K|>c\alpha_n\log n$ with some constants $c_1\ge 1$ and $c>4c_1(1+2^{-1}c_0)$.
	\end{itemize}
\end{cond}

\begin{cond}\label{as7}
	It holds that $|d_1|=O(|d_K|)$, $|d_K|\sigma_{\min}/\alpha_n\rightarrow \infty$,  $\|\bbv_k\|_{\infty}^2/\sigma_{\min}\rightarrow 0$, $\alpha_n^4\|\bbv_k\|_{\infty}^4/\\
	(\sqrt n\sigma_{\min}^2)\rightarrow 0$, and $ \sigma_{\min}^2n\rightarrow \infty$, where $\sigma_{\min}= \{\min_{1\leq i,j\leq n,\,i\neq j}\mathbb{E} w_{ij}^2\}^{1/2}$.
\end{cond}
Conditions \ref{as1}--\ref{as3} are needed in all our Theorems \ref{0619-1}--\ref{0518-1} and imposed for our general model \eqref{model}, including the specific case of sparse models. In contrast, condition \ref{as7} is required \textit{only} for Theorem \ref{1107-1} under some specific models with dense structures such as the stochastic block models with or without overlapping communities.

Condition \ref{as1} restricts essentially the sparsity level of the random matrix (e.g., given by a network). Note that it follows easily from $\max_{1\le i,j\le n}|w_{ij}|\le 1$ that $\alpha_n\le n^{1/2}$. It is a rather mild condition that can be satisfied by very sparse networks. For example, if $\mathbb{E}w^2_{12}=\cdots=\mathbb{E}w^2_{1\lfloor\log n \rfloor}=1/2$ and the other $w_{1j}$'s are equal to zero, then we have $\alpha_n^2\ge 2^{-1} \log n\rightarrow \infty$. Many network models in the literature satisfy this condition; see, for example, \cite{JKL17}, \cite{L16}, and \cite{YEJ15}.

Condition \ref{as3} requires that the spiked population eigenvalues of the mean matrix $\bH$ (in the diagonal matrix $\bbD$) are simple and there is enough gap between the eigenvalues. The constant $c_0$ can be replaced by some $o(1)$ term and our theoretical results can still be proved with more delicate derivations. This requirement ensures that we can obtain higher order expansions of the general linear combination for each empirical eigenvector precisely. Otherwise if there exist some eigenvalues such that $d_i=d_{i+1}$, then $\widehat \bbv_i$ and $\widehat \bbv_{i+1}$ are generally no longer identifiable so we cannot derive clear asymptotic expansions for them; see also \cite{AFWZ17} for related discussions. Condition \ref{as3} also requires a gap between $\alpha_n$ and $|d_K|$. Since parameter $\alpha_n$ reflects the strength of the noise matrix $\bbW$, it requires essentially the signal part $\bbH$ to dominate the noise part $\bbW$ with some asymptotic rate. Similar  condition is used commonly in the network literature; see, for instance, \cite{AFWZ17} and \cite{JKL17}.

Condition \ref{as7} restricts our attention to some specific dense network models. In particular, $|d_1|=O(|d_K|)$ assumes that the eigenvalues in  $\bbD$ share the same order. The other assumptions in Condition \ref{as7} require essentially that the minimum variance of the off-diagonal entries of $\bbW$ cannot tend to zero too fast, which is used only to establish a more simplified theory under the more restrictive model; see Theorem \ref{1107-1}.

\subsection{Asymptotic distributions of spiked eigenvalues} \label{Sec3.2}

We first present the asymptotic expansions and CLT for the spiked empirical eigenvalues $\lambda_1, \cdots, \lambda_K$.  For each $1 \leq k \leq K$, denote by $t_k$ the solution to equation
\begin{align}\label{0515.3.1}
	f_k(z)&=1+d_k\Big\{\mathcal{R}(\bbv_k,\bbv_k,z)-\mathcal{R}(\bbv_k,\bbV_{-k},z)\left[(\bbD_{-k})^{-1}+\mathcal{R}(\bbV_{-k},\bbV_{-k},z)\right]^{-1} \nonumber \\
	&\quad \times \mathcal{R}(\bbV_{-k},\bbv_k,z)\Big\} =0
\end{align}
when restricted to the interval $z\in [a_k, b_k ]$,
where $$a_k=\begin{cases}
d_k/(1+2^{-1}c_0), &d_k>0\cr (1+2^{-1}c_0)d_k, &d_k<0
\end{cases} \quad \text{ and  } \quad b_k=\begin{cases}
(1+2^{-1}c_0)d_k, &d_k>0\cr d_k/(1+2^{-1}c_0), &d_k<0
\end{cases}.$$
The following lemma characterizes the properties of the population quantity $t_k$'s defined in (\ref{0515.3.1}): It is unique and the asymptotic mean of $\lambda_k$.

\begin{lem}\label{lem: define-t}
	Equation \eqref{0515.3.1} has a unique solution in the interval $z\in [a_k, b_k ]$
	and thus $t_k$'s are well defined. Moreover, for each $1 \leq k \leq K$  we have $t_k/d_k \to 1$ as $n\rightarrow \infty$.
\end{lem}

It is seen from Lemma \ref{lem: define-t} that when the matrix size $n$ is large enough, the values of $t_k$ and $d_k$ are very close to each other. The following theorem establishes the asymptotic expansions and CLT for $\lambda_k$ and reveals that $t_k$ is in fact its asymptotic mean.

\begin{thm}\label{0619-1}
	Under Conditions \ref{as1}--\ref{as3}, for each $1 \leq k \leq K$ we have
	\begin{eqnarray}\label{0927.3h}
		\lambda_k-t_k=\bbv_k^T\bbW\bbv_k+O_p(\alpha_n d_k^{-1}).
	\end{eqnarray}
	Moreover, if $\var(\bbv_k^T\bbW\bbv_k)\gg \alpha^2_n d^{-2}_k$ and the pair of vectors $(\bv_k, \bv_k)$ satisfies the $\bW^1$-CLT condition,
	then we have 
	{\begin{eqnarray}\label{0408.1}
			\frac{\lambda_k-t_k-\mathbb{E}\bbv_k^T\bbW\bbv_k}{\left[\var(\bbv_k^T\bbW\bbv_k)\right]^{1/2}}\toD N(0,1).
	\end{eqnarray}}
\end{thm}

\cite{CDF12} and \cite{AJ14} established the joint distribution of the spiked eigenvalues for the deformed Wigner matrix in different settings than ours. \cite{CDF12} assumed that $\mathbb{E}w_{ii}^2=1/2$ and $\mathbb{E}w_{ij}^2=1$ for $i\neq j$, 
while \cite{AJ14} assumed that $\mathbb{E}w_{ij}^2=1$ for all $(i,j)$. Under their model settings, the smallest spiked eigenvalue $|d_K|$ and the noise level $\alpha_n$ are of the same order, and as a result, their asymptotic distributions depend on the distributions of the Wigner matrix. In contrast, our Theorem \ref{0619-1} is proved in the setting of diverging spikes. Thanks to the stronger signal-to-noise ratio, the noise matrix contributes to the distributions of the spiked eigenvalues in Theorem \ref{0619-1} in a global way, allowing for more heterogeneity in the variances of entries of the noise matrix $\bW$.

Theorem \ref{0619-1}   requires that $(\bv_k, \bv_k)$ satisfies the  $\bW^1$-CLT condition and  $\var(\bbv_k^T\bbW\bbv_k)\gg \alpha^2_n d^{-2}_k$.  To gain some insights into these two conditions, we will provide some sufficient conditions for such assumptions. Let us consider the specific case of $\sigma_{\min}>0$, that is, the generalized Wigner matrix $\bW$ is nonsparse.
We will show that as long as
\begin{align}\label{e001}
	\|\bbv_k\|_{\infty}^2 \sigma_{\min}^{-1} \rightarrow 0 \  \text{ and } \  \sigma_{\min}\gg \alpha_n |d_k|^{-1},
\end{align}
the aforementioned two conditions in Theorem \ref{0619-1} hold. We first verify the $\bW^1$-CLT condition.	By Lemma \ref{0524-1}, a sufficient condition for $(\bv_k, \bv_k)$ to satisfy the $\bW^1$-CLT condition is that
\begin{align}
	\|\bv_k\|_\infty^2 \ll \left[\var(\bbv_k^T\bbW\bbv_k)\right]^{1/2}=\left[\mathbb{E}(\bbv_k^T\bbW\bbv_k-\mathbb{E}\bbv_k^T\bbW\bbv_k)^2\right]^{1/2}. \label{eq003}
\end{align}
Observe that it follows from $\sum_{1\leq i \leq n}(\bv_k)_i^2 = \bv_k^T\bv_k = 1$ and $\sum_{1 \leq i \leq n}(\bv_k)_i^4 \leq \|\bv_k\|_\infty^2  \leq 1 $ that
\begin{align}
	\nonumber &\left[\mathbb{E}(\bbv_k^T\bbW\bbv_k-\mathbb{E}\bbv_k^T\bbW\bbv_k)^2\right]^{1/2} \geq\left[2\sum_{1\leq i,j\leq n,\,i\neq j} \sigma_{ij}^2(\bbv_k)_i^2(\bbv_k)_j^2\right]^{1/2}\\
	\nonumber
	&\quad \geq \sigma_{\min}\left[2\sum_{1\leq i,j\leq n,\,i\neq j}(\bbv_k)_i^2(\bbv_k)_j^2\right]^{1/2} = \sigma_{\min}\left[2-2\sum_{1\leq i\leq n}(\bbv_k)_i^4\right]^{1/2} \\
	& \quad \geq \sigma_{\min}\left(2-2\|\bv_k\|_\infty^2\right)^{1/2}, \label{eq002}
\end{align}
where $(\bbv_k)_i$ stands for the $i$th component of vector $\bv_k$.	The assumption $\|\bv_k\|_\infty^2 \sigma_{\min}^{-1}\rightarrow 0$ in \eqref{e001} together with \eqref{eq002} ensures \eqref{eq003}, which consequently entails that $(\bv_k, \bv_k)$ satisfies the $\bW^1$-CLT condition.

We next check the condition $\var(\bbv_k^T\bbW\bbv_k)=\mathbb{E}(\bbv_k^T\bbW\bbv_k-\mathbb{E}\bbv_k^T\bbW\bbv_k)^2\gg \alpha^2_n d^{-2}_k$. It follows directly from \eqref{eq002} that this condition holds under \eqref{e001}. In fact, since Condition \ref{as3} guarantees that $\alpha_n/|d_k|$ asymptotically vanishes, 
the assumption $\sigma_{\min}\gg \alpha_n |d_k|^{-1}$ can be very mild.
In particular, for the Wigner matrix $\bW$ with $\sigma_{ij}\equiv 1$ for all $1\leq i,j\leq n$, it holds that
\begin{align}\label{1203.1}
	{\mathbb{E}(\bbv_k^T\bbW\bbv_k-\mathbb{E}\bbv_k^T\bbW\bbv_k)^2} =2.
\end{align}
Thus the condition of  ${\mathbb{E}(\bbv_k^T\bbW\bbv_k-\mathbb{E}\bbv_k^T\bbW\bbv_k)^2}\gg \alpha_n^2 d_k^{-2}$ reduces to that of $\alpha_n^2 d_k^{-2} \ll 1$, which is guaranteed to hold under Condition \ref{as3}.

We also would like to point out that one potential application of the new results in Theorem \ref{0619-1}  is determining the number of spiked eigenvalues, which in the network models reduces to determining  the number of non-overlapping (or possibly overlapping) communities or clusters.

\subsection{Asymptotic distributions of spiked eigenvectors} \label{Sec3.3}

We now present the asymptotic distributions of the spiked empirical eigenvectors $\widehat{\bv}_k$ for $1 \leq k \leq K$. To this end, we will first establish the asymptotic expansions and CLT for the bilinear form
$$
\bbx^T\widehat\bbv_k\widehat\bbv_k^T\bby
$$
with $1 \leq k \leq K$, where $\bbx, \bby \in \mathbb{R}^n$ are two arbitrary non-random unit vectors. Then by setting {$\by = \bv_k$}, we can establish the asymptotic expansions and CLT for the general linear combination $\bx^T\widehat\bv_k$. Although the limiting distribution of  the bilinear form $\bbx^T\widehat\bbv_k\widehat\bbv_k^T\bby$ is the theoretical foundation for establishing the limiting distribution of the general linear combination $\bx^T\widehat\bv_k$, due to the technical complexities we will defer the theorems summarizing the limiting distribution of $\bbx^T\widehat\bbv_k\widehat\bbv_k^T\bby$ to a later technical section (i.e., Section \ref{Sec5}), and present only the results for $\bx^T\widehat\bv_k$ in this section.  This should not harm the flow of the paper. For readers who are also interested in our technical proofs, they can refer to Section \ref{Sec5} for more technical details; otherwise it is safe to skip that technical section. For each $1 \leq k \leq K$, let us choose the direction of $\widehat \bbv_k$ such that $\bbv_k^T\widehat\bbv_k\ge 0$ for the theoretical derivations, which is always possible after a sign change when needed.

\begin{thm}\label{0609-1}
Under Conditions \ref{as1}--\ref{as3},  for each $1 \leq k \leq K$ we have the following properties:

	1)	If the unit vector $\bbu$ satisfies that $|\bbu^T\bbv_k|\in [0,1)$  and
	$\alpha_n^{-2}d_k^2\var[(\bbb^T_{\sbu,k,t_k}-\bbu^T\bbv_k\bbv_k^T)\bbW\bbv_k]\rightarrow \infty$, then it holds that
	{\begin{align}\label{0609.1}\nonumber
			& t_k\left(\bbu^T\widehat\bbv_k+A_{\sbu,k,t_k}\mathcal{\widetilde P}_{k,t_k}^{1/2}\right)=(\bbb^T_{\sbu,k,t_k}-\bbu^T\bbv_k\bbv_k^T)\bbW\bbv_k\\
			&\quad+o_p\left(\left\{\var[(\bbb^T_{\sbu,k,t_k}-\bbu^T\bbv_k\bbv_k^T)\bbW\bbv_k]\right\}^{1/2}\right),
	\end{align}}
	where the asymptotic mean has the expansion $A_{\sbu,k,t_k}\mathcal{\widetilde P}_{k,t_k}^{1/2}=-\bbu^T\bbv_k+O(\alpha_n^2 d_k^{-2})$. Furthermore,  if $(\bbb_{\sbu,k,t_k}-\bbv_k\bbv_k^T\bbu, \bbv_k)$ satisfies the $\bW^1$-CLT condition, then it holds that  {$$\frac{t_k\left(\bbu^T\widehat\bbv_k+A_{\sbu,k,t_k}\mathcal{\widetilde P}_{k,t_k}^{1/2}\right)-\mathbb{E}\left[(\bbb^T_{\sbu,k,t_k}-\bbu^T\bbv_k\bbv_k^T)\bbW\bbv_k\right]}{\left\{\var\left[(\bbb^T_{\sbu,k,t_k}-\bbu^T\bbv_k\bbv_k^T)\bbW\bbv_k\right]\right\}^{1/2}}
		\toD N(0, 1).
		$$}
	
	2)	If
	$(\alpha_n^{-4}d_k^2+1)\var(\bbv_k^T\bbW^2\bbv_k)\rightarrow \infty$, then it holds that
	\begin{equation}\label{0609.2}
		2t_k^2\left(\bbv_k^T\widehat\bbv_k+A_{\sbv_k,k,t_k}\mathcal{\widetilde P}_{k,t_k}^{1/2}\right)=-\bbv_k^T\left(\bbW^2-\mathbb{E}\bbW^2\right)\bbv_k + o_p\left\{\left[\emph{var}(\bbv_k^T\bbW^2\bbv_k)\right]^{1/2}\right\},
	\end{equation}
	where the asymptotic mean has the expansion $A_{\sbv_k,k,t_k}\mathcal{\widetilde P}_{k,t_k}^{1/2}=-1+2^{-1}t_k^{-2}\bbv_k^T\mathbb{E}\bbW^2\bbv_k+O(\alpha_n^3 d_k^{-3})$. Furthermore,  if $(\bv_k, \bbv_k)$ satisfies the $\bW^2$-CLT condition, then it holds that  $$\frac{2t_k^2\left(\bbv_k^T\widehat\bbv_k+A_{\sbv_k,k,t_k}\mathcal{\widetilde P}_{k,t_k}^{1/2}\right)}{\left[\var\left(\bbv_k^T\bbW^2\bbv_k\right)\right]^{1/2}}\toD N(0, 1).$$
\end{thm}

The two parts of Theorem \ref{0609-1} correspond to two different cases when $\var(\bu^T\widehat{\bv}_k)$ can be of different magnitude. To understand this, note that for large enough matrix size $n$, we have $|t_K| \gg \alpha_n$ by Condition \ref{as3} and Lemma \ref{lem: define-t}. In view of \eqref{0609.2}, the asymptotic variance of $\bv_k^T\widehat\bv_k$ is equal to
$
\var(2^{-1} t^{-2}_k\bbv_k^T\bbW^2\bbv_k).
$
In contrast, in light of \eqref{0609.1}, the asymptotic variance of $\bu^T\widehat\bv_k$ with $|\bu^T\bv_k| \in [0,1)$ is equal to
$
\var\left[t_k^{-1}(\bbb^T_{\bu,k,t_k}-\bbu^T\bbv_k\bbv_k^T)\bbW\bbv_k\right].
$
Let us consider a specific case when  $\var[(\bbb^T_{\bu,k,t_k}-\bbu^T\bbv_k\bbv_k^T)\bbW\bbv_k]\sim 1$.  By Lemma \ref{0426-1} in Section \ref{Sec5},
we have
$$
\var\left(2^{-1} t^{-2}_k\bbv_k^T\bbW^2\bbv_k \right)=O\left(\alpha_n^2t^{-4}_k\right)\ll \var\left[t_k^{-1}(\bbb^T_{u,k,t_k}-\bbu^T\bbv_k\bbv_k^T)\bbW\bbv_k\right]=O\left(t_k^{-2}\right).
$$
This shows that the above two cases can be very different in the magnitude for the asymptotic variance of $\bu^T\widehat{\bv}_k$ and thus should be analyzed separately.

To gain some insights into why $\bv_k^T\widehat\bv_k$ has smaller variance, let us consider the simple case  of $K=1$. Then in view of our technical arguments, it holds that
\begin{eqnarray}\label{0615.3}
	\bbv_1^T\widehat\bbv_1\widehat\bbv_1^T\bbv_1&=&-\frac{1}{2\pi i}\oint_{\Omega_1}\frac{\bbv_1^T(\bbW-z\bbI)^{-1}\bbv_1}{1+d_1\bbv_1^T(\bbW-z\bbI)^{-1}\bbv_1} dz\non
	&=&-\frac{1}{2\pi i}\oint_{\Omega_1}\frac{1}{\left[\bbv_1^T(\bbW-z\bbI)^{-1}\bbv_1\right]^{-1}+d_1} dz,
\end{eqnarray}
where $i = \sqrt{-1}$  is associated with the complex integrals represents the imaginary unit and the line integrals are taken over the contour $\Omega_1$ that is centered at $(a_1+b_1)/2$  with radius $c_0|d_1|/2$. Then we can see that the population eigenvalue $d_1$ is enclosed by the  contour $\Omega_1$.
By the Taylor expansion, we can show that with significant probability, $$\left[\bbv_1^T(\bbW-z\bbI)^{-1}\bbv_1\right]^{-1}=-z-\bbv_1^T\bbW\bbv_1+O\left(|z|^{-1} \alpha^2_n\log n\right).$$
Substituting the above expansion into (\ref{0615.3}) results in
\begin{align}\label{0615.4}
	&-\frac{1}{2\pi i}\oint_{\Omega_1}\frac{1}{\left[\bbv_1^T(\bbW-z\bbI)^{-1}\bbv_1\right]^{-1}+d_1} dz=-\frac{1}{2\pi i}\oint_{\Omega_1}\frac{1}{d_1-z-\bbv_1^T\bbW\bbv_1+O(|z|^{-1} \alpha^2_n\log n)} dz\non
	&\quad=-\frac{1}{2\pi i}\oint_{\Omega_1}\frac{1}{d_1-z}- \frac{1}{2\pi i}\oint_{\Omega_1}\frac{\bbv_1^T\bbW\bbv_1}{(d_1-z)^2}dz+O\left({d_1^{-2}\alpha^2_n\log n}\right)\non
	&\quad=-\frac{1}{2\pi i}\oint_{\Omega_1}\frac{1}{d_1-z}+O\left(d_1^{-2}\alpha^2_n\log n\right)
\end{align}
with significant probability. Thus the asymptotic distribution of $\bv_1^T\widehat\bv_1\widehat\bbv_1^T\bbv_1$ is determined by $O_p(d_1^{-2}\alpha^2_n\log n)$, which has no contribution from  $\bbv_1^T\bbW\bbv_1$. On the other hand, our technical analysis for $\bu^T\widehat\bv_1\widehat\bv_1\bbv_1$ (which is much more complicated and can be found in the technical proofs section) reveals that the dominating term is $\bbu^T\bbW\bbv_1$ when $\bu\neq \bv_1$ or $-\bv_1$. This explains why we need to treat differently the two cases of $\bu$ close to or far away from $\bbv_1$.

\subsection{A more specific structure and an application} \label{Sec3.4}

Theorem \ref{0609-1} in Section \ref{Sec3.3} provides some general sufficient conditions to ensure the asymptotic normality for the spiked empirical eigenvectors. Under some simplified but stronger assumptions in Condition \ref{as7}, the same results on the empirical eigenvectors and  eigenvalues continue to hold. Note that the stochastic block models with non-overlapping or overlapping communities can both be included as specific cases of our theoretical analysis. As mentioned before, we choose the direction of $\widehat \bbv_k$ such that $\bbv_k^T\widehat\bbv_k\ge 0$ for each  $1 \leq k \leq K$.

\begin{thm}\label{1107-1}
	Under Conditions \ref{as1}--\ref{as7}, for each $1 \leq k \leq K$ we have the following properties:
	
	1) 	(Eigenvalues) It holds that
	{$$	\frac{\lambda_k-t_k-\mathbb{E}\bbv_k^T\bbW\bbv_k}{\left[\mathbb{E}(\bbv_k^T\bbW\bbv_k-\mathbb{E}\bbv_k^T\bbW\bbv_k)^2\right]^{1/2}}\toD N(0,1).$$}
	
	2) (Eigenvectors) If the unit vector $\bu$ satisfies that $\sigma_{\min}^{-1}\|\bbv_k(\bbb^T_{\sbu,k,t_k}-\bbu^T\bbv_k\bbv_k^T)\|_{\infty}\rightarrow 0$ and $|\bbu^T\bbv_k|\in [0,1-\ep]$ for some positive constant $\ep$, then it holds that
	\begin{eqnarray}\label{2rv10}
		{\frac{t_k\left(\bbu^T\widehat\bbv_k+A_{\sbu,k,t_k}\mathcal{\widetilde P}_{k,t_k}^{1/2}\right)-\mathbb{E}\left[(\bbb^T_{\sbu,k,t_k}-\bbu^T\bbv_k\bbv_k^T)\bbW\bbv_k\right]}{\left\{\var\left[(\bbb^T_{\sbu,k,t_k}-\bbu^T\bbv_k\bbv_k^T)\bbW\bbv_k\right]\right\}^{1/2}}
			\toD N(0, 1).} \label{0609.1h} 	\end{eqnarray}
	Moreover, it also holds that
	\begin{eqnarray}\label{0927.4h}
		\frac{2t_k^2\left(\bbv_k^T\widehat\bbv_k+A_{\sbv_k,k,t_k}\mathcal{\widetilde P}_{k,t_k}^{1/2}\right)}{\left[\var(\bbv_k^T\bbW^2\bbv_k)\right]^{1/2}}&\toD& N(0,1). \label{0609.2h}
	\end{eqnarray}
\end{thm}

Theorem \ref{0609-1} also gives the asymptotic expansions for the asymptotic mean term $A_{\sbu,k,t_k}\mathcal{\widetilde P}_{k,t_k}^{1/2}$. It is seen that if $|d_K|$ diverges to infinity much faster than $\alpha_n^2$, then the $O(\cdot)$ terms in the asymptotic expansions of the mean become smaller order terms and thus the following corollary follows immediately from Theorem \ref{1107-1}.
\begin{coro}\label{1109-1}
	Assume that Conditions \ref{as1}--\ref{as7} hold. For each $1 \leq k \leq K$, if the unit vector $\bu$ satisfies that  $|\bbu^T\bbv_k|\in [0,1-\ep]$ for some positive constant $\ep$ and $\alpha_n^{-4}d_k^2\var[(\bbb^T_{\sbu,k,t_k}-\bbu^T\bbv_k\bbv_k^T)\bbW\bbv_k]\rightarrow \infty$, then we have
	{\begin{eqnarray}\label{2rv6} \frac{t_k\left(\bbu^T\widehat\bbv_k-\bbu^T\bbv_k\right)-\mathbb{E}\left[(\bbb^T_{\sbu,k,t_k}-\bbu^T\bbv_k\bbv_k^T)\bbW\bbv_k\right]}{\left\{\var\left[(\bbb^T_{\sbu,k,t_k}-\bbu^T\bbv_k\bbv_k^T)\bbW\bbv_k\right]\right\}^{1/2}}&\toD& N(0,1).
	\end{eqnarray}}
	Moreover, if  $\alpha_n^{-6}d_k^3\var(\bbv_k^T\bbW^2\bbv_k)\rightarrow \infty$ then we have
	\begin{eqnarray}\label{2rv1}
		\frac{2t_k^2\left(\bbv_k^T\widehat\bbv_k-1\right)+\bbv_k^T\mathbb{E}\bbW^2\bbv_k}{\left[\var(\bbv_k^T\bbW^2\bbv_k)\right]^{1/2}}&\toD& N(0,1).
	\end{eqnarray}
\end{coro}

Theorem \ref{1107-1} includes the stochastic block model as a specific case. If $\bbX$ is the affinity matrix from a stochastic block model with $K$ non-overlapping communities and the size of each  community is of the same order $O(n)$, then it holds that  $\|\bbv_k\|_{\infty}=O(n^{-1/2})$, $d_K = O(n)$, {$\alpha_n\leq n^{1/2}$, and $\alpha_n\|\bbv_k\|_\infty \leq O(1)$}. Thus Condition \ref{as7} can be satisfied as long as $\sigma_{\min} \gg n^{-1/4}$, leading to the asymptotic normalities in Theorem \ref{1107-1}.

Our Theorem \ref{1107-1} also covers the stochastic block models with overlapping communities.  For example, the following network model was considered in \cite{YEJ15}
\begin{equation}\label{0104.1}
	\mathbb{E}\bbX  =\bTheta\bPi\bbP\bPi^T\bTheta^T,
\end{equation}
where $\bTheta$ is an $n \times n$ diagonal degree heterogeneity matrix, $\bPi$ is an $n \times K$ community membership matrix, and $\bbP$ is a $K\times K$ nonsingular irreducible matrix with unit diagonal entries. Observe that the above model has low-rank mean matrix and thus can be connected to our general form of eigendecomposition $\mathbb{E}\bX = \bH=\bbV\bbD\bbV^T$. If the spiked eigenvalues and spiked eigenvectors satisfy that $|d_k| = O(n)$  and $\|\bbv_k\|_\infty = O(n^{-1/2})$ for all $1 \leq k \leq K$, then Condition \ref{as7} can be satisfied when $\sigma_{\min} \gg n^{-1/4}$. Consequently, the asymptotic normalities in Theorem \ref{1107-1} can hold.



\subsection{Proofs architecture} \label{Sec3.5}

The key mathematical tools are from complex analysis and random matrix theory. 
At a high level, our technical proofs consist of four steps. First, we apply Cauchy's residue theorem to represent the desired bilinear form $\bbx^T\widehat\bbv_k\widehat\bbv_k^T\bby$ with $1\leq k\leq K$ as a complex integral over a contour for a functional of the Green function associated with the original random matrix $\bX = \bH + \bW$. It is worth mentioning that such an approach was used before to study the asymptotic distributions for linear combinations of eigenvectors in the setting of covariance matrix estimation for the case of i.i.d. Gaussian random matrix coupled with linear dependency. Second, we reduce the problem to one that involves a functional of the new Green function associated with only the noise part $\bW$ by extracting the spiked part. Such a step enables us to conduct precise high order asymptotic expansions. Third, we conduct delicate high order Taylor expansions for the noise part using new Green function corresponding to the noise part. In this step, we apply the asymptotic expansion directly to the evaluated complex integral over the contour instead of an expansion of the integrand. Such a new way of asymptotic expansion is crucial to our study. Fourth, we bound the variance of $\bbx^T(\bbW^l-\mathbb{E}\bbW^l)\bby$ using delicate random matrix techniques. In contrast to just counting the number of certain paths in a graph used in classical random matrix theory literature, we need to carefully bound the individual contributions toward the quantity $\alpha_n = \|\mathbb{E}(\bbW-\mathbb{E}\bbW)^2\|^{1/2}$; otherwise simple counting leads to rather loose upper bound.

\subsection{Comparisons with the statistics literature} \label{Sec3.6}

In a related work, \cite{T18} established the CLT for the entries of eigenvectors of a random adjacency matrix. Our work differs significantly from theirs in at least four important aspects. First,  \cite{T18} assumed a prior distribution on the mean adjacency matrix, while we assume a deterministic mean matrix. {As a result, the asymptotic variance in \cite{T18} is determined by the prior distribution and is the same for each entry of an eigenvector,} while in our paper the CLT for different entries of an eigenvector can be different and the asymptotic variance depends on all entries of the eigenvector. 
While \cite{T18} also provided the conditional CLT under the setting of the stochastic block model, their result conditions on just one node.
Second, our model is much more general than that in \cite{T18} in that the spiked eigenvalues can have different orders and different signs.
Third, we establish the CLT for the general linear combinations of the components of normalized eigenvectors and the CLT for eigenvalues, while \cite{T18} proved the CLT for the rows of $\bLam^{1/2} \widehat\bbV^T$, where $\bLam \in \mathbb{R}^{K\times K}$ is the diagonal matrix formed by $K$ spiked eigenvalues of the adjacency matrix and $\widehat\bbV=(\widehat\bbv_1,\cdots,\widehat\bbv_K)$ is the matrix collecting the corresponding eigenvectors of the adjacency matrix.
Fourth, through a dedicated analysis of the higher order expansion for the general linear combination $\bu^T\widehat\bbv_k$, we uncover an interesting phase transition phenomenon that the  limiting distribution of  $\bu^T\widehat\bbv_k$ is different when the deterministic weight vector $\bu$ is close to or far away from $\bv_k$ (modulo the sign), which is new to the literature.

\cite{WF17} proved the asymptotic distribution of the linear form $\bbv_i^T\widehat \bbv_k$ with $1 \leq i,k \leq K$, where $\bbv_i$'s and $\widehat \bbv_k$'s are the spiked population and empirical eigenvectors for some \textit{covariance matrix}, respectively. Their asymptotic normality results cover the case of $\bbv_1^T\widehat\bbv_1$ when $K=1$, and $\bbv_i^T\widehat\bbv_k$ for $1 \leq i, k \leq K$ with $i\neq k$ when $K\geq 2$.
Similarly, \cite{KK16} considered the sample \textit{covariance matrix} under the Gaussian distribution assumption, and derived the asymptotic expansion of the bilinear form $\bbx^T\widehat\bbv_k\widehat\bbv_k^T\bby$, where $\bx, \by$ are two deterministic unit vectors. They also obtained the asymptotic distribution of $\bbx^T\widehat\bbv_k$. Different from \cite{WF17} and \cite{KK16}, in this paper we establish the asymptotic distribution for the \textit{general} linear combination $\bu^T\widehat \bbv_k$ for the large structured symmetric random matrix from model \eqref{model}  under fairly weak regularity conditions. Our proof techniques differ from those in \cite{WF17} and \cite{KK16}, and are also distinct from most of existing ones in the literature.


\section{Statistical applications} \label{new.Sec4}

The new asymptotic expansions and asymptotic distributions of spiked eigenvectors and eigenvalues established in Section \ref{Sec3} have many natural statistical applications. Next we discuss three specific ones. See also \cite{FFHL19b} for another application on testing the node membership profiles in network models. 

\subsection{Detecting the existence of clustering power} \label{new.Sec4.1}

One potential application of Theorem \ref{1107-1} is to improve the results on community detection under model setting \eqref{0104.1}.  Spectral methods have been used popularly in the literature for recovering the memberships of nodes in network models. For example, applying the $K$-means clustering algorithm to the $K$ spiked eigenvectors calculated from the adjacency matrix has been a prevalent method for inferring the memberships of nodes. However, it may not be true that all these $K$ eigenvectors are useful for clustering. For example, if eigenvector $\bbv_k={\bf1}/\sqrt{n}$, then it has zero clustering power and should be dropped in the $K$-means clustering algorithm. This is especially important in large networks because including a useless high-dimensional eigenvector may significantly increase the noise in clustering. Theorem  \ref{1107-1} suggests that we can test the hypothesis $H_0: \bbv_k ={\bf1}/\sqrt{n}$ using the test statistic $\widehat \bbv_k^T{\bf1}$. Then with the aid of Theorem \ref{1107-1}, the asymptotic null distribution can be established and the critical value can be calculated. This naturally suggests a method for selecting important eigenvectors in community detection. 


\subsection{Detecting the existence of denser subgraph} \label{sec: app2}

Another application of Theorem \ref{1107-1} is to detect the existence of a denser community in a given random graph, the same problem as studied in \cite{arias2014community} and \cite{verzelen2015community}. 
Specifically, assume that the data matrix $\bbX = (x_{ij})$ is a symmetric adjacency matrix with independent Bernoulli entries on and above the diagonal. Let $\bbH = \mathbb E[\bbX]$ be the mean adjacency matrix. Consider the following null and alternative hypotheses
$$
H_0: \bbH = p\bone\bone^T \quad \text{ vs. } \quad  H_1: \bbH = p\bone\bone^T + (q-p)\bbell\bbell^T,
$$
where $\bbell$ is the vector with the first $n_1$ entries being 1 and all remaining entries being 0, and $q\in (p,1]$. It can be seen that under the alternative hypothesis, there is a denser subgraph and $q$ measures the connectivity of nodes within it. \cite{arias2014community} and \cite{verzelen2015community} proposed tests for the above hypothesis in the setting of $n_1=o(n)$. We focus on the same setting and in addition assume that $n^{-1} \ll p< q $ and $q \sim p$. We next discuss how to exploit our Theorem \ref{1107-1} to test the same hypothesis.

Under the null hypothesis,  a natural estimator of $p$ is given by $\widehat p=\frac{1}{n(n-1)}\sum_{1\le i\le j\le n}x_{ij}$. Moreover, direct calculations show that
\begin{equation}\label{2rv2}
\bbv_1^T\mathbb{E}\bbW^2\bbv_1=np(1-p) \ \text{ and } \  \var(\bbv_1^T\bbW^2\bbv_1)=p(1-p)\left[2(n-1)+p^3+(1-p)^3\right].
\end{equation}
Thus the mean and variance of $\bbv_1^T\bbW^2\bbv_1$ in \eqref{2rv2} can be estimated as
\begin{equation}\label{2rv4-a}
n\widehat p(1-\widehat p) \ \text{ and } \ \widehat p(1-\widehat p)\left[2(n-1)+\widehat{p}^3+(1-\widehat p)^3\right],
\end{equation}
receptively.
In view of \eqref{2rv1} in Corollary \ref{1109-1}, since $\bbv_1 = n^{-1/2}\bone$ under the null hypothesis $H_0: \bbH = p\bone\bone^T$, a natural test statistic for testing $H_0: \bbH = p\bone\bone^T$ is given by 
$$
T_n = \frac{2\lambda_1^2\left(n^{-1/2}\bone^T\widehat\bbv_1-1\right)+n\widehat p(1-\widehat p)}{\left[\widehat p(1-\widehat p)\left[2(n-1)+\widehat{p}^3+(1-\widehat p)^3\right]\right]^{1/2}}.
$$

It can be seen that since $\lambda_1 \approx t_1$ (see Lemma \ref{lem: define-t}), the asymptotic null distribution of $T_n$ is expected to be $N(0,1)$ by resorting to \eqref{2rv1} in Corollary \ref{1109-1}. On the other hand, under the alternative hypothesis, since  the leading eigenvector differs from $n^{-1/2}\bone$, the term $n^{-1/2}\bone^T\widehat\bbv_1-1$ in the numerator of $T_n$ is expected to take some negative value, and thus $T_n$ is expected to have different asymptotic behavior than $N(0,1)$. In fact, we provide the proof sketch in Section \ref{sec:prof-app2} of  Supplementary Material on the asymptotic null and alternative distributions. In particular, we show that the asymptotic null distribution of $T_n$ is $N(0,1)$, and if  $\frac{n_1^2(q-p)^2}{np}+\frac{n_1^2(q-p)}{n}\gg 1$, then $T_n \rightarrow -\infty$ with asymptotic probability one under the alternative hypothesis.

\subsection{Rank inference} \label{subsec:rankInf}

Our theory can also be applied to statistical testing on the true rank $K$ of the mean matrix $\bbH$. Rank inference is an important problem in many high-dimensional network applications. See, for example, \cite{L16}, \cite{CL18}, and \cite{li2020network}, and the importance of the problem discussed therein. Consider the following hypotheses
$$H_0: K=K_0 \quad \text{ vs. } \quad H_1:K>K_0,$$
where $K_0$ is some prespecified positive integer satisfying $K_0\leq K$.
Define
\begin{align}\label{2rv7}
\widehat w_{ij}& = 
x_{ij}-\sum_{k=1}^{K_0}\lambda_k\bbe_i^T\widehat\bbv_k\widehat\bbv_k^T\bbe_j\non
&=w_{ij}-\sum_{k=1}^K\left[\lambda_k\bbe_i^T\widehat\bbv_k\widehat\bbv_k^T\bbe_j-d_k\bbe_i^T\bbv_k\bbv_k^T\bbe_j\right] + \sum_{k=K_0+1}^{K}\lambda_k\bbe_i^T\widehat\bbv_k\widehat\bbv_k^T\bbe_j.
\end{align}
Under the null hypothesis $H_0: K=K_0$, the last term in \eqref{2rv7} disappears and we can obtain the asympttoic expansion of $\widehat w_{ij}$ around $w_{ij}$ explicitly by an application of Theorems \ref{0619-1} and \ref{0609-1}. Then under some additional regularity conditions, it is expected that $\widehat w_{ij}$ is close to $w_{ij}$.  By the independence of $w_{ii}$, $i=1,\cdots, n$, it holds that
$$
\frac{\sum_{i=1}^n w_{ii}}{\sqrt{\sum_{i=1}^n w^2_{ii}}}\toD N(0,1) \ \text{ as } n\rightarrow \infty.
$$
Since $\widehat w_{ii} \approx w_{ii}$ under the null hypothesis, the following asymptotic distribution is expected to hold as well
\begin{equation}\label{test-Tn}
T_n := \frac{\sum_{i=1}^n\widehat w_{ii}}{\sqrt{\sum_{i=1}^n\widehat w^2_{ii}}}\toD N(0,1).
\end{equation}
This naturally suggests a statistical test based on statistic $T_n$ for testing $H_0: K=K_0$. Under the alternative hypothesis, since $\widehat w_{ij}$ contains  the smallest $K-K_0$ spiked eigenvalues and the corresponding eigenvectors, its asymptotic behavior is expected to be different, and consequently, the test can have nontrivial power. In fact, a more sophisticated version of this test constructed based on the off-diagonal entries of $\widehat \bbW$ was investigated recently in \cite{han2019universal}.  

The above asymptotic distribution can also be used to construct confidence intervals for the rank $K$. To understand this, note that $T_n$ defined in \eqref{test-Tn} is a function of $K_0$. Thus an immediate idea for the $100(1-\alpha)\%$ confidence interval construction is to identify all $K_0$ such that the corresponding $T_n$ falls into the range of $[-\Phi^{-1}(1-\alpha), \Phi^{-1}(1-\alpha)]$, where $\Phi^{-1}(\cdot)$ is the inverse distribution function of the standard normal. Similar ideas can also be exploited to construct confidence intervals for other parameters in network models.

\section{Simulation studies}  \label{Sec4}

In this section, we use simulation studies to verify the validity of our theoretical results.  We consider the stochastic block model with $K = 2$ communities. Assume that the number of nodes is $n$, the first $n/2$ nodes belong to the first community, and the rest belong to the second one. Then the adjacency matrix $\bbX$ has the mean structure $\mathbb{E}\bbX=\bbH=\bbA\bbR\bbA^T$, where $\bbR$ is a $2\times 2$ matrix of the connectivity probabilities,  and $\bbA=(\bba_1,\bba_2)\in \mathbb{R}^{n\times 2}$ with $\bba_1=n^{-1/2}(\textbf{1}^T,\textbf{0}^T)^T$ and $\bba_2=n^{-1/2}(\textbf{0}^T,\textbf{1}^T)^T$, where $\textbf{0}, \textbf{1}\in \mathbb{R}^{n/2}$ are vectors of zeros and ones, respectively. It is worth mentioning that $\bbA\bbR\bbA^T$ is \textit{not} the eigendecomposition of the mean matrix $\bbH$, which is why we use different notation than that in model \eqref{model}.

For the connectivity probability matrix $\bR$, we consider the structure
\[ \bbR= r\left(
\begin{array}{ccc}
2 & 1 \\
1 & 2\\
\end{array}
\right), \]
where parameter $r$ takes $6$ different values $0.02$, $0.05$, $0.1$, $0.2$, $0.3$, and $0.4$. A similar model was considered in  \cite{AFWZ17} and \cite{L16}. For the connectivity matrix $\bX$, we simulate its entries on and above the diagonal as independent Bernoulli random variables with means given by the corresponding entries in the mean matrix $\bH$, and set the entries below the diagonal to be the same as the corresponding ones above the diagonal. We choose the number of nodes as $n=3000$ and repeat the simulations for $10,000$ times.

\begin{figure*}[t!]
	\centering
	\begin{subfigure}[b]{0.5\textwidth}
		\centering
		\includegraphics[scale=0.4]{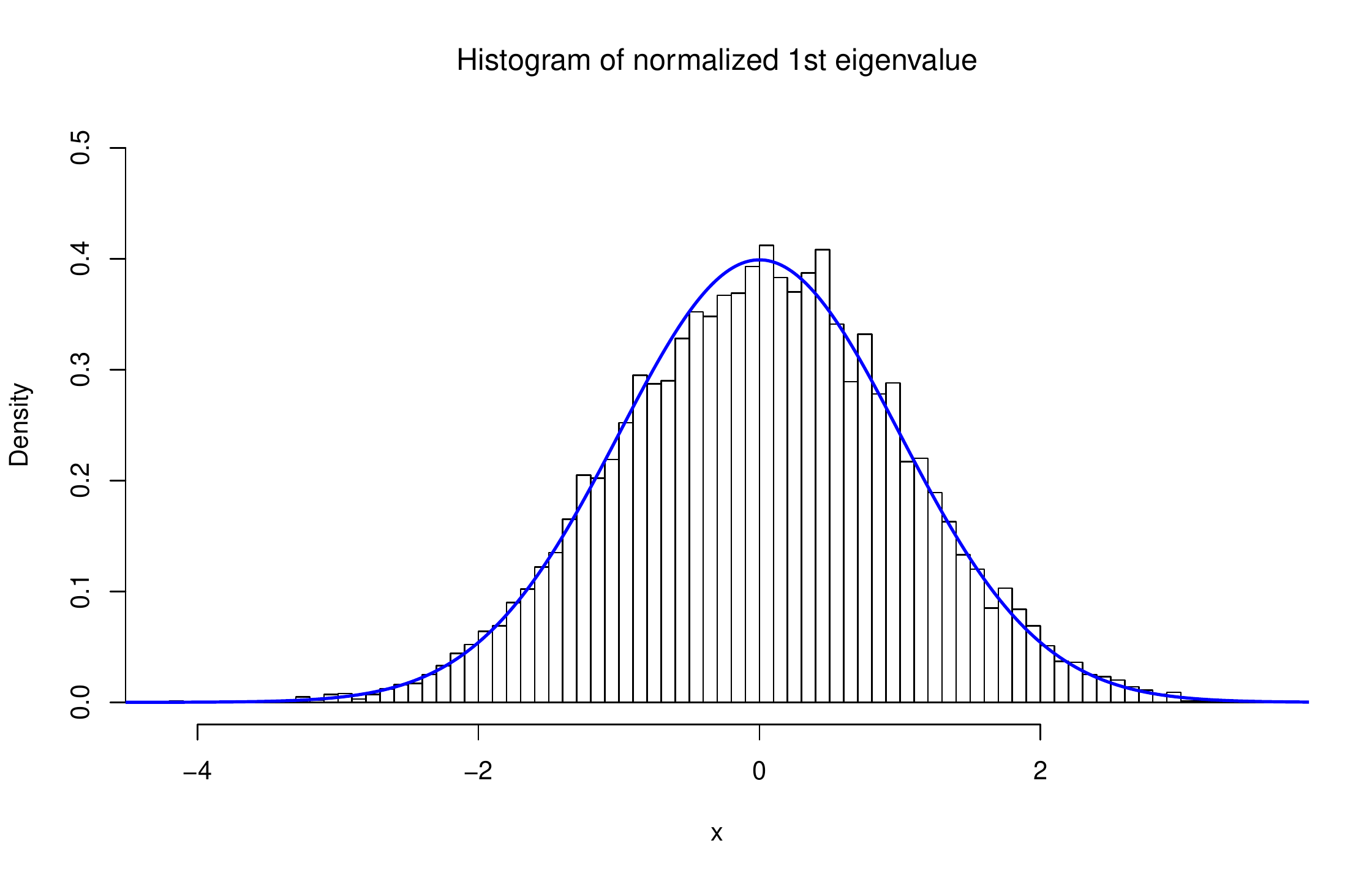}
	\end{subfigure}%
	~
	\begin{subfigure}[b]{0.5\textwidth}
		\centering
		\includegraphics[scale=0.4]{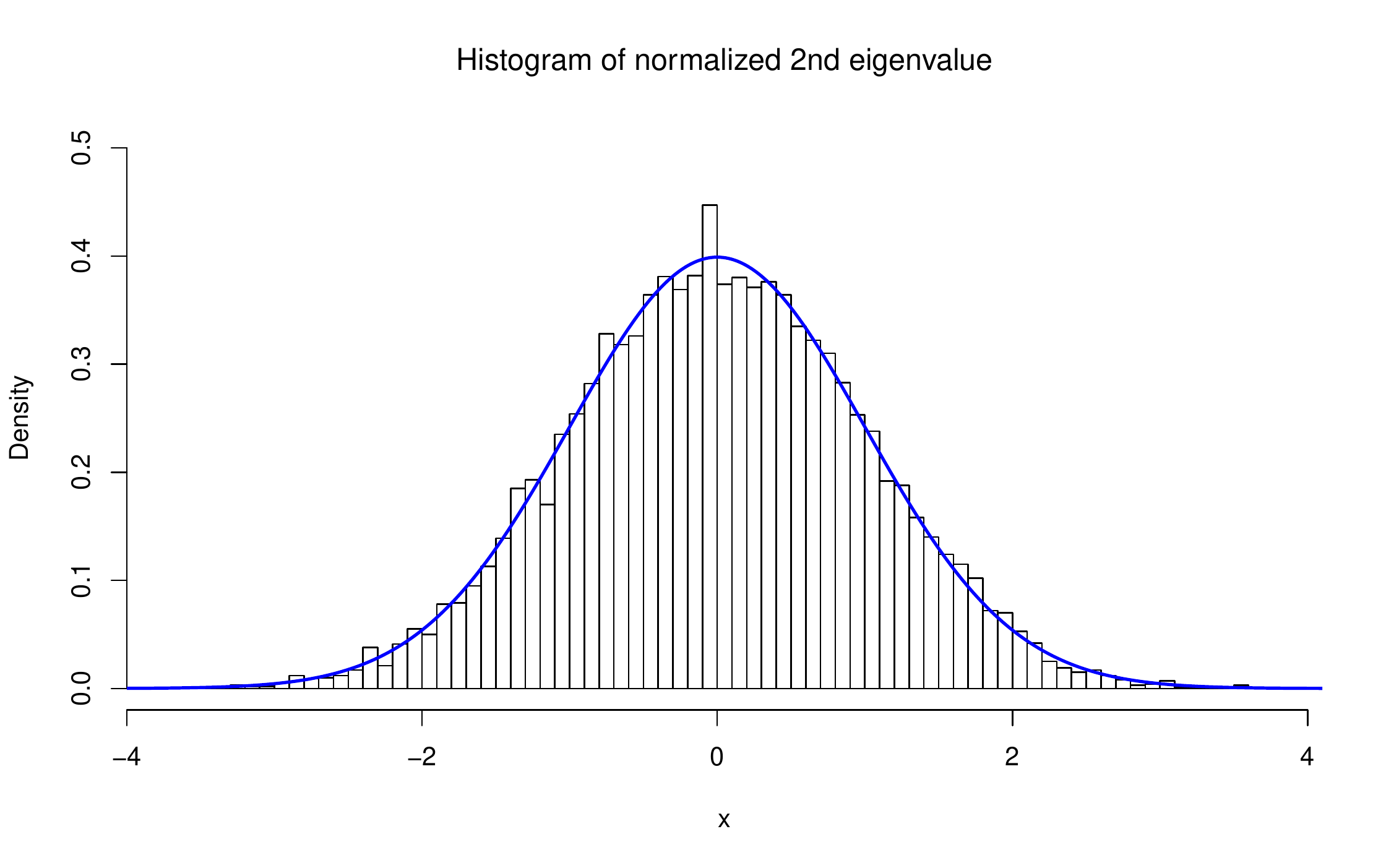}
	\end{subfigure}
	\caption{Histograms of the normalized eigenvalues \eqref{eq: stand-eigenvalue} when $r = 0.4$, with the blue curves representing the standard normal density. Left panel: the first eigenvalue; right panel: the second eigenvalue.}\label{fig:eig-val}
\end{figure*}

To verify our theoretical results, for each simulated connectivity matrix $\bX$ we calculate its eigenvalues and corresponding eigenvectors.
For the eigenvalues, we compare the empirical distribution of
\begin{eqnarray}\label{eq: stand-eigenvalue}
	\frac{\lambda_k-t_k}{\left[\var(\bbv_k^T\bbW\bbv_k)\right]^{1/2}}
\end{eqnarray}
with the standard normal distribution, where $t_k$ is the solution to equation \eqref{0515.3.1}. The exact expression of $\mathcal{R}(\bbv_k,\bbV_{-k},z)[(\bbD_{-k})^{-1}+\mathcal{R}(\bbV_{-k},\bbV_{-k},z)]^{-1}\mathcal{R}(\bbV_{-k},\bbv_k,z)$ in \eqref{0515.3.1} is complicated. Since this term is much smaller than $\mathcal{R}(\bbv_k,\bbv_k,z)$, we can calculate an approximation of $t_k$ by solving the equation
\begin{equation}\label{0408.2}
	1+d_k\mathcal{R}(\bbv_k,\bbv_k,z)=0
\end{equation}
using the Newton--Raphson method. Guided by the theoretical derivations, we use $L=4$ in the asymptotic expansion of  $\mathcal{R}(\bbx,\bby,t)$ in \eqref{0619.0} for all of our simulation examples. Tables \ref{tab1}--\ref{tab5} summarize the means and standard deviations of \eqref{eq: stand-eigenvalue} with $k=1$ and $2$ calculated from the 10,000 repetitions as well as the p-values from the  Anderson--Darling (AD) test for the normality. Figure \ref{fig:eig-val} presents the histograms of the normalized first and second eigenvalues (i.e., \eqref{eq: stand-eigenvalue} when $k=1$ and 2) from the 10,000 repetitions.

\begin{figure}
	\centering
	\includegraphics[width=16.5cm,height=10cm]{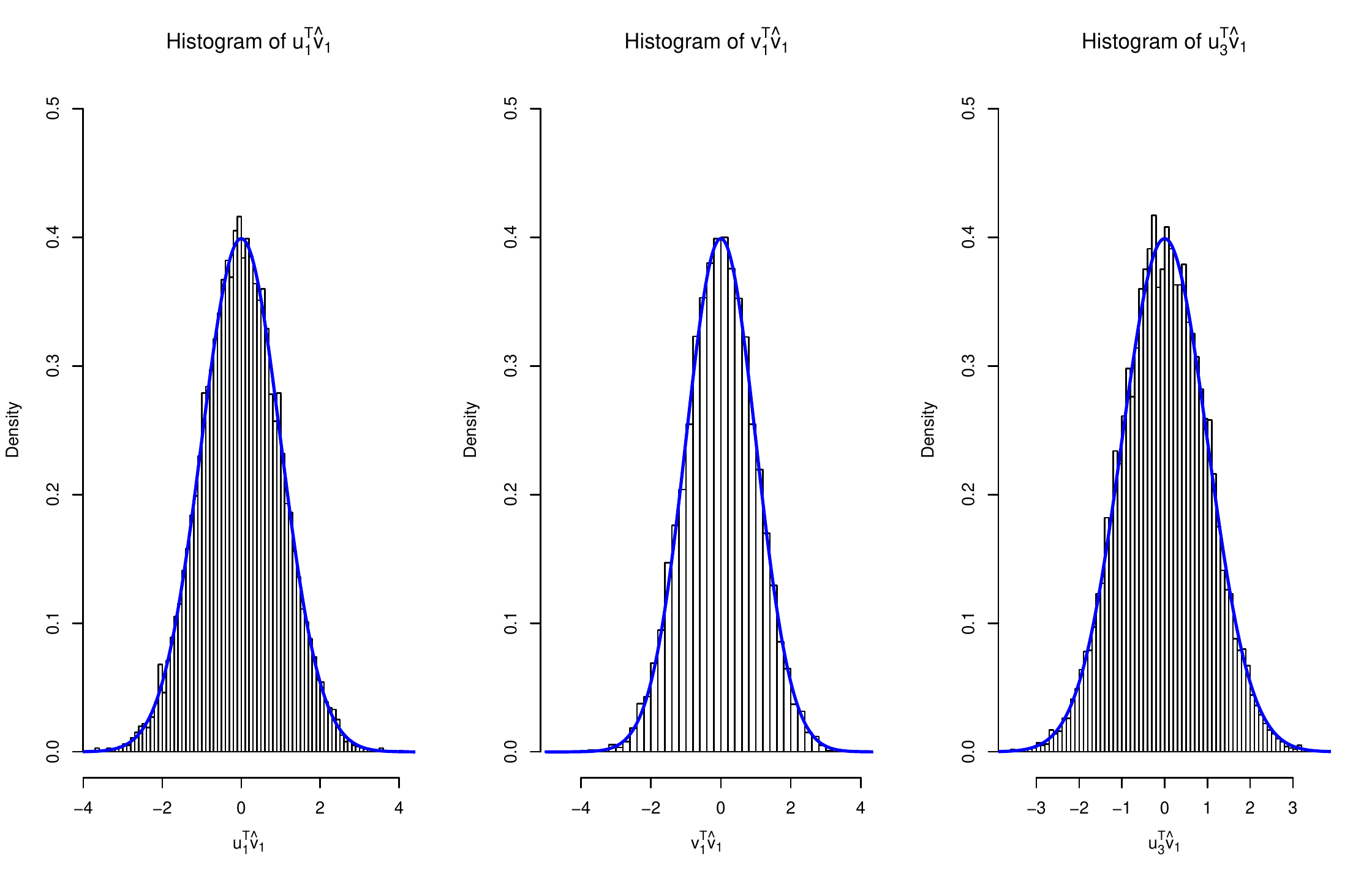}
	\caption{Histograms corresponding to the first eigenvector $\widehat\bbv_1$ when $r = 0.4$, with the blue curves representing the standard normal density. Left panel: $\bu_1^T\widehat\bv_1$; middle panel: $\bv_1^T\widehat\bv_1$; right panel: $\bu_3^T\widehat\bv_1$, where $\bu_1 = \ba_1$ and $\bu_3 = (1, 0, \cdots, 0)^T$.}\label{fig: 1st-eigen-vec}
\end{figure}

\begin{figure}
	\centering
	\includegraphics[width=16.5cm,height=10cm]{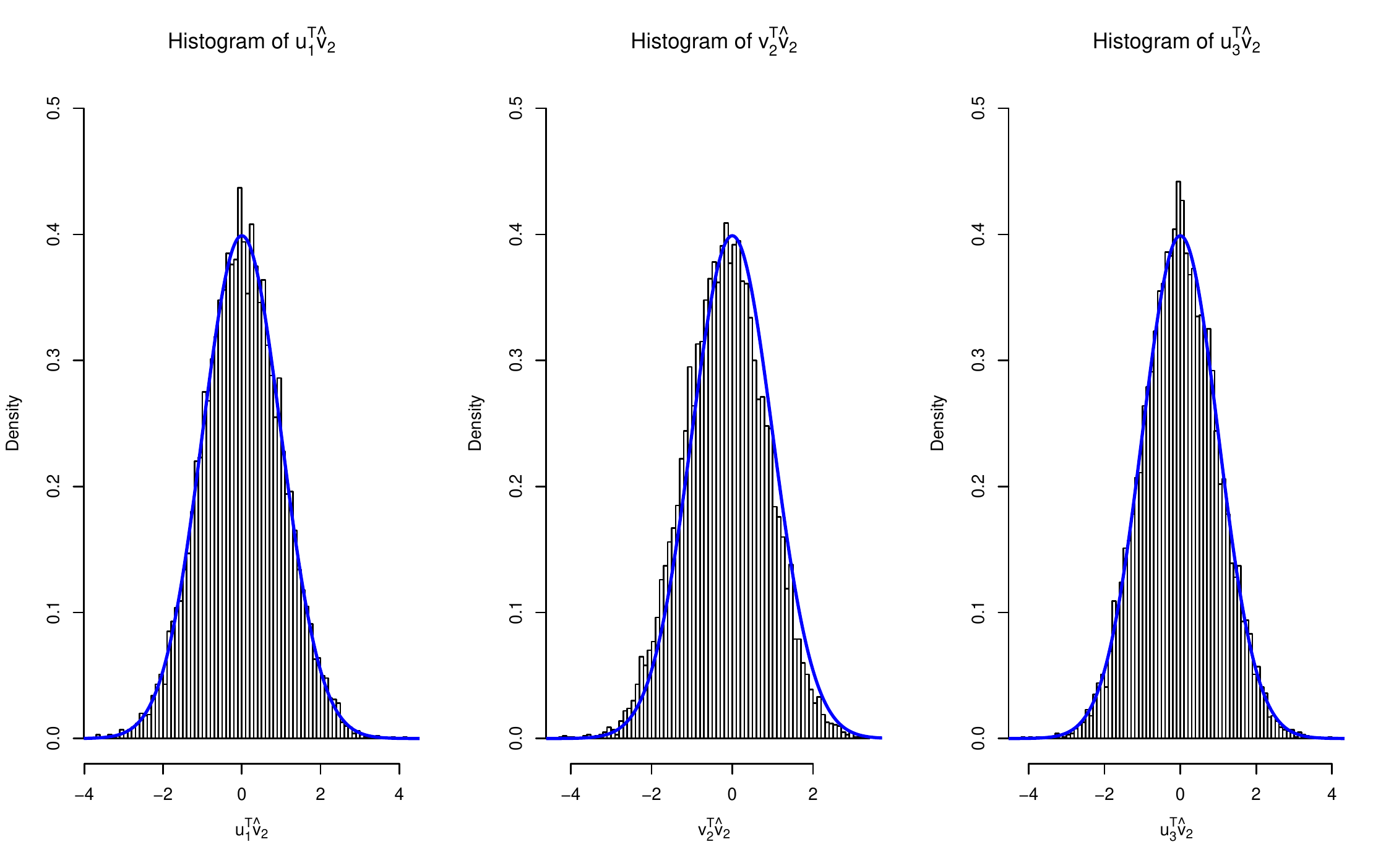}
	\caption{Histograms corresponding to the second eigenvector $\widehat\bbv_2$ when $r = 0.4$, with the blue curves representing the standard normal density. Left panel: $\bu_1^T\widehat\bv_2$; middle panel: $\bv_2^T\widehat\bv_2$; right panel: $\bu_3^T\widehat\bv_1$, where $\bu_1 = \ba_1$ and $\bu_3 = (1, 0, \cdots, 0)^T$.
	}\label{fig: 2nd-eigen-vec}
\end{figure}

For the eigenvectors, we evaluate the asymptotic normality of the linear combination $\bbu^T\widehat\bbv_k$ with $k=1$ and $2$. We experiment with three different values for $\bu$: $\bba_1$, $(1,0,\cdots,0)^T$, and $\bbv_k$.
When $\bu = \bba_1$ or $(1,0, \cdots,0)^T$, we calculate the normalized statistic
$$\frac{t_k\left(\bbu^T\widehat\bbv_k+A_{\sbu,k,t_k}\mathcal{\widetilde P}_{k,t_k}^{1/2}\right)}{\left\{\var\left[(\bbb^T_{\sbu,k,t_k}-\bbu^T\bbv_k\bbv_k^T)\bbW\bbv_k\right]\right\}^{1/2}}$$
using the 10,000 simulated data sets, while when $\bu = \bv_k$ we calculate the normalized statistic  $$\frac{2t^2_k\left(\bbv_k^T\widehat\bbv_k+A_{\sbv_k,k,t_k}\mathcal{\widetilde P}_{k,t_k}^{1/2}\right)}{\left[\var(\bbv_k\bbW^2\bbv_k)\right]^{1/2}}$$ instead. In either of the two cases above, the variance in the denominator is calculated as the sample variance from 2,000 simulated independent copies of the noise matrix $\bbW$. We compare the empirical distributions of the above two normalized statistics  with the standard normal distribution. The simulation results are summarized in Tables \ref{tab2}--\ref{tab8} and Figures \ref{fig: 1st-eigen-vec}--\ref{fig: 2nd-eigen-vec}.

\begin{table}[htbp]
	\centering
	\caption{Simulation results for $(\lambda_1-t_1)/[\var(\bbv_1^T\bbW\bbv_1)]^{1/2}$}
	\begin{tabular}{lrrrrrr}
		\toprule
		$r$    & 0.02 & 0.05 & 0.1 & 0.2 & 0.3 & 0.4 \\
		\midrule
		Mean  & 0.0719  & 0.0149  & -0.0068  & -0.0080  & -0.0024  & 0.0124  \\
		Standard deviation & 1.0107  & 1.0085  & 0.9927  & 1.0115  & 1.0023  & 1.0125  \\
		AD.p-value & 0.0725  & 0.5387  & 0.6263  & 0.2342  & 0.9243  & 0.2010  \\
		\bottomrule
	\end{tabular}%
	\label{tab1}%
\end{table}%

\begin{table}[htbp]
	\centering
	\caption{Simulation results for $(\lambda_2-t_2)/[\var(\bbv_2^T\bbW\bbv_2)]^{1/2}$}	
	\begin{tabular}{lrrrrrr}
		\toprule
		$r$    & 0.02 & 0.05 & 0.1 & 0.2 & 0.3 & 0.4 \\
		\midrule
		Mean  & 1.0761  & 0.2552  & 0.0681  & 0.0272  & 0.0093  & 0.0052  \\
		Standard deviation & 0.9630  & 0.9820  & 0.9872  & 1.0100  & 1.0057  & 1.0005  \\
		AD p-value & 0.5349  & 0.6722  & 0.8406  & 0.1806  & 0.0535  & 0.8341  \\
		\bottomrule
	\end{tabular}%
	
	\label{tab5}%
\end{table}%

Our simulation results in Figure \ref{fig:eig-val} and Tables \ref{tab1}--\ref{tab5} suggest that
the normalized spiked eigenvalues have distributions very close to standard Gaussian which supports our results in Theorem \ref{0619-1}.  Indeed, such a large p-value is extremely impressive given the ``sample size'' (the number of simulations is 10,000).  In general, the simulation results for the eigenvectors support our theoretical findings in Section \ref{Sec3}. However, the results corresponding to the first spiked eigenvector $\widehat\bbv_1$ (Tables \ref{tab2}--\ref{tab4}) are better than those for the second spiked eigenvector $\widehat\bbv_2$ (Tables \ref{tab3}--\ref{tab8}). This is reasonable since for the larger spiked eigenvalue, the negligible terms that we dropped in the proofs of the  asymptotic normality become relatively smaller and thus have smaller finite-sample effects on the asymptotic distributions. For the linear form  $\bu^T\widehat\bv_k$, when $\bu=\bv_k$ the convergence  to standard normal is slower when compared to the case of $\bu\neq \bv_k$. This again supports our theoretical findings in Section \ref{Sec3} and explains why we need to separate the cases of $\bu = \bv_k$ and $\bu\neq\bv_k$. Such effect is especially prominent for $\bv_2^T\widehat\bv_2$, whose sample mean is $-11.8020$ when $r=0.02$ as shown in Table \ref{tab7}.  However, it is seen from the same table (and other tables) that as the spiked eigenvalue increases with $r$, the distribution gets closer and closer to standard Gaussian.	

\begin{table}[htbp]
	\centering
	\caption{Simulation results for $\bbu^T\widehat \bbv_1$ with $\bbu = \ba_1$}	
	\begin{tabular}{lrrrrrr}
		\toprule
		$r$    & 0.02 & 0.05 & 0.1 & 0.2 & 0.3 & 0.4 \\
		\midrule
		Mean  & -0.0573  & -0.0140  & -0.0023  & -0.0045  & -0.0071  & -0.0069  \\
		Standard deviation & 1.0335  & 1.0244  & 1.0011  & 1.0001  & 1.0214  & 1.0016  \\
		AD.p-value & 0.7879  & 0.4012  & 0.2417  & 0.5300  & 0.9482  & 0.9935  \\	
		\bottomrule
	\end{tabular}%
	\label{tab2}%
\end{table}%

\begin{table}[htbp]
	\centering
	\caption{Simulation results for $\bbv_1^T\widehat \bbv_1$}
	\begin{tabular}{lrrrrrr}
		\toprule
		$r$    & 0.02 & 0.05 & 0.1 & 0.2 & 0.3 & 0.4 \\
		\midrule
		Mean  & -1.3288  & -0.4817  & -0.1900  & -0.0742  & -0.0409  & -0.0186  \\
		Standard deviation & 1.0940  & 1.0545  & 1.0338  & 0.9749  & 1.0030  & 1.0005  \\
		AD.p-value & 0.0582  & 0.4251  & 0.0251  & 0.0225  & 0.3312  & 0.2912  \\
		\bottomrule
	\end{tabular}%
	\label{tab3}%
\end{table}%

\begin{table}[htbp]
	\centering
	\caption{Simulation results for $\bbu^T\widehat \bbv_1$ with  $\bbu = (1,0,\cdots, 0)^T$}
	\begin{tabular}{lrrrrrr}
		\toprule
		$r$    & 0.02 & 0.05 & 0.1 & 0.2 & 0.3 & 0.4 \\
		\midrule
		Mean  & 0.0025  & 0.0021  & 0.0003  & 0.0105  & 0.0061  & -0.0122  \\
		Standard deviation & 1.0432  & 1.0354  & 0.9871  & 1.0016  & 1.0205  & 0.9898  \\
		AD.p-value & 0.0044  & 0.4877  & 0.3752  & 0.1514  & 0.1304  & 0.3400  \\
		\bottomrule
	\end{tabular}%
	\label{tab4}%
\end{table}%

\begin{table}[htbp]
	\centering
	\caption{Simulation results for $\bbu^T\widehat \bbv_2$ with $\bbu = \ba_1$}	
	\begin{tabular}{lrrrrrr}
		\toprule
		$r$    & 0.02 & 0.05 & 0.1 & 0.2 & 0.3 & 0.4 \\
		\midrule
		Mean  & 4.2611  & 1.0129  & 0.3067  & 0.0745  & 0.0219  & 0.0037  \\
		Standard deviation & 1.2384  & 1.0952  & 1.0294  & 1.0098  & 1.0280  & 1.0044  \\
		AD p-value & 0.3829  & 0.7535  & 0.3759  & 0.4105  & 0.9129  & 0.9873  \\
		\bottomrule
	\end{tabular}%
	\label{tab6}%
\end{table}%

\begin{table}[htbp]
	\centering
	\caption{Simulation results for $\bbv_2^T\widehat \bbv_2$}
	\begin{tabular}{lrrrrrr}
		\toprule
		$r$    & 0.02 & 0.05 & 0.1 & 0.2 & 0.3 & 0.4 \\
		\midrule
		Mean  & -11.8020  & -4.3274  & -2.0057  & -0.7447  & -0.3526  & -0.1650  \\
		Standard deviation & 1.3775  & 1.1192  & 1.0980  & 1.0343  & 1.0104  & 1.0089  \\
		AD p-value & 0.0000  & 0.0011  & 0.0422  & 0.3964  & 0.4980  & 0.1186  \\
		\bottomrule
	\end{tabular}%
	\label{tab7}%
\end{table}%

\begin{table}[htbp]
	\centering
	\caption{Simulation results for $\bbu^T\widehat \bbv_2$ with  $\bbu = (1,0,\cdots, 0)^T$}
	\begin{tabular}{lrrrrrr}
		\toprule
		$r$    & 0.02 & 0.05 & 0.1 & 0.2 & 0.3 & 0.4 \\
		\midrule
		Mean  & 0.0622  & 0.0204  & 0.0018  & -0.0074  & -0.0119  & -0.0049  \\
		Standard deviation & 1.1221  & 1.0537  & 1.0272  & 1.0022  & 1.0088  & 0.9933  \\
		AD p-value & 0.0003  & 0.5853  & 0.0930  & 0.6011  & 0.2423  & 0.4385  \\
		\bottomrule
	\end{tabular}%
	\label{tab8}%
\end{table}%

\section{A more general asymptotic theory} \label{Sec5}

As mentioned before, the asymptotic theory on the spiked eigenvectors in terms of the general linear combination and on the spiked eigenvalues presented in Section \ref{Sec3} is in fact a consequence of a more general asymptotic theory on the spiked eigenvectors in terms of the bilinear form. In this section, we focus our attention on such a more general asymptotic theory for the bilinear form $\bbx^T\widehat\bbv_k\widehat\bbv_k^T\bby$ with $1 \leq k \leq K$, where $\bx$ and $\by$ are two arbitrary $n$-dimensional unit vectors. See Sections \ref{Sec3.5} and \ref{Sec3.6} for detailed discussions on the technical innovations of our novel ATE theoretical framework and comparisons with the existing literature on the asymptotic distributions of eigenvectors.


For technical reasons, we will break our main results on the asymptotic distributions of the bilinear form  $\bbx^T\widehat\bbv_k\widehat\bbv_k^T\bby$ down to two theorems, where we consider in Theorem \ref{0511-1} the case when either vector $\bx$ or vector $\by$ is sufficiently further away from the population eigenvector $\bv_k$, and then we study in Theorem \ref{0518-1} the case when both vectors $\bx$ and $\by$ are very close to $\bv_k$ 
 The technical treatments for these two cases are different since in the latter scenario, the first order term which determines the asymptotic distribution in Theorem \ref{0511-1} vanishes, and thus we need to consider higher order expansions to obtain the asymptotic distribution in Theorem \ref{0518-1}. Let $\bbJ_{\sbx,\sby,k,t_k}$, $\bbL_{\sbx,\sby,k,t_k}$, and $\bbQ_{\sbx,\sby,k,t_k}$ be the three rank one matrices given in (\ref{0605.2})--(\ref{0605.4}), respectively, in the proof of Theorem \ref{0518-1} in Section \ref{SecA.6}. Denote by $\sigma_k^2=\var[\tr(\bbW\bbJ_{\sbx,\sby,k,t_k})]$ and
\begin{equation}\label{0619.5}
\widetilde\sigma_k^2=\var\left\{\tr\left[\bbW\bbJ_{\sbx,\sby,k,t_k} - \left(\bbW^2-\mathbb{E}\bbW^2\right)\bbL_{\sbx,\sby,k,t_k}\right]+\tr\left(\bbW\bbv_k\bbv_k^T\right)\tr\left(\bbW\bbQ_{\sbx,\sby,k,t_k}\right)\right\}.
\end{equation}
Both of the quantities above play an important role in our more general asymptotic theory.

\begin{thm}\label{0511-1}
Assume that Conditions \ref{as1}--\ref{as3} hold and $\bx$ and $\by$ are two $n$-dimensional unit vectors. Then for each $1 \leq k \leq K$, if
$\sigma_k^2\gg t^{-4}_k\alpha^2_n (|A_{\sbx,k,t_k}|+|A_{\sby,k,t_k}|)^2+t^{-4}_k$ we have the asymptotic expansion
	\begin{eqnarray}\label{0511.2}
	\bbx^T\widehat \bbv_k\widehat \bbv^T_k\bby=a_k+ \tr(\bbW\bbJ_{\sbx,\sby,k,t_k})+O_p \left\{t^{-2}_k\alpha_n (|A_{\sbx,k,t_k}|+|A_{\sby,k,t_k}|)+t^{-2}_k\right\},
	\end{eqnarray}
	where the quantity $a_k=A_{\sbx,k,t_k}A_{\sby,k,t_k}\mathcal{\widetilde P}_{k,t_k}$.
\end{thm}

The assumption of $\sigma_k^2\gg t^{-4}_k\alpha^2_n (|A_{\sbx,k,t_k}|+|A_{\sby,k,t_k}|)^2+t^{-4}_k$ in Theorem \ref{0511-1} requires the variance of random variable $\tr(\bbW\bbJ_{\sbx,\sby,k,t_k})$ not too small, which at high level, requires that either vector $\bbx$ or $\bby$ is sufficiently faraway from the population eigenvector $\bbv_k$. If $\sigma_{ij}\sim 1$ for each $(i,j)$ pair, then such an assumption restricts essentially that  $\|\bbJ_{\sbx,\sby,k,t_k}\|$ should not be too close to zero. This in turn ensures that the first order expansion is sufficient for deriving the asymptotic normality of  $\bbx^T\widehat \bbv_k\widehat \bbv^T_k\bby$. Theorem \ref{0426-1} also entails that a simple upper bound for $\widetilde \sigma_k$ as defined in (\ref{0619.5}) can be shown to be $O(t_k^{-2} \alpha_n)$.

\begin{thm}\label{0518-1}
Assume that Conditions \ref{as1}--\ref{as3} hold and $\bx$ and $\by$ are two $n$-dimensional unit vectors. Then for each $1 \leq k \leq K$, if $\sigma_k^2=O(\widetilde \sigma_k^2)$ and
	$\widetilde \sigma_k^2\gg t^{-6}_k\alpha^4_n (|A_{\sbx,k,t_k}|+|A_{\sby,k,t_k}|)^2+t^{-6}_k$ we have the asymptotic expansion
\begin{align}\label{0518.3h}
\bbx^T\widehat \bbv_k\widehat \bbv^T_k\bby&=a_k+\tr\left[\bbW\bbJ_{\sbx,\sby,k,t_k} - \left(\bbW^2-\mathbb{E}\bbW^2\right)\bbL_{\sbx,\sby,k,t_k}\right]+\tr\left(\bbW\bbv_k\bbv_k^T\right)\tr\left(\bbW\bbQ_{\sbx,\sby,k,t_k}\right)\nonumber\\
&\quad+O_p \left\{|t_k|^{-3}\alpha^2_n (|A_{\sbx,k,t_k}|+|A_{\sby,k,t_k}|)+|t_k|^{-3}\right\},
\end{align}
where the quantity $a_k$ is given in (\ref{0511.2}).
\end{thm}

The ATE theoretical framework for the more general asymptotic theory established in Theorems \ref{0511-1} and \ref{0518-1} is empowered by the following two technical lemmas.

\begin{lem}\label{0426-1}
For any $n$-dimensional unit vectors $\bbx$ and  $\bby$, we have
	\begin{equation}\label{0506.1}
	\bbx^T(\bbW^l-\mathbb{E}\bbW^l)\bby=O_p(\min\{\alpha_n^{l-1},d_{\bx}\alpha_n^l, d_{\by}\alpha_n^l\})
	\end{equation}
	with $l \geq 1$ some bounded positive integer and $d_{\bbx} = \|\bbx\|_\infty$. 
\end{lem}

\begin{lem}\label{1212-1h}
For any $n$-dimensional unit vectors $\bbx$ and  $\bby$, we have  $\mathbb{E}\bbx^T\bbW^l\bby=O(1)$ and
	\begin{equation}\label{1212.3h}
	\mathbb{E}\bbx^T\bbW^l\bby=O(\alpha_n^{l})
	\end{equation}
with $l\geq 2$ some bounded positive integer.
\end{lem}

The detailed proofs of Lemmas \ref{0426-1} and \ref{1212-1h} are provided in Sections \ref{SecB.4} and \ref{SecB.6} of Supplementary Material. Our delicate technical arguments therein establish useful refinements to the classical idea of counting the number of nonzero terms from the random matrix theory. In particular, Lemma \ref{0426-1} is the key building block for high order Taylor expansions that involve polynomials of quantities in the lemma with different choices of $(\bx, \by, l)$.

\section{Discussions} \label{Sec6}

In contrast to the immense literature on the asymptotic distributions for eigenvalues of large spiked random matrices, the counterpart asymptotic theory for eigenvectors has remained largely underdeveloped in statistics literature for years. Yet such a theory is much desired for understanding the precise asymptotic properties of various statistical and machine learning algorithms that build upon the spectral information of the eigenspace constructed from observed data matrix. Our work in this paper provides a first attempt with a general ATE theoretical framework for underpinning the precise asymptotic expansions and asymptotic distributions for spiked eigenvectors and spiked eigenvalues of large spiked random matrices with diverging spikes. Our results complement existing ones in the RMT literature as well as the networks literature.

The family of models in our ATE framework includes many popularly used ones for large-scale applications including network analysis and text analysis such as the stochastic block models with or without overlapping communities and the topic models. Our general asymptotic theory for eigenvectors can be exploited to develop new useful tools for precise statistical inference in these applications. It would be interesting to investigate the problem of reproducible large-scale inference as in  \cite{BarberCandes2015, CandesFanJansonLv2018, LuFanLvNoble2018, FanDemirkayaLiLv2018}; \cite{FanLvSharifvaghefiUematsu2018} in these model settings. It would also be interesting to develop a general method to determine the rank and provide robust rank inference in such high-dimensional low-rank models. These extensions are beyond the scope of the current paper and will be interesting topics for future research.

\appendix
\section{Proofs of main results} \label{SecA}
Recall that Condition \ref{as3} involves two scenarios of the spike strength. We will first prove all the results under scenario i). Then in Section \ref{sec:cond2ii} of Supplementary Material, we will adapt the proofs to show that the same results also hold under scenario ii).
We provide the proofs of Theorems \ref{0619-1}--\ref{0518-1} and Corollary \ref{1109-1} in this appendix. Additional technical details including the proofs of all the lemmas and further discussions on when the asymptotic normality can hold for the asymptotic expansion in Theorem \ref{0518-1} are contained in the Supplementary Material.

\subsection{Proof of Theorem \ref{0619-1}} \label{SecA.1}

The results on the asymptotic distributions of spiked eigenvalues in Theorem \ref{0619-1} are in fact a consequence of those on the asymptotic expansions and asymptotic distributions for the spiked eigenvectors, where a more general asymptotic theory of the eigenvectors is presented in Theorems \ref{0511-1}--\ref{0518-1} in Section \ref{Sec5}. Let us define a matrix-valued function that is referred to as the Green function associated with only the noise part $\bW$
\begin{equation} \label{neweq001}
\bbG(z) = (\bbW-z\bbI)^{-1}
\end{equation}
for $z$ in the complex plane $\mathbb{C}$, where $\bbI$ stands for the identity matrix of size $n$. Recall that $\lambda_1, \cdots , \lambda_n$ are the eigenvalues of matrix $\bX$ and $\widehat\bbv_1, \cdots, \widehat{\bbv}_n$ are the corresponding eigenvectors. By Weyl's inequality, it holds that $\max|\lambda_i-d_i|\le \|\bW\|$. Thus, in view of Condition \ref{as3}  and Lemma \ref{0505-1} in Supplementary Material, all the spiked eigenvalues $\lambda_k$ with $1 \leq k \leq K$ of the observed random matrix $\bX$ have magnitudes of larger order than the eigenvalues of the noise matrix $\bbW$ with significant probability as the matrix size $n$ increases. This entails that with significant probability, matrices $\bbG(\lambda_k)$ with $1 \leq k \leq K$ are well defined and nonsingular. For the rest of this proof, we restrict all the derivations on such an event that holds with asymptotic probability one.

It follows from the definition of the eigenvalue, the representation $\bX = \bH + \bW = \bV \bD \bV^T + \bW$, (\ref{neweq001}), and the properties of the determinant function $\det(\cdot)$ that for each $1 \leq k \leq K$,
\begin{align*}
0 &= \det(\bX - \lambda_k \bbI) =  \det(\bbW-\lambda_k\bbI+\bbV\bbD\bbV^T)=\det[\bbG^{-1}(\lambda_k)+\bbV\bbD\bbV^T]\non
&=\det[\bbG^{-1}(\lambda_k)] \det[\bbI+ \bbG(\lambda_k)\bbV\bbD\bbV^T],
\end{align*}
which leads to $\det[\bbI+ \bbG(\lambda_k)\bbV\bbD\bbV^T] = 0$ since $\det[\bbG^{-1}(\lambda_k)] = \det[\bbG(\lambda_k)]^{-1}$ is nonzero. Using the identity $\det(\bbI+\bbA\bbB)=\det(\bbI+\bbB\bbA)$ for matrices $\bA$ and $\bB$, we obtain for each $1 \leq k \leq K$,
\begin{equation} \label{0523.1}
0 = \det[\bbI+ \bbG(\lambda_k)\bbV\bbD\bbV^T] = \det[\bbI+ \bbD\bbV^T\bbG(\lambda_k)\bbV],
\end{equation}
where the second $\bbI$ represents an identity matrix of size $K$ and we slightly abuse the notation for simplicity. Since the diagonal matrix $\bD$ is nonsingular by assumption, it follows from (\ref{0523.1}) that
\begin{eqnarray}\label{0314.5h}
\det[d_k\bbV^T\bbG(\lambda_k)\bbV+d_k\bbD^{-1}] = d_k \det(\bbD^{-1}) \det[\bbI+ \bbD\bbV^T\bbG(\lambda_k)\bbV] = 0
\end{eqnarray}
for each $1 \leq k \leq K$.

By the asymptotic expansions in (\ref{0508.5}), Lemmas \ref{0426-1} and \ref{1212-1h}, and Weyl's inequality $\max|\lambda_k-d_k|\leq \|\bW\|$, we have for $j\neq \ell$,
$
d_k\bv_j^T\bbG(\lambda_k)\bv_{\ell} = -d_kO_p(\lambda_k^{-2})= O_p(1/|d_k|).
$
Thus, we can see that all off diagonal entries of matrix $d_k\bbV^T\bbG(\lambda_k)\bbV+d_k\bbD^{-1}$ in (\ref{0314.5h}) are of order $O_p(1/|d_k|)$.   
For $j\neq k$, the $j$th diagonal entry of $d_k\bbV^T\bbG(\lambda_k)\bbV+d_k\bbD^{-1}$ equals $d_k\bbv_j^T\bbG(\lambda_k)\bbv_j+d_k/d_j$. By (\ref{0508.4}) and Lemma \ref{0426-1}, we have $d_k\bbv_j^T\bbG(\lambda_k)\bbv_j+1=o_p(1)$. Moreover, by  Condition \ref{as3} $|d_k/d_j-1|\ge c$ for some positive constant $c$.  Hence, all these diagonal entries but the $k$th one are of order at least $O_p(1)$.  Thus the matrix $(d_k\bbv_i^T\bbG(\lambda_k)\bbv_j+\delta_{ij} d_k/d_i)_{1\le i,j\le K,\, i,j\neq k}$ is invertible with significant probability, where $\delta_{ij} = 1$ when $i = j$ and $0$ otherwise. Recall the determinant identity for block matrices from linear algebra
$$\det\left(
\begin{array}{ccc}
\bA_{11} &\bA_{12}\\
\bA_{21} & \bA_{22}\\
\end{array}
\right)=\det(\bA_{22})\det(\bA_{11}-\bA_{12}\bA_{22}^{-1}\bA_{21})$$
when the lower right block matrix $\bA_{22}$ is nonsingular. Treating the $k$th diagonal entry of $d_k\bbV^T\bbG(\lambda_k)\bbV+d_k\bbD^{-1}$ as the first block,  we have with significant probability
\begin{eqnarray} \label{0523.5}
\det[d_k\bbV^T\bbG(\lambda_k)\bbV+d_k\bbD^{-1}]=0
\end{eqnarray}
entailing $d_k\bbv_k^T\bbG(\lambda_k)\bbv_k+1=d_k\bbv_k^T\bbF_k(\lambda_k)\bbv_k$, where $\bbF_k(z)=\bbG(z)\bbV_{-k}[\bbD_{-k}^{-1}+\bbV_{-k}^T\bbG(z)\bbV_{-k}]^{-1}\\ \cdot \bbV_{-k}^T \bbG(z)$ and $\bA_{- k}$ denotes the submatrix of matrix $\bA$ by removing the $k$th column. In light of (\ref{0523.5}) and the solution $\widehat t_k$ to equation (\ref{0523.6}) in the proof of Theorem \ref{0511-1} in Section \ref{SecA.5}, it holds from the uniqueness of $\widehat t_k$ that
\begin{equation} \label{neweq002}
\lambda_k = \widehat t_k.
\end{equation}
Therefore, combining equality (\ref{neweq002}) with asymptotic expansion
$\widehat t_k-t_k=\bbv_k^T\bbW\bbv_k+O_p(\alpha_n/t_k)$ obtained in (\ref{0516.15}) completes the proof of Theorem \ref{0619-1}.

\subsection{Proof of Theorem \ref{0609-1}} \label{SecA.2}

The results on the asymptotic distributions of spiked eigenvectors in Theorem \ref{0609-1} are also an implication of a more general asymptotic theory of the eigenvectors presented in Theorems \ref{0511-1}--\ref{0518-1} in Section \ref{Sec5} on the delicate asymptotic expansions and asymptotic distributions for the spiked eigenvectors. Recall that $\widehat\bbV=(\widehat \bbv_1, \cdots,\widehat \bbv_K)$ with $\widehat\bbv_k$ for $1 \leq k \leq K$ the empirical spiked eigenvectors of the observed random matrix $\bX$. Without loss of generality, let us choose the direction of eigenvector $\widehat \bbv_k$ such that $\widehat\bbv_k^T\bbv_k\ge 0$. Clearly, fixing the direction of $\widehat\bbv_k$ does not affect the distribution of $\bbx^T\widehat\bbv_k\widehat\bbv_k^T\bby$; that is, its distribution stays the same when $-\widehat\bbv_k$ is chosen as the eigenvector. We will separately consider the two cases of $\bbv_k^T\widehat{\bbv}_k$ and $\bbu^T\widehat{\bbv}_k$ with $\bbu\neq \bbv_k$, where the former relies on the second order expansion given in (\ref{0517.3}) in the proof of Theorem \ref{0518-1} in Section \ref{SecA.6}, and the latter utilizes the first order expansion given in \eqref{0516.1a} in the proof of Theorem \ref{0511-1} in Section \ref{SecA.5}.

We first consider $\bbv_k^T\widehat{\bbv}_k$. Choosing $\bx=\by=\bv_k$ in Theorem \ref{0511-1} gives $a_k=A^2_{\sbv_k,k,t_k} \mathcal{\widetilde P}_{k,t_k}$. 
  By Lemma \ref{1212-1h}, it holds that
\begin{equation}\label{1016.1h}
\mathcal{P}(\bbv_k,\bbv_k,t_k)=-\sum_{l=0,l\neq 2}^L\frac{\bbv_k^T\mathbb{E}\bbW^l\bbv_k}{t_k^l}=-1+O(\alpha_n^2/t_k^2)
\end{equation}
and
\begin{equation}\label{1016.2h}
\|\mathcal{P}(\bbv_k,\bbV_{-k},t_k)\|=\|-\sum_{l=0,l\neq 2}^L\frac{\bbv_k^T\mathbb{E}\bbW^l\bbV_{-k}}{t_k^l}\|=O(\alpha_n^2/t_k^2).
\end{equation}
Moreover, recalling the definition of $A_{\bu,k,t_k}$ in (\ref{0619.2}), $A_{\bu,k,t_k}$ can be rewritten as
$$A_{\bv_k,k,t_k}=\mathcal{P}(\bbv_k,\bbv_k,t)-t_k^{-1}\mathcal{P}(\bbv_k,\bbV_{-k},t_k)\Big(\bbD_{-k}^{-1}+\mathcal{R}(\bbV_{-k},\bbV_{-k},t_k)\Big)^{-1}\mathcal{P}(\bbV_{-k},\bbv_k,t_k).$$
Therefore, by (\ref{1016.1h})--(\ref{1016.2h}), (\ref{1010.1h}), and (\ref{1011.5h}), we have
\begin{equation}\label{0927.1h}
A_{\bv_k,k,t_k}=-1+O(\alpha_n^2/t_k^2) \ \text{ and } \  \mathcal{\widetilde P}_{k,t_k}=1+O(\alpha_n^2/t_k^2).
\end{equation}

Now recall the second order expansion of $\bbx^T\widehat{\bbv}_k\widehat{\bbv}_k\bby$ given in (\ref{0517.3}) in the proof of Theorem \ref{0518-1}. We next calculate the orders of each term in the expansion (\ref{0517.3}). First, we consider $\bbb^T_{\sbv_k,k,t_k}$. By (\ref{1016.2h}), (\ref{1010.1h}), and the definition (\ref{0619.1}),  we have
	\begin{equation}\label{1017.2h}
	\|\bbb^T_{\sbv_k,k,t_k}-\bbv^T_k\|=\left\|\mathcal{R}(\bbv_k,\bbV_{-k},t_k)\Big((\bbD_{-k})^{-1}+\mathcal{R}(\bbV_{-k},\bbV_{-k},t_k)\Big)^{-1}\bbV_{-k}^T\right\|=O(\alpha_n^2/t_k^2).
	\end{equation}
This together with \eqref{0927.1h}  entails that
	\begin{equation}\label{1017.3h}
	\|\bbb^T_{\sbv_k,k,t_k}+A_{\sbv_k,k,t_k}\mathcal{\widetilde P}_{k,t_k}\bbv_k^T\|=O(\alpha_n^2/t_k^2).
	\end{equation}
  It follows from Lemma \ref{0426-1} and \eqref{0927.1h} that
	$$A_{\sbv_k,k,t_k} \mathcal{\widetilde P}_{k,t_k}(\bbb^T_{\sbv_k,k,t_k}+A_{\sbv_k,k,t_k}\mathcal{\widetilde P}_{k,t_k}\bbv_k^T)\bbW\bbv_k / t_k=O_p(\alpha^2_n/|t_k|^3), $$
		\begin{align*} 
	&\mathcal{\widetilde P}_{k,t_k}t^{-2}_k\left[2\mathcal{\widetilde P}_{k,t_k}\left(A_{\sbv_k,k,t_k}\bbb^T_{\sbv_k,k,t_k}+A_{\sbv_k,k,t_k}\bbb^T_{\sbv_k,k,t_k}\right)\bbW\bbv_k\bbv_k^T+\bbb^T_{\sbv_k,k,t_k}\bbW\bbv_k\bbb^T_{\sbv_k,k,t_k}\right]\bbW\bbv_k\nonumber\\
	&\quad+2A_{\sbv_k,k,t_k}A_{\sby,k,t_k}t^{-2}_k\left(\bbv_k^T\bbW\bbv_k\right)^2\nonumber \\
	&\quad+A_{\sbv_k,k,t_k}\mathcal{\widetilde P}_{k,t_k}\Big\{t^{-2}_k\bbx^T\bbW\bbv_k\bbv^T_k\bbW\bbv_k-t^{-2}_k\bbv^T_k\bbW\bbv_k\mathcal{R}(\bbv_k,\bbV_{-k},t)\nonumber\\
	&\quad\times\left[\bbD_{-k}^{-1}+\mathcal{R}(\bbV_{-k},\bbV_{-k},t_k)\right]^{-1}\bbV_{-k}^T\bbW\bbv_k\Big\}\nonumber \\
	&\quad+A_{\sbv_k,k,t_k}\mathcal{\widetilde P}_{k,t_k}\Big\{t^{-2}_k\bbv_k^T\bbW\bbv_k\bbv^T_k\bbW\bbv_k-t^{-2}_k\bbv^T_k\bbW\bbv_k\mathcal{R}(\bbv_k,\bbV_{-k},t)\nonumber\\
	&\quad \times\left[\bbD_{-k}^{-1}+\mathcal{R}(\bbV_{-k},\bbV_{-k},t_k)\right]^{-1}\bbV_{-k}^T\bbW\bbv_k\Big\}\nonumber \\
	&\quad+A_{\sbv_k,k,t_k}\mathcal{\widetilde P}_{k,t_k}t^{-2}_k\mathcal{R}(\bbv_k,\bbV_{-k},t_k)\left[\bbD_{-k}^{-1}+\mathcal{R}(\bbV_{-k},\bbV_{-k},t_k)\right]^{-1}\bbV_{-k}^T(\bbW^2-\mathbb{E}\bbW^2)\bbv_k\nonumber \\
	&\quad+A_{\sbv_k,k,t_k}\mathcal{\widetilde P}_{k,t_k}t^{-2}_k\mathcal{R}(\bbv_k,\bbV_{-k},t_k)\left[\bbD_{-k}^{-1}+\mathcal{R}(\bbV_{-k},\bbV_{-k},t_k)\right]^{-1}\bbV_{-k}^T(\bbW^2-\mathbb{E}\bbW^2)\bbv_k \nonumber\\
	&=O_p(\frac{1}{|t_k|^2}+\frac{\alpha_n}{|t_k|^3}),
	\end{align*}
	and
	\begin{align*}
	\mathcal{\widetilde P}_{k,t_k}&t^{-2}_k(A_{\sbv_k,k,t_k}\bbx^T+A_{\sbv_k,k,t_k}\bby^T)(\bbW^2-\mathbb{E}\bbW^2)\bbv_k
	+3t^{-2}_kA_{\sbv_k,k,t_k}A_{\sbv_k,k,t_k}\mathcal{\widetilde P}_{k,t_k}\bbv_k^T(\bbW^2-\mathbb{E}\bbW^2)\bbv_k\\
	&=\frac{\bbv_k^T(\bbW^2-\mathbb{E}\bbW^2)\bbv_k }{t^2_k}+O_p(\alpha_n^3/t_k^4).
	\end{align*}
Substituting the above equations into (\ref{0517.3}) results in
\begin{equation} \label{neweq005}
\bbv^T_k\widehat\bbv_k\widehat\bbv_k^T\bbv_k-A^2_{\sbv_k,k,t_k}\mathcal{\widetilde P}_{k,t_k}=-\bbv_k^T(\bbW^2-\mathbb{E}\bbW^2)\bbv_k / t^2_k+O_p(|t_k|^{-2}+\alpha^2_n/|t_k|^3),
\end{equation}
where the leading term of the asymptotic expansion now depends on the second moments of the noise matrix $\bW$. Recall that $\bv_k^T\widehat\bv_k\ge 0$. By (\ref{0927.1h}) and (\ref{neweq005}) we have
\begin{align} \label{neweq007}
\bbv_k^T\widehat\bbv_k+A_{\sbv_k,k,t_k} \mathcal{\widetilde P}_{k,t_k}^{1/2}&=-\frac{\bbv_k^T(\bbW^2-\mathbb{E}\bbW^2)\bbv_k}{2A_{\sbv_k,k,t_k}\mathcal{\widetilde P}_{k,t_k}^{1/2} t^2_k}+O_p(\alpha^2_n/t_k^3) \nonumber \\
& =-\frac{\bbv_k^T(\bbW^2-\mathbb{E}\bbW^2)\bbv_k}{2 t^2_k}+O_p(|t_k|^{-2}+\alpha^2_n/|t_k|^3).
\end{align}

We now consider an arbitrary unit vector $\bbu \in \mathbb{R}^n$ with $|\bu^T\bv_k| \in [0,1) $ for investigating the asymptotic distributions of the general linear combinations $\bbu^T\widehat\bbv_k$. Recall the first order expansion given in \eqref{0516.1a} in the proof of Theorem \ref{0511-1} and \eqref{1017.3h} that
	\begin{align} \label{neweq004}
		\bbu^T\widehat\bbv_k\widehat\bbv_k^T\bbv_k-A_{\sbu,k,t_k}A_{\sbv_k,k,t_k}\mathcal{\widetilde P}_{k,t_k} & =-A_{\sbv_k,k,t_k} \mathcal{\widetilde P}_{k,t_k}(\bbb^T_{\sbu,k,t_k}+A_{\sbu,k,t_k}\mathcal{\widetilde P}_{k,t_k}\bbv_k^T)\bbW\bbv_k / t_k \nonumber \\
		& \quad +O_p(\alpha_n/t^2_k).
\end{align}
Then dividing (\ref{neweq004}) by $\bv_k^T\widehat\bv_k$ and using (\ref{0927.1h}) and (\ref{neweq007}), we can deduce that
\begin{align} \label{neweq006}
\bbu^T\widehat\bbv_k+A_{\sbu,k,t_k} \mathcal{\widetilde P}_{k,t_k}^{1/2} & =\mathcal{\widetilde P}_{k,t_k}^{1/2} (\bbb^T_{\sbu,k,t_k}+A_{\sbu,k,t_k}\mathcal{\widetilde P}_{k,t_k}\bbv_k^T)\bbW\bbv_k / t_k+O_p(\alpha_n / t^2_k) \nonumber \\
& = (\bbb^T_{\sbu,k,t_k}-\bbu^T\bbv_k\bbv_k^T)\bbW\bbv_k/t_k+O_p(\alpha_n/t^2_k).
\end{align}

In view of the asymptotic expansions in (\ref{neweq006}) and (\ref{neweq007}), we can see that the desired asymptotic normalities in the two parts of Theorem \ref{0609-1} follow from the conditions of Lemmas \ref{0524-1} or \ref{0524-1h}. More specifically, for (\ref{neweq006}) if $\alpha_n^{-2}d_k^2\var[(\bbb^T_{\sbu,k,t_k}-\bbu^T\bbv_k\bbv_k^T)\bbW\bbv_k]\rightarrow \infty$, then we have $\alpha_n/t^2_k\ll \{\var[(\bbb^T_{\sbu,k,t_k}-\bbu^T\bbv_k\bbv_k^T)\bbW\bbv_k]\}^{1/2}$ and thus the first part of Theorem \ref{0609-1} in (\ref{0609.1}) holds in view of (\ref{neweq006}). Furthermore, if $(\bbb^T_{\sbu,k,t_k}-\bbu^T\bbv_k\bbv_k^T)\bbW\bbv_k$ is $\bW^1$-CLT, then $(\bbb_{\sbu,k,t_k}-\bbv_k\bbv_k^T\bbu, \bbv_k)$ is also $\bW^1$-CLT and thus we have  $$\frac{t_k\left(\bbu^T\widehat\bbv_k+A_{\sbu,k,t_k}\mathcal{\widetilde P}_{k,t_k}^{1/2}\right)-\mathbb{E}\left[(\bbb^T_{\sbu,k,t_k}-\bbu^T\bbv_k\bbv_k^T)\bbW\bbv_k\right]}{\left\{\var\left[(\bbb^T_{\sbu,k,t_k}-\bbu^T\bbv_k\bbv_k^T)\bbW\bbv_k\right]\right\}^{1/2}}
\toD N(0, 1). $$
Similarly, the second part of Theorem \ref{0609-1} in (\ref{0609.2}) also holds under the condition ($\alpha_n^{-4}d_k^2+1)\var[\bbv_k^T(\bbW^2-\mathbb{E}\bbW^2)\bbv_k]\rightarrow \infty$ and the CLT holds if $(\bv_k, \bbv_k)$ is $\bW^2$-CLT. This concludes the proof of Theorem \ref{0609-1}.

\subsection{Proof of Theorem \ref{1107-1}} \label{SecA.3}
The results on the asymptotic distributions of spiked eigenvalues and spiked eigenvectors in Theorem \ref{1107-1} are an application of those in Theorems \ref{0619-1} and \ref{0609-1} for a more specific structure of the low rank model \eqref{model}, including the stochastic block model with both non-overlapping and overlapping communities as special cases. 

First, note that (\ref{eq002}) implies that the condition of Lemma \ref{0524-1} holds for $\bbv^T_k\bbW\bbv_k$ under   Condition \ref{as7}.  Consequently, $(\bbv_k, \bbv_k)$ is $\bbW^1$-CLT. In addition, \eqref{eq002} ensures that $\mathbb{E}(\bbv_k^T\bbW\bbv_k-\mathbb{E}\bbv_k^T\bbW\bbv_k)^2\gg \alpha^2_n/d^2_k$ under Condition \ref{as7}. Therefore, it follows from  Theorem \ref{0619-1} that the first result of Theorem \ref{1107-1} holds. Recall that in  (\ref{0614.3}), $s_{\sbx,\sby}$ is defined as the expected value of the conditional variance of $\bbv_k^T(\bbW^2-\mathbb{E}\bbW^2)\bbv_k$. By definition, we have $\var[\bbv_k^T(\bbW^2-\mathbb{E}\bbW^2)\bbv_k]\geq s_{\sbx,\sby}\ge c\sigma_{\min}^2n$. Thus the condition $(\alpha_n^{-4}d_k^2+1)\var[\bbv_k^T(\bbW^2-\mathbb{E}\bbW^2)\bbv_k]\rightarrow \infty$ in Theorem \ref{0609-1} is ensured by the assumptions
$$\sigma_{\min}^2n\rightarrow    \infty, \  \frac{|d_K|\sigma_{\min}}{\alpha_n}\rightarrow \infty, \  \alpha_n\le n^{1/2}$$
in Condition \ref{as7}. Moreover, by (\ref{neweq025}) we can see that the conditions of Lemma \ref{0524-1h} are satisfied for $\bbv_k^T(\bbW^2-\mathbb{E}\bbW^2)\bbv_k$ under Condition \ref{as7}. Thus $(\bbv_k, \bbv_k)$ is $\bbW^2$-CLT. Therefore, (\ref{0927.4h}) holds by an application of \eqref{0609.2} in Theorem \ref{0609-1}.

It remains to show that the condition
\begin{equation} \label{neweq008}
\var[(\bbb^T_{\sbu,k,t_k}-\bbu^T\bbv_k\bbv_k^T)\bbW\bbv_k]\gg \alpha_n^2/d_k^2
\end{equation}
in Theorem \ref{0609-1} can be guaranteed by Condition \ref{as7}. Then the expansion in \eqref{0609.1} holds.  Moreover, the condition $\sigma_{\min}^{-1}\big\|\bbv_k[\bbb^T_{\sbu,k,t_k}-\bbu^T\bbv_k\bbv_k^T]\big\|_{\infty}\rightarrow 0$ ensures that $(\bbb_{\sbu,k,t_k}-\bbv_k\bbv_k^T\bbu, \bbv_k)$ is $\bbW^1$-CLT. Combining these results entails that the asymptotic normality \eqref{0609.1h} holds.  Now we proceed to verify \eqref{neweq008}.
Consider an arbitrary unit vector $\bu \in \mathbb{R}^n$ satisfying $|\bbu^T\bbv_k|\in [0,1-\ep]$ for some positive constant $\ep$. Recalling the definition of $\bbb_{\sbu,k,t}$ in (\ref{0619.3}), we have $\bbb^T_{\sbu,k,t}\bbv_k=\{\bbu^T-\mathcal{R}(\bbu,\bbV_{-k},t)[(\bbD_{-k})^{-1}+\mathcal{R}(\bbV_{-k},\bbV_{-k},t)]^{-1}\bbV_{-k}^T\}\bbv_k=\bbu^T\bbv_k$. Thus it holds that $\bbb^T_{\sbu,k,t_k}-\bbu^T\bbv_k\bbv_k^T = \bbb^T_{\sbu,k,t_k} - \bbb^T_{\sbu,k,t_k}\bbv_k\bbv_k^T =  \bbb^T_{\sbu,k,t_k}(\bI-\bv_k\bv_k^T)$. Moreover, similar to \eqref{eq002} we can show that
 \begin{align}\label{0927.5h} \left[{\mathbb{E}(\bbu^T\bbW\bbv_k-\mathbb{E}\bbu^T\bbW\bbv_k)^2}\right]^{1/2} 
	\geq \sigma_{\min}(2-2\|\bv_k\|_\infty^2)^{1/2}.
	\end{align}
	 This ensures that there exists some positive constant $c_1$ such that
\begin{eqnarray} \label{0608.1}
&&\var[(\bbb^T_{\sbu,k,t_k}-\bbu^T\bbv_k\bbv_k^T)\bbW\bbv_k]^2
\ge \sigma_{\min}^2 (2-2\|\bv_k\|_\infty^2)\|\bbb^T_{\sbu,k,t_k}-\bbu^T\bbv_k\bbv_k^T\|^2\non
&& \ge \sigma_{\min}^2 (2-2\|\bv_k\|_\infty^2)\|\bbb^T_{\sbu,k,t_k}(\bI-\bbv_k\bbv_k^T)\|^2  \non
&&\ge c_1\sigma_{\min}^2[-(\bbu^T\bbv_k)^2+\bbb^T_{\sbu,k,t_k}\bbb_{\sbu,k,t_k}],
\end{eqnarray}
where we have applied $\bbb^T_{\sbu,k,t}\bbv_k= \bbu^T\bbv$ again in the last step.

Let $\widetilde \bbV=(\bbv_{K+1}, \cdots,\bbv_n)$ be an $n\times(n-K)$ matrix such that $(\bbV,\widetilde \bbV)$ is an orthogonal matrix of size $n$. Then the $n$-dimensional unit vector $\bu$ can be represented as
$\bbu=\sum_{i=1}^na_i\bbv_i$ for some scalars $a_i$'s. For each $1 \leq k \leq K$, by the definition of $\mathcal R$ in (\ref{0619.0}) and  Lemma \ref{1212-1h} we can show that
\begin{equation}\label{1017.1h}
\|\mathcal{R}(\bbV_{-k},\bbV_{-k},t_k)+t_k^{-1}\bbI\|=O\Big(\frac{\alpha_n^2}{|t_k|^3}\Big) \ \text{ and }  \  \|\mathcal{R}(\bbu,\bbV_{-k},t_k)+t_k^{-1}\bbu^T\bbV_{-k}\|=O\Big(\frac{\alpha_n^2}{|t_k|^3}\Big).
\end{equation}
Therefore it holds that
\begin{align} \label{neweq009}
\Big\| & \mathcal{R}(\bbu,\bbV_{-k},t_k)  [\bbD_{-k}^{-1}+\mathcal{R}(\bbV_{-k},\bbV_{-k},t_k)]^{-1}\bbV_{-k}^T + \sum_{1 \leq i\neq k \leq K}a_i(t_kd_i^{-1}-1)^{-1}\bbv_i^T\Big\| \nonumber\\
& \quad =O(\alpha^2_n/t^2_k).
\end{align}
Then it follows from (\ref{neweq009}) and (\ref{0619.3}) that
$$\bbb_{\sbu,k,t_k}^T=\sum_{i=1}^na_i\bbv_i^T-\mathcal{R}(\bbu,\bbV_{-k},t_k)  [\bbD_{-k}^{-1}+\mathcal{R}(\bbV_{-k},\bbV_{-k},t_k)]^{-1}\bbV_{-k}^T$$
and
\begin{eqnarray} \label{0528.1}
\left\|\bbb_{\sbu,k,t_k}- a_k\bbv_k-\sum_{1 \leq i\neq k \leq K}a_i[1+(t_kd_i^{-1}-1)^{-1}]\bbv_i-\sum_{i=K+1}^na_i\bbv_i\right\|=O(\alpha^2_n / t^2_k).
\end{eqnarray}

We denote by $\bbc_k = a_k\bbv_k+\sum_{1 \leq i\neq k \leq K}a_i[1+(t_kd_i^{-1}-1)^{-1}]\bbv_i+\sum_{i=K+1}^na_i\bbv_i$. By (\ref{0528.1}), we can obtain
\begin{align} \label{1107.5}
- & (\bbu^T\bbv_k)^2+\bbb^T_{\sbu,k,t_k}\bbb_{\sbu,k,t_k} = -a_k^2+\|\bbc_k\|^2+\left\|\bbb_{\sbu,k,t_k}- \bbc_k\right\|^2+2(\bbb_{\sbu,k,t_k}- \bbc_k)^T\bbc_k\non
& =\sum_{1 \leq i \neq k \leq K}a^2_i [1+(t_kd_i^{-1}-1)^{-1}]^2+\sum_{i=K+1}^n a_i^2+O(\alpha^2_n / t^2_k) \nonumber \\
& \quad {+ \text{ some small order term}},
\end{align}
where the small order term takes a rather complicated form and thus we omit its expression for simplicity. Since by assumption $|\bbu^T\bbv_k|\in [0,1-\ep]$, $\bbu=\sum_{i=1}^na_i\bbv_i$ is a unit vector, and $(\bv_1, \cdots, \bv_n)$ is an orthogonal matrix, it holds that
\begin{equation} \label{neweq010}
\sum_{1 \leq i\neq k \leq n}^n a_i^2\ge 1-(1-\ep)^2.
\end{equation}
Moreover, Condition \ref{as7} and Lemma \ref{lem: define-t} together entail that $|t_kd_i^{-1}|$ is bounded away from 0 and 1. Thus there exists some positive constant $c_2 < 1$ such that
\begin{equation} \label{neweq011}
[1+(t_kd_i^{-1}-1)^{-1}]^2\ge c_2
\end{equation}
for each $1 \leq i \neq k \leq K$.   Therefore, combining (\ref{0608.1}) and (\ref{1107.5})--(\ref{neweq011}), and by the assumption $\sigma_{\min} \gg \alpha_n/t_k$, we can obtain the desired claim in (\ref{neweq008}), which completes the proof of Theorem \ref{1107-1}.

\subsection{Proof of Corollary \ref{1109-1}} \label{SecA.4}
The conclusions of Corollary \ref{1109-1} follow directly from the results of Theorem \ref{1107-1}.

\subsection{Proof of Theorem \ref{0511-1}} \label{SecA.5}
The more general asymptotic theory in Theorem \ref{0511-1} focuses on the first order asymptotic expansion for the bilinear form $\bbx^T\widehat\bbv_k\widehat\bbv_k^T\bby$ with $\bx$ and $\by$ two arbitrary $n$-dimensional unit vectors, while that in Theorem  \ref{0518-1} further establishes the higher order (which is second order) asymptotic expansion for the same bilinear form. We begin with the analysis for the first order asymptotic expansion. The main ingredients of the proof are as follows. First, we represent $\bbx^T\widehat\bbv_k\widehat\bbv_k^T\bby$ as an integral which is a functional of $\bbX = \bbH+\bbW$. By doing so we can deal with the matrix  $\bbH+\bbW$ instead of the eigenvectors. Second, for the functional of $\bbH+\bbW$ obtained in the previous step we extract the $\bbH$ part from $\bbH+\bbW$ and further obtain a functional of $\bbW$. Roughly speaking, we can get an explicit function of form  $f((\bbW-t\bbI)^{-1})$ with $|t|\gg \|\bbW\|$. Third, by the matrix series expansion $(\bbW-t\bbI)^{-1}=-\sum_{l=0}^{\infty}t^{-(l+1)}\bbW^l$, the function $f((\bbW-t\bbI)^{-1})$ can be approximated by $f(-\sum_{l=0}^{L}t^{-(l+1)}\bbW^l)$ for some positive integer $L$.  Fourth, we can then calculate the first (second or higher) order expansion of $f(-\sum_{l=0}^{L}t^{-(l+1)}\bbW^l)$ since we have an explicit expression of function $f$.

To facilitate our technical derivations, let us recall some basic matrix identities from the Sherman--Morrison--Woodbury formula. For any matrices $\bbA$, $\bbB$, $\bbC$, and $\bbF$ of appropriate dimensions and any vectors $\bba$ and $\bbb$ of appropriate dimensions, it holds that
\begin{eqnarray} \label{1101.11}
(\bbA+\bbB \bbF  \bbC)^{-1}=\bbA^{-1}-\bbA^{-1}\bbB(\bbF^{-1}+\bbC\bbA^{-1}\bbB)^{-1}\bbC\bbA^{-1}
\end{eqnarray}
and
\begin{eqnarray} \label{1101.11h}
(\bbC+\bba\bbb^T)^{-1} \ba=\frac{\bbC^{-1}\bba}{1+\bbb^T\bbC^{-1}\bba}
\end{eqnarray}
when the corresponding matrices for matrix inversion are nonsingular.

To illustrate the main ideas of our proof, we first consider the simple case of $K=1$ and $\bbx=\bby=\bbv_1$. The general case of $K\geq 1$ and arbitrary unit vectors will be discussed later. Let  $\Omega_1$ be a contour centered at $(a_1+b_1)/2$  with radius $|b_1-a_1|/2$, where the quantities $a_k$ and $b_k$ with $1 \leq k \leq K$ are defined in Section \ref{Sec3.2}. Then it is seen that $d_1$ is enclosed by $\Omega_1$.   In view of Condition \ref{as3},  Lemma \ref{0505-1}, and Weyl's inequality, we have
\[ |\lambda_1-d_1|\le \|\bW\|< \min\{|d_1-a_1|,|d_1-b_1|\} \]
and
\[ |\lambda_j-d_1|\ge |d_1|- \|\bW\|> \max\{|d_1-a_1|,|d_1-b_1|\}, \ j\ge 2 \]
with significant probability. We can see that the contour $\Omega_1$ does not enclose any other eigenvalues $\lambda_j$ with $j \neq 1$. Thus, by Cauchy's residue theorem from complex analysis, we have with significant probability
\[ -\frac{1}{2\pi i}\oint_{\Omega_1}\frac{1}{\lambda_1-z} dz=1 \ \text{ and } \  -\frac{1}{2\pi i}\oint_{\Omega_1}\frac{1}{\lambda_j-z} dz=0, \ j\ge 2, \] where $i$ associated with the complex integrals represents the imaginary unit $(-1)^{1/2}$ and the line integrals are taken over the contour $\Omega_1$. Noticing that $(\bbX-z\bbI)^{-1}=\sum_{i=1}^n(\lambda_j-z)^{-1}\widehat\bv_j\widehat\bv_j^T$, we can then obtain an integral representation of the desired bilinear form that with significant probability 
\begin{align} \label{0410.1}
\nonumber \bbv^T_1\widehat \bbv_1\widehat \bbv^T_1\bbv_1&=-\frac{\bbv^T_1\widehat \bbv_1\widehat \bbv^T_1\bbv_1}{2\pi i}\oint_{\Omega_1}\frac{1}{\lambda_1-z}dz =-\frac{1}{2\pi i}\oint_{\Omega_1}\bbv^T_1\Big(\sum_{j=1}^n\frac{\widehat \bbv_j\widehat \bbv^T_j}{\lambda_j-z}\Big)\bbv_1 dz \nonumber\\
&=-\frac{1}{2\pi i}\oint_{\Omega_1}\bbv^T_1\widetilde \bbG(z)\bbv_1 dz,
\end{align}
where the matrix-valued function $\widetilde \bbG(z)=(\bbX-z\bbI)^{-1}$ for $z$ in the complex plane $\mathbb{C}$ is referred to as the Green function associated with the original random matrix $\bX = \bH + \bW$.

Note that by (\ref{model}) and $K = 1$ for the simple case, we have  $\bbX=\bbH+\bbW=d_1\bbv_1\bbv^T_1+\bbW$. Thus the line integral in
(\ref{0410.1}) can be rewritten as
\begin{eqnarray} \label{0410.1h}
 \bbv^T_1\widehat \bbv_1\widehat \bbv^T_1\bbv_1=-\frac{1}{2\pi i}\oint_{\Omega_1}\bbv^T_1(\bbW-z\bbI+d_1\bbv_1\bbv_1^T)^{-1}\bbv_1 dz.
\end{eqnarray}
With the aid of (\ref{1101.11}) and (\ref{1101.11h}), the line integral in (\ref{0410.1h}) can be further represented as
\begin{eqnarray} \label{0410.2}
 \bbv^T_1\widehat \bbv_1\widehat \bbv^T_1\bbv_1=-\frac{1}{2\pi i}\oint_{\Omega_1}\frac{\bbv^T_1(\bbW-z\bbI)^{-1}\bbv_1}{1+d_1\bbv^T_1(\bbW-z\bbI)^{-1}\bbv_1} dz.
\end{eqnarray}
To analyze the integrand of the line integral on the right hand side of (\ref{0410.2}), we first consider the term $(\bbW-z\bbI)^{-1}$. Such a term admits the matrix series expansion
\begin{eqnarray} \label{0410.31}
(\bbW-z\bbI)^{-1}=-\sum_{l=0}^{\infty}z^{-(l+1)}\bbW^l.
\end{eqnarray}
Let $L$ be the smallest positive integer such that \begin{equation} \label{eq011}
\frac{\alpha_n^{L+1}(\log n)^{(L+1)/2}}{|d_K|^{L-2}}\rightarrow 0.
	\end{equation}
Such an integer $L$ always exists since $|d_K|/(n^{\ep}\alpha_n)\rightarrow \infty$ for small positive constant $\ep$ by Condition \eqref{as3} and $\alpha_n \leq n^{1/2}$ by definition. Since we consider $z$ on the contour $\Omega_1$, it follows that $|z|\ge c|d_1|$ for some positive constant $c$. Thus,
by (\ref{0410.31}), Condition \ref{as1}, and Lemma \ref{0505-1} in Section \ref{SecB.5} of Supplementary Material, with the above choice of $L$ in (\ref{eq011}) we have with probability tending to one that
	\begin{equation} \label{0426.1}
	\left\|\sum_{l=L+1}^{\infty}z^{-(l+1)}\bbW^l\right\|\le \sum_{l=L+1}^{\infty}\frac{C^l\alpha_n^l(\log n)^{l/2}}{|z|^{l+1}}=\frac{O\{C^{L+1}\alpha_n^{L+1}(\log n)^{(L+1)/2}\}}{|z|^{L+2}}=\frac{O(1)}{|z|^4},
	\end{equation}
	where $C$ is some positive constant.
In light of (\ref{0410.31}) and (\ref{0426.1}), we can obtain the asymptotic expansion
\begin{align} \label{neweq012}
\bbv^T_1(\bbW-z\bbI)^{-1}\bbv_1 & =-\sum_{l=0}^{L}z^{-(l+1)}\bbv^T_1\bbW^l\bbv_1-\sum_{l=L+1}^{\infty}z^{-(l+1)}\bbv^T_1\bbW^l\bbv_1 \nonumber\\
&=-\sum_{l=0}^{L}z^{-(l+1)}\bbv^T_1\bbW^l\bbv_1+\frac{O_p(1)}{d_1^4}
\end{align}
for $z$ on the contour $\Omega_1$.

Directly working with the line integral in \eqref{0410.1} or \eqref{0410.2} is challenging in deriving the CLT for the bilinear form $\bv_1^T\widehat\bv_1\widehat\bv_1^T\bv_1$. Next we introduce some simple facts about Cauchy's residue theorem. Assume that a complex  function $f(z)$ is a holomorphic function inside $\Omega_1$ except at one point $t$. Then it holds that
\[ \frac{1}{2\pi i}\oint_{\Omega_1}f(z)dz=\text{Res}(f,t), \]
where $\text{Res}(f,t)$ represents the residue of function $f$ at point $t$. In addition, assume that the Laurent series expansion of $f$ around point $t$ is given by
$$f(z)=\sum_{j=-\infty}^{\infty}a_j(z-t)^j $$
with $a_j$ some constants. Then we have $\text{Res}(f,t)=(2\pi i)^{-1}\oint_{\Omega_1}f(z)dz=a_{-1}$. Furthermore, if $\lim_{z\rightarrow t}(z-t)f(z)$ exists then the Laurent series expansion of $f$ entails that
\begin{equation}\label{0927.6h}
\lim_{z\rightarrow t}(z-t)f(z)=a_{-1}.
 \end{equation}
Now let us consider the line integral in (\ref{0410.2}). Observe that the only singular point of function $\bbv^T_1(\bbW-z\bbI)^{-1}\bbv_1/[1+d_1\bbv^T_1(\bbW-z\bbI)^{-1}\bbv_1]$ inside $\Omega_1$ is the solution to equation
$$1+d_1\bbv^T_1(\bbW-z\bbI)^{-1}\bbv_1=0, 
$$
which we denote as $\widehat t_1$. Let us use $[(\bbW-\widehat t_1\bbI)^{-1}]'$ as a shorthand notation for $h'(\widehat t_1)$ with $h(t)=(\bbW-t\bbI)^{-1}$.  Then by Cauchy's residue theorem and in view of \eqref{0410.2}, we have
\begin{align*}
\bbv^T_1\widehat \bbv_1\widehat \bbv^T_1\bbv_1 & = -\frac{1}{2\pi i}\oint_{\Omega_1}\frac{\bbv^T_1(\bbW-z\bbI)^{-1}\bbv_1}{1+d_1\bbv^T_1(\bbW-z\bbI)^{-1}\bbv_1} dz=-\lim_{z\rightarrow \widehat t_1}\frac{(z-\widehat t_1)\bbv^T_1(\bbW-z\bbI)^{-1}\bbv_1}{1+d_1\bbv^T_1(\bbW-z\bbI)^{-1}\bbv_1}\\
&
=-\frac{\bbv^T_1(\bbW-\widehat t_1\bbI)^{-1}\bbv_1}{d_1\bbv^T_1[(\bbW-\widehat t_1\bbI)^{-1}]'\bbv_1}.
\end{align*}
Therefore, an application of the Taylor expansion to function  $\bbv^T_1(\bbW-\widehat t_1\bbI)^{-1}\bbv_1/\{d_1\bbv^T_1[(\bbW-\widehat t_1\bbI)^{-1}]'\bbv_1\}$ yields
\begin{equation}\label{eq010}
	-\frac{\bbv^T_1(\bbW-\widehat t_1\bbI)^{-1}\bbv_1}{d_1\bbv^T_1[(\bbW-\widehat t_1\bbI)^{-1}]'\bbv_1}=\frac{\sum_{l=0}^{L}\widehat t_1^{-(l+1)}\bbv^T_1\bbW^l\bbv_1+ d_1^{-4} O_p(1)}{d_1\sum_{l=0}^{L}(l+1)\widehat t_1^{-(l+2)}\bbv^T_1\bbW^l\bbv_1+d_1^{-4} O_p(1)}.
	\end{equation}
Note that $\widehat t_1$ is a random variable that depends on random matrix $\bbX$. 
In fact, from (\ref{0516.15}) we can see that the asymptotic expansion of $\widehat t_1$ is a polynomial of $\bbv^T_1\bbW^l\bbv_1$. Thus the asymptotic expansion of \eqref{eq010} is also a polynomial function of $\bbv^T_1\bbW^l\bbv_1$. Therefore, controlling the variance of  $\bbv^T_1\bbW^l\bbv_1$ can facilitate us in identifying the leading term of the asymptotic expansion. So far we have laid out the major steps in deriving the asymptotic expansion for $\bbv^T_1\widehat \bbv_1\widehat \bbv^T_1\bbv_1$.  This can shed light on the detailed proof for the general case of $\bbx^T\widehat \bbv_k\widehat \bbv^T_k\bby$ with $K\geq 1$.

We now move on to the general case of $K \geq 1$ and arbitrary $n$-dimensional unit vectors $\bx$ and $\by$. The technical arguments for the general case are similar to those for the simple case of $K=1$ and $\bbx=\bby=\bbv_1$ presented above, but with more delicate technical derivations. Similarly as in (\ref{0410.1}), it follows from Cauchy's residue theorem, the definitions of the eigenvalue and eigenvector, and (\ref{model}) that the bilinear form $\bbx^T\widehat\bbv_k\widehat\bbv_k^T\bby$ for each $1 \leq k \leq K$ admits a natural integral representation; that is, with significant probability,
\begin{align} \label{0410.3}
& \bbx^T \widehat \bbv_k\widehat \bbv^T_k\bby=-\frac{1}{2\pi i}\oint_{\Omega_k}\bbx^T\widetilde \bbG(z)\bby dz=-\frac{1}{2\pi i}\oint_{\Omega_k}\bbx^T\Big(\bbW-z\bbI+\sum_{j=1}^Kd_j\bbv_j\bbv^T_j\Big)^{-1}\bby dz\nonumber \\
&=\frac{1}{2\pi i}\oint_{\Omega_k} \frac{d_k \bbx^T\Big(\bbW-z\bbI+\sum\limits_{1 \leq j \neq k \leq K} d_j\bbv_j\bbv^T_j\Big)^{-1}\bbv_k\bbv_k^T\Big(\bbW-z\bbI+\sum\limits_{1 \leq j \neq k \leq K} d_i\bbv_i\bbv^T_i\Big)^{-1}\bby}{1+d_k\bbv_k^T\Big(\bbW-z\bbI+\sum\limits_{1 \leq j \neq k \leq K} d_j\bbv_j\bbv^T_j\Big)^{-1}\bbv_k} dz,
\end{align}
where the Green function $\widetilde \bbG(z)$ associated with the original random matrix $\bX$ is defined in (\ref{0410.1}) and the line integral is taken over a contour $\Omega_k$ that is centered at $(a_k+b_k)/2$ with radius $|b_k-a_k|/2$.  Then the contour $\Omega_k$ encloses the population eigenvalue $d_k$ of the latent mean matrix $\bH$. Note that in the representation above, we have used the results, which can be derived from Condition \ref{as3},  Lemma \ref{0505-1}, and Weyl's inequality, that for each $j=1,\cdots, K$,
	$$|\lambda_k-d_k|\le \|\bW\|< \min\{|d_k-a_k|,|d_k-b_k|\},$$
	$$|\lambda_j-d_k|\ge |d_j-d_k|- |\lambda_j-d_j| \ge |d_j-d_k|- \|\bW\|> \max\{|d_k-a_k|,|d_k-b_k|\} $$
	for $j\neq k$ with significant probability; that is, the contour $\Omega_k$ encloses $\lambda_k$ but not any other eigenvalues with high probability.

An application of (\ref{1101.11}) leads to
\begin{equation} \label{neweq014}
\Big(\bbW-z\bbI+\sum_{1\leq j \neq k \leq K} d_j\bbv_j\bbv^T_j\Big)^{-1}=\bbG(z)-\bbG(z)\bbV_{-k}\left[\bbD_{-k}^{-1}+\bbV_{-k}^T\bbG(z)\bbV_{-k}\right]^{-1}\bbV_{-k}^T\bbG(z),
\end{equation}
where the Green function $\bbG(z)$ associated with only the noise part $\bW$ is defined in (\ref{neweq001}). To simplify the expression, let
\begin{equation}\label{eq013}
\bbF_k(z) = \bbG(z)\bbV_{-k}[\bbD_{-k}^{-1}+\bbV_{-k}^T\bbG(z)\bbV_{-k}]^{-1}\bbV_{-k}^T\bbG(z)
\end{equation}
 Then in view of (\ref{neweq014}), the last line integral in (\ref{0410.3}) can be further represented as
\begin{align} \label{0410.4}
\bbx^T \widehat \bbv_k\widehat \bbv^T_k\bby=\frac{1}{2\pi i} \oint_{\Omega_k}\frac{d_k\bbx^T\left[\bbG(z)-\bbF_k(z)\right]\bbv_k\bbv_k^T\left[\bbG(z)-\bbF_k(z)\right]\bby}{1+d_k\bbv_k^T\left[\bbG(z)-\bbF_k(z)\right]\bbv_k} dz. 
\end{align}
 It is challenging to analyze the terms in (\ref{0410.4}) since the expression of $\bbF_k(z)$ is complicated and  we need to study the asymptotic expansion of $\bbF_k(z)$ carefully. In the proof below, we will see that Lemma  \ref{0426-1} in Section \ref{Sec5} is a key ingredient of the technical arguments; see Section \ref{SecB.4} of Supplementary Material for the proof of this lemma.

We will conduct detailed calculations for the asymptotic expansion of $\bbF_k(z)$. Let us choose  $L$  as the same positive integer as in \eqref{eq011}.
 Then we have $\sum_{l={L+1}}^{\infty} z^{-(l+1)} \bbx^T\bbW^l\bby=O_p(|z|^{-4})$ for $z$ on the contour $\Omega_k$.
  It follows from Lemma \ref{0426-1} and Condition \ref{as3} that
   \begin{align*}
  	\sum_{l=2}^L & z^{-(l+1)} \bbx^T(\bbW^l-\mathbb{E}\bbW^l)\bby =O_p\left\{\alpha_n |z|^{-3}+\alpha_n^2 |z|^{-4}+\cdots+\alpha_n^{L-1} |z|^{-(L+1)}\right\}\\
  	&\quad=O_p(\alpha_n |z|^{-3}).
  	\end{align*}
Therefore, similar to \eqref{neweq012} we can show that
\begin{align} \label{0508.2}
\bbx^T\bbG(z)\bby & =-z^{-1} \bbx^T\bby - z^{-2}\bbx^T\bbW\bby -\sum_{l=2}^L z^{-(l+1)} \bbx^T\mathbb{E}\bbW^l\bby-\sum_{l={L+1}}^\infty z^{-(l+1)} \bbx^T\bbW^l\bby  \non
&-\sum_{l=2}^L z^{-(l+1)} \bbx^T(\bbW^l-\mathbb{E}\bbW^l)\bby \nonumber \\
&=-z^{-1} \bbx^T\bby - z^{-2}\bbx^T\bbW\bby -\sum_{l=2}^L z^{-(l+1)} \bbx^T\mathbb{E}\bbW^l\bby  +O_{p}\left(|z|^{-4}+\alpha_n |z|^{-3}\right).
\end{align}
Moreover, since for $z\in \Omega_k$ we have $|z|^{-4}\le \alpha_n |z|^{-3}$ by Condition \ref{as1}, we can further obtain
\begin{equation} \label{0512.1}
\bbx^T\bbG(z)\bby=-z^{-1} \bbx^T\bby-z^{-2}\bbx^T\bbW\bby -\sum_{l=2}^L z^{-(l+1)}\bbx^T\mathbb{E}\bbW^l\bby +O_{p}(\alpha_n |z|^{-3}).
\end{equation}
In fact, the probabilistic event associated with the small order term $O_{p}(\alpha_n |z|^{-3})$ in (\ref{0512.1}) holds uniformly over $z$ since the term $O_{p}(\alpha_n |z|^{-3})$ is simply $|z|^{-3} O_{p}(\alpha_n)$.

To simplify the technical presentation, hereafter we use the generic notation $\bu$ to denote either $\bx$ or $\by$ unless specified otherwise, which means that the corresponding derivations and results hold when $\bu$ is replaced by $\bx$ and $\by$. Since $\bx$ and $\by$ can be chosen as any unit vectors,
we can obtain from (\ref{0512.1}) the following asymptotic expansions by different choices of $\bx$ and $\by$
\begin{align}
	\bbu^T\bbG(z)\bbv_k & =-z^{-1}\bbu^T\bbv_k-z^{-2}\bbu^T\bbW\bbv_k-\sum_{l=2}^L z^{-(l+1)}\bbu^T\mathbb{E}\bbW^l\bbv_k +O_p(\alpha_n |z|^{-3}), \label{0508.3} \\
	\bbv^T_k\bbG(z)\bbv_k & =-z^{-1}-z^{-2}\bbv^T_k\bbW\bbv_k -\sum_{l=2}^Lz^{-(l+1)}\bbv^T_k\mathbb{E}\bbW^l\bbv_k+O_p(\alpha_n |z|^{-3}), \label{0508.4}\\
	\bbv^T_k\bbG(z)\bbV_{-k}&=-z^{-2}\bbv^T_k\bbW\bbV_{-k}-\sum_{l=2}^Lz^{-(l+1)}\bbv^T_k\mathbb{E}\bbW^l\bbV_{-k} +O_p(\alpha_n |z|^{-3}), \label{0508.5}\\
	\bbu^T\bbG(z)\bbV_{-k} & =-z^{-1} \bbu^T\bbV_{-k}-z^{-2}\bbu^T\bbW\bbV_{-k}-\sum_{l=2}^Lz^{-(l+1)}\bbu^T\mathbb{E}\bbW^l\bbV_{-k} \nonumber\\
	&\quad +O_p(\alpha_n |z|^{-3}), \label{0508.6}\\
	\bbV^T_{-k}\bbG(z)\bbV_{-k}&=-z^{-1}\bbI-z^{-2}\bbV^T_{-k}\bbW\bbV_{-k} -\sum_{l=2}^Lz^{-(l+1)}\bbV^T_{-k}\mathbb{E}\bbW^l\bbV_{-k}+O_p(\alpha_n |z|^{-3}). \label{0508.7}
	\end{align}
Thus it follows from (\ref{0512.1})--(\ref{0508.7}) that
	\begin{align} \label{0508.8}
	\bbu^T & \bbF_k(z)\bbv_k 
	=\mathcal{R}(\bbu,\bbV_{-k},z)\left[\bbD_{-k}^{-1}+\mathcal{R}(\bbV_{-k},\bbV_{-k},z)\right]^{-1}\mathcal{R}(\bbV_{-k},\bbv_k,z)\nonumber\\
	&\quad-z^{-2}\mathcal{R}(\bbu,\bbV_{-k},z)\left[\bbD_{-k}^{-1}+\mathcal{R}(\bbV_{-k},\bbV_{-k},z)\right]^{-1}\bbV_{-k}^T\bbW\bbv_k\nonumber\\
	&\quad-z^{-2}\bbu^T\bbW\bbV_{-k}\left[\bbD_{-k}^{-1}+\mathcal{R}(\bbV_{-k},\bbV_{-k},z)\right]^{-1}\mathcal{R}(\bbV_{-k},\bbv_k,z)\nonumber\\
	&\quad+z^{-2}\mathcal{R}(\bbu,\bbV_{-k},z)\left[\bbD_{-k}^{-1}+\mathcal{R}(\bbV_{-k},\bbV_{-k},z)\right]^{-1}\bbV_{-k}^T\bbW\bbV_{-k}\nonumber\\
	&\quad\times\left[\bbD_{-k}^{-1}+\mathcal{R}(\bbV_{-k},\bbV_{-k},z)\right]^{-1}\mathcal{R}(\bbV_{-k},\bbv_k,z)+O_p(\alpha_n |z|^{-3}),\nonumber\\
	&=\mathcal{R}(\bbu,\bbV_{-k},z)\left[\bbD_{-k}^{-1}+\mathcal{R}(\bbV_{-k},\bbV_{-k},z)\right]^{-1}\mathcal{R}(\bbV_{-k},\bbv_k,z) \nonumber\\
	&\quad-z^{-2}\mathcal{R}(\bbu,\bbV_{-k},z)\left[\bbD_{-k}^{-1}+\mathcal{R}(\bbV_{-k},\bbV_{-k},z)\right]^{-1}\bbV_{-k}^T\bbW\bbv_k+O_p(\alpha_n |z|^{-3})
	\end{align}
and
\begin{align} \label{0508.9}
\bbv^T_k & \bbF_k(z)\bbv_k  =\bbv^T_k\bbG(z)\bbV_{-k}\left[\bbD_{-k}^{-1}+\bbV_{-k}^T\bbG(z)\bbV_{-k}\right]^{-1}\bbV_{-k}^T\bbG(z)\bbv_k \nonumber\\
&=\mathcal{R}(\bbv_k,\bbV_{-k},z)\left[\bbD_{-k}^{-1}+\mathcal{R}(\bbV_{-k},\bbV_{-k},z)\right]^{-1}\mathcal{R}(\bbV_{-k},\bbv_k,z) \nonumber\\
&\quad-z^{-2}\mathcal{R}(\bbv_k,\bbV_{-k},z)\left[\bbD_{-k}^{-1}+\mathcal{R}(\bbV_{-k},\bbV_{-k},z)\right]^{-1}\bbV_{-k}^T\bbW\bbv_k+O_p(\alpha_n |z|^{-3})\nonumber\\
&=\mathcal{R}(\bbv_k,\bbV_{-k},z)\left[\bbD_{-k}^{-1}+\mathcal{R}(\bbV_{-k},\bbV_{-k},z)\right]^{-1}\mathcal{R}(\bbV_{-k},\bbv_k,z)+O_p(\alpha_n |z|^{-3}),
\end{align}
where $\bbF_k(z)$ is defined in \eqref{eq013} and $\mathcal{R}$ is defined in \eqref{0619.0}.

With all the technical preparations above, we are now ready to analyze the terms in representation (\ref{0410.4}). Specifically let us consider the ratio  $\{d_k\bbx^T[\bbG(z)-\bbF_k(z)]\bbv_k\bbv_k^T[\bbG(z)-\bbF_k(z)]\bby\}/\{1+d_k\bbv_k^T\left[\bbG(z)-\bbF_k(z)\right]\bbv_k\}$ that appears as the integrand on the left hand side of (\ref{0410.4}). Similar to (\ref{0508.2}), taking the derivative of $\bbG(z)$ we have
	\begin{align}\label{1011.3h}
	\bbx^T\bbG'(z)\bby&=\bbx^T(\bbW-z\bbI)^{-2}\bby=\sum_{l=0}^{\infty}(l+1)z^{-(l+2)}\bbx^T\bbW^l\bby\non
	&=\mathcal{R}'(\bbx,\bby,z)+2z^{-3}\bbx^T\bbW\bby+z^{-4}O_p(\alpha_n).
	\end{align}
	It follows from Lemmas \ref{0426-1}--\ref{1212-1h} that
	\begin{align}\label{1011.4h}
	& \mathcal{R}'(\bbv_k,\bbV_{-k},z)=O(\alpha_n^2/z^4), \  \mathcal{R}'(\bbv_k,\bbv_k,z)-\frac{1}{z^2}=O(\alpha_n^2/z^4), \nonumber\\
	&  \left\|\mathcal{R}'(\bbV_{-k},\bbV_{-k},z)-z^{-2}\bbI\right\|=O(\alpha_n^2/z^4).
	\end{align}
	By (\ref{0508.5}) and Lemmas \ref{0426-1}--\ref{1212-1h}, we can conclude that
	\begin{align}\label{1011.1h}
	\bbv^T_k\bbG(z)\bbV_{-k}=z^{-2}O_p(1)+|z|^{-3}O_p(\alpha^2_n).
	\end{align}
	Moreover, by (\ref{0508.6}) and (\ref{1010.1h}) we have
	\begin{align}\label{1010.3h}
	\Big\| & \left\{\left[\bbD_{-k}^{-1}+\bbV_{-k}^T\bbG(z)\bbV_{-k}\right]^{-1}-\left[\bbD_{-k}^{-1}+\mathcal{R}(\bbV_{-k},\bbV_{-k},z)\right]^{-1}\right\}'\Big\|\non
	&=\Big\|\left[\bbD_{-k}^{-1}+\bbV_{-k}^T\bbG(z)\bbV_{-k}\right]^{-1}\bbV_{-k}^T\bbG'(z)\bbV_{-k} \left[\bbD_{-k}^{-1}+\bbV_{-k}^T\bbG(z)\bbV_{-k}\right]^{-1}\non
	&-\left[\bbD_{-k}^{-1}+\mathcal{R}(\bbV_{-k},\bbV_{-k},z)\right]^{-1}\mathcal{R}'(\bbV_{-k},\bbV_{-k},z)\left[\bbD_{-k}^{-1}+\mathcal{R}(\bbV_{-k},\bbV_{-k},z)\right]^{-1}\Big\|\non
	&=O\left\{\left\|\bbV_{-k}^T\bbG'(z)\bbV_{-k}-\mathcal{R}'(\bbV_{-k},\bbV_{-k},z)\right\| \left\|\left[\bbD_{-k}^{-1}+\bbV_{-k}^T\bbG(z)\bbV_{-k}\right]^{-1}\right\|^2\right\}\non
	&+O\Big\{\left\|\left[\bbD_{-k}^{-1}+\bbV_{-k}^T\bbG(z)\bbV_{-k}\right]^{-1}-\left[\bbD_{-k}^{-1}+\mathcal{R}(\bbV_{-k},\bbV_{-k},z)\right]^{-1}\right\| \nonumber \\
	&
	\quad \cdot \left\|\left[\bbD_{-k}^{-1}+\bbV_{-k}^T\bbG(z)\bbV_{-k}\right]^{-1}\mathcal{R}'(\bbV_{-k},\bbV_{-k},z)\right\|\Big\}\non
	&=|z|^{-1}O_p(1)+z^{-2}O_p(\alpha_n)
	\end{align}
	and
	\begin{align}\label{1010.6h}
	\Big\| & \left\{\left[\bbD_{-k}^{-1}+\mathcal{R}(\bbV_{-k},\bbV_{-k},z)\right]^{-1}\right\}'\Big\|\nonumber \\
	&=\left\|\left[\bbD_{-k}^{-1}+\mathcal{R}(\bbV_{-k},\bbV_{-k},z)\right]^{-1}\mathcal{R}'(\bbV_{-k},\bbV_{-k},z)\left[\bbD_{-k}^{-1}+\mathcal{R}(\bbV_{-k},\bbV_{-k},z)\right]^{-1}\right\| \nonumber\\
	&=O(1).
	\end{align}
	
	Note that in light of (\ref{1011.3h})--(\ref{1010.3h}), we can obtain
	\begin{align}\label{1011.2h}
	\bbv_k^T\bbF'_k(z)\bbv_k & =2\bbv_k^T\bbG'(z)\bbV_{-k}\left[\bbD_{-k}^{-1}+\bbV_{-k}^T\bbG(z)\bbV_{-k}\right]^{-1}\bbV_{-k}^T\bbG(z)\bbv_k\non
	&\quad +\bbv_k^T\bbG(z)\bbV_{-k}\left\{\left[\bbD_{-k}^{-1}+\bbV_{-k}^T\bbG(z)\bbV_{-k}\right]^{-1}\right\}'\bbV_{-k}^T\bbG(z)\bbv_k\non
	&= 2\mathcal{R}'(\bbv_k,\bbV_{-k},z)\left[\bbD_{-k}^{-1}+\mathcal{R}(\bbV_{-k},\bbV_{-k},z)\right]^{-1}\mathcal{R}(\bbV_{-k},\bbv_k,z)\non
	&\quad+\mathcal{R}(\bbv_k,\bbV_{-k},z)\left\{\left[\bbD_{-k}^{-1}+\mathcal{R}(\bbV_{-k},\bbV_{-k},z)\right]^{-1}\right\}'\mathcal{R}(\bbV_{-k},\bbv_k,z) \nonumber\\
	&\quad+z^{-4}O_p(1)+z^{-6}O_p(\alpha_n^3).
	\end{align}
	Combining the above result with \eqref{1011.3h} leads to
\begin{align} \label{0515.2c}
d_k\bbv_k^T\left[\bbG'(z)-\bbF'_k(z)\right]\bbv_k=\frac{d_k}{z^2\mathcal{\widetilde P}_{k,z}}+2z^{-3}d_k\bbv_k^T\bbW\bbv_k+z^{-4}O_p(|d_k|\alpha_n)
\end{align}
for $z\in [a_k, b_k]$.  
Further, recalling the definition in (\ref{0619.1}) and by (\ref{1010.6h}), it holds that
\begin{align}\label{1011.5h}
\frac{1}{z^2\mathcal{\widetilde P}_{k,z}}&=\left(\frac{A_{\sbv_k,k,z}}{z}\right)'=\mathcal{R}'(\bbv_k,\bbv_k,z)- 2\mathcal{R}'(\bbv_k,\bbV_{-k},z)\left[\bbD_{-k}^{-1}+\mathcal{R}(\bbV_{-k},\bbV_{-k},z)\right]^{-1} \nonumber\\
&\quad \times \mathcal{R}(\bbV_{-k},\bbv_k,z) -\mathcal{R}(\bbv_k,\bbV_{-k},z)\left\{\left[\bbD_{-k}^{-1}+\mathcal{R}(\bbV_{-k},\bbV_{-k},z)\right]^{-1}\right\}'\mathcal{R}(\bbV_{-k},\bbv_k,z) \nonumber\\
&=z^{-2}+O(\alpha_n^2/z^4).
\end{align}
Plugging this into \eqref{0515.2c} and by Lemmas \ref{0426-1}--\ref{1212-1h}, we have for all $z\in [a_k, b_k]$,
\begin{align}\label{0515.2b}
d_k & \bbv_k^T \left[\bbG'(z)-\bbF'_k(z)\right]\bbv_k= d_k z^{-2}+2z^{-3}d_k\bbv_k^T\bbW\bbv_k+z^{-4}O_p(|d_k|\alpha_n^2)\non
& = d_kz^{-2}\left[1+ O_p(|z|^{-1} + |z|^{-2}\alpha_n^2)\right] =  d_kz^{-2}\left[1+o_p(1)\right].
\end{align}
Thus $1+d_k\bbv_k^T\left[\bbG(z)-\bbF_k(z)\right]\bbv_k$ is a monotone  function with probability tending to one.

Further, in light of expressions (\ref{0508.4}) and (\ref{0508.9}) we can obtain the asymptotic expansion
\begin{align} \label{0515.2}
1 & +d_k\bbv_k^T\left[\bbG(z)-\bbF_k(z)\right]\bbv_k = f_k(z) - d_kz^{-2}\bbv_k^T\bbW\bbv_k + z^{-2}O_p(\alpha_n)
\end{align}
for all $z\in [a_k, b_k]$, where $f_k(z)$ is defined in (\ref{0515.3.1}).
Note that $f_k(a_k)= O(1)$, $f_k(b_k)= O(1)$, and $f_k(a_k)f_k(b_k)<0$ as shown in the proof of Lemma \ref{lem: define-t} in Section \ref{SecB.3} of  Supplementary Material. These results together with (\ref{0515.2b}), which gives the order for the derivative of $1 +d_k\bbv_k^T\left[\bbG(z)-\bbF_k(z)\right]\bbv_k$, entail that
there exists a unique solution $\widehat t_k$ to the equation
\begin{eqnarray} \label{0523.6}
1+d_k\bbv_k^T \left[\bbG(z)-\bbF_k(z) \right]\bbv_k=0
\end{eqnarray}
for $z$ in the interval $[a_k, b_k]$. Using Lemma \ref{0426-1}, we can further show that \eqref{0515.2} becomes
		\begin{equation} \label{0603.1}
1+		d_k\bbv_k^T\left[\bbG(z)-\bbF_k(z)\right]\bbv_k-f_k(z)=-\frac{d_k}{z^2}\bbv_k^T\bbW\bbv_k+O_p(|z|^{-2}\alpha_n)
		=O_p(|z|^{-1})
		\end{equation}
	for $z\in [a_k, b_k ]$. Note that $f_k(z)$ is a monotone function over $z\in [a_k, b_k ]$ as shown in the proof of Lemma \ref{lem: define-t} and (\ref{neweq026}). Thus it follows from (\ref{0523.6}) and (\ref{0603.1}) that
	\begin{equation} \label{neweq017}
	\widehat t_k-t_k=O_p(1).
	\end{equation}
	
	In fact, we can obtain a more precise order of $\widehat t_k-t_k$ than the initial one in (\ref{neweq017}).  	In view of (\ref{0515.2}) and the definition of $t_k$, we have
	\begin{equation} \label{neweq019}
	1+d_k\bbv_k^T\left[\bbG(t_k)-\bbF_k(t_k)\right]\bbv_k=-d_kt^{-2}_k\bbv^T_k\bbW\bbv_k+O_p(\alpha_n t^{-2}_k).
	\end{equation}
	By (\ref{0515.2b}) and (\ref{neweq019}), an application of the mean value theorem yields
		\begin{align} \label{neweq018}
		0 & = 1 +d_k\bbv_k^T\left[\bbG(\widehat t_k)-\bbF_k(\widehat t_k)\right]\bbv_k=1+d_k\bbv_k^T\left[\bbG(t_k)-\bbF_k(t_k)\right]\bbv_k \nonumber\\
		&\quad +d_k \widetilde{t}^{-2}_k\left[1+ O_p(|d_k|^{-1} + |d_k|^{-2}\alpha_n^2)\right](\widehat t_k-t_k),
		\end{align}
		where $\widetilde{t}_k$ is some number between $t_k$ and $\widehat t_k$. The asymptotic expansions in (\ref{neweq018}) and (\ref{neweq019}) entail further that
	\begin{equation} \label{0516.15}
	\widehat t_k-t_k=\frac{\mathfrak{t}_k^2}{t_k^2}\bbv^T_k\bbW\bbv_k+O_{p}(\alpha_nt_k^{-1})=\bbv^T_k\bbW\bbv_k+O_{p}(\alpha_nt_k^{-1}).
	\end{equation}

Now by the similar arguments as for obtaining (\ref{0927.6h}), the integral (\ref{0410.4}) can be evaluated as
\begin{align}
\nonumber \bbx^T \widehat \bbv_k\widehat \bbv^T_k\bby & = \frac{1}{2\pi i} \oint_{\Omega_k}\frac{d_k\bbx^T\left[\bbG(z)-\bbF_k(z)\right]\bbv_k\bbv_k^T\left[\bbG(z)-\bbF_k(z)\right]\bby}{1+d_k\bbv_k^T\left[\bbG(z)-\bbF_k(z)\right]\bbv_k} dz\\
&=\frac{\widehat t^2_k\bbx^T\left[\bbG(\widehat t_k)-\bbF_k(\widehat t_k)\right]\bbv_k\bbv_k^T\left[\bbG(\widehat t_k)-\bbF_k(\widehat t_k)\right]\bby}{\widehat t^2_k\bbv_k^T\left[\bbG'(\widehat t_k)-\bbF'_k(\widehat t_k)\right]\bbv_k}.\label{e020}
\end{align}
By  \eqref{0515.2c} we have
\begin{equation}\label{0927.8h}
\frac{1}{\widehat t^2_k\bbv_k^T\left[\bbG'(\widehat t_k)-\bbF'_k(\widehat t_k)\right]\bbv_k}=\mathcal{\widetilde P}_{k,\widehat t_k}-2\widehat  t_k^{-1}\mathcal{\widetilde P}^2_{k,\widehat t_k}\bbv_k^T\bbW\bbv_k+ \widehat t^{-2}_kO_p(\alpha_n)
\end{equation}
and \eqref{e020} can be written as
\begin{align} \label{0515.1}
&\bbx^T \widehat \bbv_k\widehat \bbv^T_k\bby=\frac{\widehat t^2_k\bbx^T\left[\bbG(\widehat t_k)-\bbF_k(\widehat t_k)\right]\bbv_k\bbv_k^T\left[\bbG(\widehat t_k)-\bbF_k(\widehat t_k)\right]\bby}{\widehat t^2_k\bbv_k^T\left[\bbG'(\widehat t_k)-\bbF'_k(\widehat t_k)\right]\bbv_k}\nonumber \\
& =\left[\mathcal{\widetilde P}_{k,\widehat t_k}-2{\widehat t_k}^{-1}\bbv_k^T\bbW\bbv_k+\widehat t^{-2}_kO_p\left({\alpha_n}\right)\right]\widehat t^2_k \bbx^T\left[\bbG(\widehat t_k)-\bbF_k(\widehat t_k)\right]\bbv_k\bbv_k^T \nonumber\\
&\quad \times \left[\bbG(\widehat t_k)-\bbF_k(\widehat t_k)\right]\bby.
\end{align}

Recall the definitions in (\ref{0619.0}) and (\ref{0619.1}). Then it follows from (\ref{0508.3}), (\ref{0508.8}), and  (\ref{0516.15}) that
\begin{align} \label{0515.4}
& \widehat t_k\bbu^T\left[\bbG(\widehat t_k)-\bbF_k(\widehat t_k)\right]\bbv_k=\mathcal{P}(\bbu,\bbv_k,\widehat t_k)-\mathcal{P}(\bbu,\bbV_{-k},\widehat t_k)\left[\widehat t_k\bbD_{-k}^{-1}+\mathcal{P}(\bbV_{-k},\bbV_{-k},\widehat t_k)\right]^{-1} \nonumber \\
& \quad \times \mathcal{P}(\bbV_{-k},\bbv_k,\widehat t_k)-\widehat t_k^{-1}\bbu^T\bbW\bbv_k+\widehat t_k^{-1}\mathcal{R}(\bbu,\bbV_{-k},\widehat t_k)\left[\bbD_{-k}^{-1}+\mathcal{R}(\bbV_{-k},\bbV_{-k},\widehat t_k)\right]^{-1} \nonumber\\
&\quad \times \bbV_{-k}^T\bbW\bbv_k+O_p(\alpha_n \widehat t^{-2}_k)\nonumber \\
&=\mathcal{P}(\bbu,\bbv_k,t_k)-\mathcal{P}(\bbu,\bbV_{-k},t_k)\left[t_k\bbD_{-k}^{-1}+\mathcal{P}(\bbV_{-k},\bbV_{-k},t_k)\right]^{-1} \nonumber \\
& \quad \times \mathcal{P}(\bbV_{-k},\bbv_k,t_k)-t_k^{-1}\bbu^T\bbW\bbv_k+t_k^{-1}\mathcal{R}(\bbu,\bbV_{-k},t_k)\left[\bbD_{-k}^{-1}+\mathcal{R}(\bbV_{-k},\bbV_{-k},t_k)\right]^{-1} \nonumber \\
& \quad \times \bbV_{-k}^T\bbW\bbv_k+O_p(\alpha_n t^{-2}_k),\non
&=A_{\sbu,k,t_k}-t_k^{-1}\bbb^T_{\sbu,k,t_k}\bbW\bbv_k+O_p(\alpha_n t^{-2}_k),
\end{align}
where $\bu$ stands for both $\bx$ and $\by$ as mentioned before.
Furthermore, by Lemma \ref{1212-1h} and (\ref{0516.15}) we can conclude that
\begin{equation}\label{0928.1h}
\mathcal{\widetilde P}_{k,\widehat t_k}=\mathcal{\widetilde P}_{k,t_k}+O_p(\alpha_n^2t_k^{-3}).
\end{equation}
Combining the representation (\ref{0515.1}) and asymptotic expansions (\ref{0515.4})--(\ref{0928.1h}), by Lemma \ref{0426-1} we can deduce that \eqref{e020}  can be further written as
\begin{align} \label{0515.5}
\bbx^T & \widehat \bbv_k\widehat \bbv^T_k\bby =\frac{\widehat t^2_k\bbx^T\left[\bbG(\widehat t_k)-\bbF_k(\widehat t_k)\right]\bbv_k\bbv_k^T\left[\bbG(\widehat t_k)-\bbF_k(\widehat t_k)\right]\bby}{\widehat t^2_k\bbv_k^T\left[\bbG'(\widehat t_k)-\bbF'_k(\widehat t_k)\right]\bbv_k}\non
&= \Big[ \mathcal{\widetilde P}_{k,t_k}-2t_k^{-1}\mathcal{\widetilde P}^2_{k,t_k}\bbv_k^T\bbW\bbv_k+O_p(\alpha_n t^{-2}_k)\Big] \Big[A_{\sbx,k,t_k}-t_k^{-1}\bbb^T_{\sbx,k,t_k}\bbW\bbv_k+O_p(\alpha_n t^{-2}_k)\Big] \nonumber\\
&\quad\times \Big[A_{\sby,k,t_k}-t_k^{-1}\bbb^T_{\sby,k,t_k}\bbW\bbv_k+O_p( \alpha_n t^{-2}_k)\Big] \nonumber\\
&=\Big[\mathcal{\widetilde P}_{k,t_k}-2t_k^{-1}\mathcal{\widetilde P}^2_{k,t_k}\bbv_k^T\bbW\bbv_k+O_p(\alpha_n t^{-2}_k)\Big]\times
\Big[A_{\sbx,k,t_k}A_{\sby,k,t_k}\nonumber\\
&\quad-t_k^{-1}\left(A_{\sbx,k,t_k}\bbb^T_{\sbx,k,t_k}+A_{\sby,k,t_k}\bbb^T_{\sby,k,t_k}\right)\bbW\bbv_k\nonumber\\
&\quad+t^{-2}_k\bbb^T_{\sbx,k,t_k}\bbW\bbv_k\bbb^T_{\sby,k,t_k}\bbW\bbv_k+O_{p}( \alpha_n c_k t^{-2}_k)\Big],
\end{align}
where $c_k=|A_{\sbx,k,t_k}|+|A_{\sby,k,t_k}|+|t_k|^{-1}$.
	
We can expand (\ref{0515.5}), or equivalently \eqref{e020}, further as
\begin{align} \label{0516.1}
&\bbx^T \widehat \bbv_k\widehat \bbv^T_k\bby  =\Big[ \mathcal{\widetilde P}_{k,t_k}-2t_k^{-1}\mathcal{\widetilde P}^2_{k,t_k}\bbv_k^T\bbW\bbv_k+O_p(\alpha_nt^{-2}_k)\Big]\times
\Big[A_{\sbx,k,t_k}A_{\sby,k,t_k} \nonumber\\
&\quad-t_k^{-1}\left(A_{\sbx,k,t_k}\bbb^T_{\sbx,k,t_k}+A_{\sby,k,t_k}\bbb^T_{\sby,k,t_k}\right)\bbW\bbv_k+t^{-2}_k\bbb^T_{\sbx,k,t_k}\bbW\bbv_k\bbb^T_{\sby,k,t_k}\bbW\bbv_k+O_{p}( \alpha_n c_k t^{-2}_k)\Big] \nonumber\\
&=\mathcal{\widetilde P}_{k,t_k}A_{\sbx,k,t_k}A_{\sby,k,t_k}-t_k^{-1} A_{\sbx,k,t_k}\mathcal{\widetilde P}_{k,t_k}\left(\bbb^T_{\sby,k,t_k}+A_{\sby,k,t_k}\mathcal{\widetilde P}_{k,t_k}\bbv_k^T\right)\bbW\bbv_k \nonumber\\
&\quad-t_k^{-1} A_{\sby,k,t_k}\mathcal{\widetilde P}_{k,t_k}\left(\bbb^T_{\sbx,k,t_k}+A_{\sbx,k,t_k}\mathcal{\widetilde P}_{k,t_k}\bbv_k^T\right)\bbW\bbv_k \nonumber\\
&\quad+t^{-2}_k\mathcal{\widetilde P}_{k,t_k}\left[2\mathcal{\widetilde P}_{k,t_k}\left(A_{\sbx,k,t_k}\bbb^T_{\sbx,k,t_k}+A_{\sby,k,t_k}\bbb^T_{\sby,k,t_k}\right)\bbW\bbv_k\bbv_k^T+\bbb^T_{\sbx,k,t_k}\bbW\bbv_k\bbb^T_{\sby,k,t_k}\right]\bbW\bbv_k \nonumber\\
&\quad-2t^{-3}_k\mathcal{\widetilde P}^2_{k,t_k}\bbb^T_{\sbx,k,t_k}\bbW\bbv_k\bbb^T_{\sby,k,t_k}\bbW\bbv_k\bbv_k^T\bbW\bbv_k+O_{p}\left\{\alpha_nc_k t^{-2}_k\right\}.
\end{align}
Therefore, we have characterized the terms involving $t_k^{-1}$ for the desired first order asymptotic expansion. That is, by (\ref{0516.1}) we have
  \begin{align}\label{0516.1a}
  &\bbx^T\widehat \bbv_k\widehat \bbv^T_k\bby=\mathcal{\widetilde P}_{k,t_k}A_{\sbx,k,t_k}A_{\sby,k,t_k}-t_k^{-1} A_{\sbx,k,t_k}\mathcal{\widetilde P}_{k,t_k}\left(\bbb^T_{\sby,k,t_k}+A_{\sby,k,t_k}\mathcal{\widetilde P}_{k,t_k}\bbv_k^T\right)\bbW\bbv_k \nonumber\\
&\quad-t_k^{-1} A_{\sby,k,t_k}\mathcal{\widetilde P}_{k,t_k}\left(\bbb^T_{\sbx,k,t_k}+A_{\sbx,k,t_k}\mathcal{\widetilde P}_{k,t_k}\bbv_k^T\right)\bbW\bbv_k+O_{p}\left\{(\alpha_nc_k+1) t^{-2}_k\right\}.
\end{align}
Thus if $\sigma_k^2=t^{-2}_k\mathcal{\widetilde P}^2_{k,t_k}\mathbb{E}[(A_{\sbx,k,t_k}\bbb^T_{\sby,k,t_k}+A_{\sby,k,t_k}\bbb^T_{\sbx,k,t_k}+2A_{\sbx,k,t_k}A_{\sby,k,t_k}\mathcal{\widetilde P}_{k,t_k}\bbv_k^T)\bbW\bbv_k]^2\gg (\alpha_nc_k+1)^2 t^{-4}_k\sim \sigma^2_n (|A_{\sbx,k,t_k}|+|A_{\sby,k,t_k}|)^2 t^{-4}_k+t^{-4}_k$ and $(A_{\sbx,k,t_k}\bbb^T_{\sby,k,t_k}+A_{\sby,k,t_k}\bbb^T_{\sbx,k,t_k}+2A_{\sbx,k,t_k}A_{\sby,k,t_k}\mathcal{\widetilde P}_{k,t_k}\bbv_k^T, \bbv_k)$ is $\bW^1$-CLT, then (\ref{0511.2}) holds, where $\sim$ means the asymptotic order. This concludes the proof of Theorem \ref{0511-1}.

\subsection{Proof of Theorem \ref{0518-1}} \label{SecA.6}
We have characterized the first order asymptotic expansion for the bilinear form $\bbx^T\widehat\bbv_k\widehat\bbv_k^T\bby$ in the proof of Theorem \ref{0511-1} in Section \ref{SecA.5}, where $\bx$ and $\by$ are two arbitrary $n$-dimensional unit vectors. We now proceed with investigating the higher order (which is second order) asymptotic expansion for the same bilinear form. More specifically, the proof of Theorem \ref{0518-1} involves further expansion for the $O_{p}\{\alpha_nc_k t^{-2}_k\}$ term given in (\ref{0516.1}).

To gain some intuition, let us recall (\ref{0508.2}) and compare with  (\ref{0508.3})--(\ref{0508.7}). By Lemma \ref{0426-1}, we see that the order $O_p(\alpha_n|z|^{-3})$  comes from the terms of form $\bbx^T(\bbW^2-\mathbb{E}\bbW^2)\bby/z^3$. Therefore, to obtain a higher order expansion we need to identify all terms of form $\bbx^T(\bbW^2-\mathbb{E}\bbW^2)\bby/z^3$. It follows from (\ref{0508.2}) and  Lemmas \ref{0426-1} and \ref{1212-1h} that
\begin{align} \label{0928.2h}
\bbx^T\bbG(z)\bby 
=& -z^{-1} \bbx^T\bby - z^{-2}\bbx^T\bbW\bby-\frac{\bbx^T(\bbW^2-\mathbb{E}\bbW^2)\bby}{z^3} \nonumber \\
&-\sum_{l=2}^L z^{-(l+1)} \bbx^T\mathbb{E}\bbW^l\bby +O_{p}\left(|z|^{-4}+\alpha_n^2 |z|^{-4}\right).
\end{align}
Moreover, using similar arguments as for proving (\ref{0927.8h}) and \eqref{0515.4} but expanding to higher orders we can obtain
\begin{align} \label{0517.1}
&\left\{\widehat t^2_k\bbv_k^T\left[\bbG'(\widehat t_k)-\bbF'_k(\widehat t_k)\right]\bbv_k\right\}^{-1}=\mathcal{\widetilde P}_{k,t_k}  \Big\{1-2 t_k^{-1} \mathcal{\widetilde P}_{k,t_k}\bbv_k^T\bbW\bbv_k-t^{-2}_k\mathcal{\widetilde P}_{k,t_k} \nonumber \\
&\quad\times \left[3\bbv_k^T(\bbW^2-\mathbb{E}\bbW^2)\bbv_k-2(\bbv_k^T\bbW\bbv_k)^2\right]\Big\} +O_p(\alpha^2_n|t_k|^{-3})
\end{align}
and
\begin{align}
&\widehat t_k\bbu^T\left[\bbG(\widehat t_k)-\bbF_k(\widehat t_k)\right]\bbv_k=A_{\sbu,k,t_k}-t_k^{-1}\bbu^T\bbW\bbv_k \nonumber\\
&\quad+t_k^{-1}\mathcal{R}(\bbu,\bbV_{-k},t_k)\left[\bbD_{-k}^{-1}+\mathcal{R}(\bbV_{-k},\bbV_{-k},t_k)\right]^{-1}\bbV_{-k}^T\bbW\bbv_k+t^{-2}_k\bbu^T\bbW\bbv_k\bbv^T_k\bbW\bbv_k \nonumber\\
&\quad-t^{-2}_k\bbv^T_k\bbW\bbv_k\mathcal{R}(\bbu,\bbV_{-k},t_k)\left[\bbD_{-k}^{-1}+\mathcal{R}(\bbV_{-k},\bbV_{-k},t_k)\right]^{-1}\bbV_{-k}^T\bbW\bbv_k \nonumber\\
&\quad+t^{-2}_k\mathcal{R}(\bbu,\bbV_{-k},t_k)\left[\bbD_{-k}^{-1}+\mathcal{R}(\bbV_{-k},\bbV_{-k},t_k)\right]^{-1}\bbV_{-k}^T(\bbW^2-\mathbb{E}\bbW^2)\bbv_k\nonumber \\
&\quad-t^{-2}_k\bbu^T(\bbW^2-\mathbb{E}\bbW^2)\bbv_k+2t_k^{-3}\bbv^T_k\bbW\bbv_k\mathcal{R}(\bbu,\bbV_{-k},t_k)\left[\bbD_{-k}^{-1}+\mathcal{R}(\bbV_{-k},\bbV_{-k},t_k)\right]^{-2}\nonumber\\
&\quad\times \bbV_{-k}^T\bbW\bbv_k+O_p(\alpha^2_n |t_k|^{-3}),\label{0517.2}
\end{align}
where $\bu$ represents both $\bx$ and $\by$ as mentioned before.

Using the representations \eqref{e020} and (\ref{0515.1}), and by the asymptotic expansions (\ref{0517.1})--(\ref{0517.2}), we can obtain the $O_p(t_k^{-2})$ term for the desired second order asymptotic expansion as follows
\begin{align} \label{0517.3}
&\bbx^T\widehat \bbv_k\widehat \bbv^T_k\bby=\frac{\widehat t^2_k\bbx^T\left[\bbG(\widehat t_k)-\bbF_k(\widehat t_k)\right]\bbv_k\bbv_k^T\left[\bbG(\widehat t_k)-\bbF_k(\widehat t_k)\right]\bby}{\widehat t^2_k\bbv_k^T\left[\bbG'(\widehat t_k)-\bbF'_k(\widehat t_k)\right]\bbv_k}\non
&= \Big(\mathcal{\widetilde P}_{k,t_k}\times \left\{1-2 t_k^{-1} \mathcal{\widetilde P}_{k,t_k}\bbv_k^T\bbW\bbv_k-t^{-2}_k\mathcal{\widetilde P}_{k,t_k}\left[3\bbv_k^T(\bbW^2-\mathbb{E}\bbW^2)\bbv_k-2(\bbv_k^T\bbW\bbv_k)^2\right]\right\} \nonumber\\
& \quad+O_p(\alpha^2_n|t_k|^{-3})\Big) \Big[\widehat t_k\bbx^T\left[\bbG(\widehat t_k)-\bbF_k(\widehat t_k)\right]\bbv_k\Big] \Big[\widehat t_k\bbv_k^T\left[\bbG(\widehat t_k)-\bbF_k(\widehat t_k)\right]\bby\Big] \nonumber\\
&=-A_{\sbx,k,t_k}\mathcal{\widetilde P}_{k,t_k}t_k^{-1}\left(\bbb^T_{\sby,k,t_k}+A_{\sby,k,t_k}\mathcal{\widetilde P}_{k,t_k}\bbv_k^T\right)\bbW\bbv_k-A_{\sby,k,t_k}\mathcal{\widetilde P}_{k,t_k}t_k^{-1} \nonumber\\
&\quad \times\left(\bbb^T_{\sbx,k,t_k}+A_{\sbx,k,t_k}\mathcal{\widetilde P}_{k,t_k}\bbv_k^T\right)\bbW\bbv_k\nonumber\\
&\quad+\mathcal{\widetilde P}_{k,t_k}t^{-2}_k\left[2\mathcal{\widetilde P}_{k,t_k}\left(A_{\sbx,k,t_k}\bbb^T_{\sbx,k,t_k}+A_{\sby,k,t_k}\bbb^T_{\sby,k,t_k}\right)\bbW\bbv_k\bbv_k^T+\bbb^T_{\sbx,k,t_k}\bbW\bbv_k\bbb^T_{\sby,k,t_k}\right]\bbW\bbv_k\nonumber\\
&\quad+2A_{\sbx,k,t_k}A_{\sby,k,t_k}\left(\bbv_k^T\bbW\bbv_k\right)^2+A_{\sby,k,t_k}\mathcal{\widetilde P}_{k,t_k}\Big\{t^{-2}_k\bbx^T\bbW\bbv_k\bbv^T_k\bbW\bbv_k-t^{-2}_k\bbv^T_k\bbW\bbv_k\mathcal{R}(\bbx,\bbV_{-k},t)\nonumber\\
&\quad\times\left[\bbD_{-k}^{-1}+\mathcal{R}(\bbV_{-k},\bbV_{-k},t_k)\right]^{-1}\bbV_{-k}^T\bbW\bbv_k\Big\}\nonumber \\
&\quad+A_{\sbx,k,t_k}\mathcal{\widetilde P}_{k,t_k}\Big\{t^{-2}_k\bby^T\bbW\bbv_k\bbv^T_k\bbW\bbv_k-t^{-2}_k\bbv^T_k\bbW\bbv_k\mathcal{R}(\bby,\bbV_{-k},t)\nonumber\\
&\quad \times\left[\bbD_{-k}^{-1}+\mathcal{R}(\bbV_{-k},\bbV_{-k},t_k)\right]^{-1}\bbV_{-k}^T\bbW\bbv_k\Big\}\nonumber \\
&\quad+A_{\sby,k,t_k}\mathcal{\widetilde P}_{k,t_k}t^{-2}_k\mathcal{R}(\bbx,\bbV_{-k},t_k)\left[\bbD_{-k}^{-1}+\mathcal{R}(\bbV_{-k},\bbV_{-k},t_k)\right]^{-1}\bbV_{-k}^T(\bbW^2-\mathbb{E}\bbW^2)\bbv_k\nonumber \\
&\quad+A_{\sbx,k,t_k}\mathcal{\widetilde P}_{k,t_k}t^{-2}_k\mathcal{R}(\bby,\bbV_{-k},t_k)\left[\bbD_{-k}^{-1}+\mathcal{R}(\bbV_{-k},\bbV_{-k},t_k)\right]^{-1}\bbV_{-k}^T(\bbW^2-\mathbb{E}\bbW^2)\bbv_k\nonumber \\
&\quad-\mathcal{\widetilde P}_{k,t_k}t^{-2}_k(A_{\sby,k,t_k}\bbx^T+A_{\sbx,k,t_k}\bby^T)(\bbW^2-\mathbb{E}\bbW^2)\bbv_k\nonumber\\
&\quad-3t^{-2}_kA_{\sbx,k,t_k}A_{\sby,k,t_k}\mathcal{\widetilde P}_{k,t_k}\bbv_k^T(\bbW^2-\mathbb{E}\bbW^2)\bbv_k+O_{p}\left\{(\alpha_n^2c_k+1) |t_k|^{-3}\right\}.
\end{align}
In contrast to the small order term $O_{p}\{\alpha_nc_k t^{-2}_k\}$ in \eqref{0516.1} from the first order asymptotic expansion, we now have the small order term $O_{p}\{(\alpha_n^2c_k+1) |t_k|^{-3}\}$ from the second order asymptotic expansion.

Let us conduct some simplifications for the expressions given in the above asymptotic expansions in \eqref{0517.3}. A combination of (\ref{0516.1}) and (\ref{0517.3}) shows that the asymptotic distribution is determined by
\begin{align} \label{0605.1}
&-A_{\sbx,k,t_k}\mathcal{\widetilde P}_{k,t_k}t_k^{-1}\left(\bbb^T_{\sby,k,t_k}+A_{\sby,k,t_k}\mathcal{\widetilde P}_{k,t_k}\bbv_k^T\right)\bbW\bbv_k\nonumber\\
& \quad-A_{\sby,k,t_k}\mathcal{\widetilde P}_{k,t_k}t_k^{-1}\left(\bbb^T_{\sbx,k,t_k}+A_{\sbx,k,t_k}\mathcal{\widetilde P}_{k,t_k}\bbv_k^T\right)\bbW\bbv_k\nonumber\\
&\quad+\mathcal{\widetilde P}_{k,t_k}t^{-2}_k\Big[2\mathcal{\widetilde P}_{k,t_k}\left(A_{\sbx,k,t_k}\bbb^T_{\sbx,k,t_k}+A_{\sby,k,t_k}\bbb^T_{\sby,k,t_k}\right)\bbW\bbv_k\bbv_k^T+\bbb^T_{\sbx,k,t_k}\bbW\bbv_k\bbb^T_{\sby,k,t_k}\Big]\bbW\bbv_k\nonumber\\
&\quad+2A_{\sbx,k,t_k}A_{\sby,k,t_k}(\bbv_k^T\bbW\bbv_k)^2+A_{\sby,k,t_k}\mathcal{\widetilde P}_{k,t_k}\Big\{t^{-2}_k\bbx^T\bbW\bbv_k\bbv^T_k\bbW\bbv_k-t^{-2}_k\bbv^T_k\bbW\bbv_k\mathcal{R}(\bbx,\bbV_{-k},t)\nonumber\\
&\quad\times\left[\bbD_{-k}^{-1}+\mathcal{R}(\bbV_{-k},\bbV_{-k},t_k)\right]^{-1}\bbV_{-k}^T\bbW\bbv_k\Big\}\nonumber\\
&\quad+A_{\sbx,k,t_k}\mathcal{\widetilde P}_{k,t_k}\Big\{t^{-2}_k\bby^T\bbW\bbv_k\bbv^T_k\bbW\bbv_k-t^{-2}_k\bbv^T_k\bbW\bbv_k\mathcal{R}(\bby,\bbV_{-k},t)\nonumber\\
&\quad\times\left[\bbD_{-k}^{-1}+\mathcal{R}(\bbV_{-k},\bbV_{-k},t_k)\right]^{-1}\bbV_{-k}^T\bbW\bbv_k\Big\}\nonumber\\
&\quad+A_{\sby,k,t_k}\mathcal{\widetilde P}_{k,t_k}t^{-2}_k\mathcal{R}(\bbx,\bbV_{-k},t) \left[\bbD_{-k}^{-1}+\mathcal{R}(\bbV_{-k},\bbV_{-k},t_k)\right]^{-1}\bbV_{-k}^T(\bbW^2-\mathbb{E}\bbW^2)\bbv_k\nonumber\\
&\quad+A_{\sbx,k,t_k}\mathcal{\widetilde P}_{k,t_k} t^{-2}_k\mathcal{R}(\bby,\bbV_{-k},t)\left[\bbD_{-k}^{-1}+\mathcal{R}(\bbV_{-k},\bbV_{-k},t_k)\right]^{-1}\bbV_{-k}^T(\bbW^2-\mathbb{E}\bbW^2)\bbv_k\nonumber\\
&\quad-\mathcal{\widetilde P}_{k,t_k}t^{-2}_k\left(A_{\sby,k,t_k}\bbx^T+A_{\sbx,k,t_k}\bby^T\right)(\bbW^2-\mathbb{E}\bbW^2)\bbv_k\nonumber\\
&\quad-3t^{-2}_kA_{\sbx,k,t_k}A_{\sby,k,t_k}\mathcal{\widetilde P}_{k,t_k}\bbv_k^T(\bbW^2-\mathbb{E}\bbW^2)\bbv_k.
\end{align}
To further simplify the notation, we define three terms
\begin{align}
\bbJ_{\sbx,\sby,k,t_k}&=-\mathcal{\widetilde P}_{k,t_k}t_k^{-1}\bbv_k\left(A_{\sby,k,t_k}\bbb^T_{\sbx,k,t_k}+A_{\sbx,k,t_k}\bbb^T_{\sby,k,t_k}+2A_{\sbx,k,t_k}A_{\sby,k,t_k}\mathcal{\widetilde P}_{k,t_k}\bbv_k^T\right), \label{0605.2}\\
\bbL_{\sbx,\sby,k,t_k}&=\mathcal{\widetilde P}_{k,t_k}t^{-2}_k\bbv_k\Big\{\left[A_{\sby,k,t_k}\mathcal{R}(\bbx,\bbV_{-k},t)+A_{\sbx,k,t_k}\mathcal{R}(\bby,\bbV_{-k},t)\right]\nonumber\\
&\quad\times\left[\bbD_{-k}^{-1}+\mathcal{R}(\bbV_{-k},\bbV_{-k},t_k)\right]^{-1}\bbV_{-k}^T+A_{\sby,k,t_k}\bbx^T+A_{\sbx,k,t_k}\bby^T\nonumber\\
&\quad+3A_{\sbx,k,t_k}A_{\sby,k,t_k}\bbv_k^T\Big\}, \label{0605.3}\\
\bbQ_{\sbx,\sby,k,t_k}&=\bbL_{\sbx,\sby,k,t_k}-\mathcal{\widetilde P}_{k,t_k}t^{-2}_kA_{\sbx,k,t_k}A_{\sby,k,t_k}\bbv_k\bbv_k^T \nonumber\\
&\quad+2\mathcal{\widetilde P}_{k,t_k}^2t_k^{-2}\bbv_k\left(A_{\sbx,k,t_k}\bbb^T_{\sbx,k,t_k}+A_{\sby,k,t_k}\bbb^T_{\sby,k,t_k}\right). \label{0605.4}
\end{align}
Note that all the three matrices defined in (\ref{0605.2})--(\ref{0605.4}) are of rank one and the identity $\bbx^T\bA\bby=\tr(\bA\bby\bbx^T)$ holds for any matrix $\bA$ and vectors $\bx$ and $\by$. Thus in view of (\ref{0605.2})--(\ref{0605.4}), the lengthy expression given in (\ref{0605.1}) can be rewritten in a compact form as
\begin{equation} \label{0605.5}
\tr\left[\bbW\bbJ_{\sbx,\sby,k,t_k} - \left(\bbW^2-\mathbb{E}\bbW^2\right)\bbL_{\sbx,\sby,k,t_k}\right]+\tr\left(\bbW\bbv_k\bbv_k^T\right)\tr\left(\bbW\bbQ_{\sbx,\sby,k,t_k}\right).
\end{equation}

So far we have shown that the second order expansion of $\bbx^T\widehat\bbv_k\widehat\bbv_k^T\bby$ is given in (\ref{0517.3}). 
Note that $\widetilde \sigma_k^2$ defined in (\ref{0619.5}) is essentially the variance of (\ref{0605.5}). Thus if $\widetilde \sigma_k^2\gg (\alpha_n^2c_k+1)^2 t_k^{-6}\sim \sigma^4_n (|A_{\sbx,k,t_k}|+|A_{\sby,k,t_k}|)^2 t^{-6}_k+t^{-6}_k$, then (\ref{0605.5}) is the leading term of  (\ref{0517.3}). Furthermore, the assumption of  $\sigma_k^2=O(\widetilde \sigma_k^2)$ entails that the first order expansion in Theorem \ref{0511-1} does not dominate the second order expansion. Therefore, we see that the asymptotic distribution in Theorem \ref{0518-1} is determined by the joint distribution of the three random variables specified in expression (\ref{0605.5}). This completes the proof of Theorem \ref{0518-1}.

\bibliographystyle{chicago}
\bibliography{references}


\newpage
\appendix
\setcounter{page}{1}
\setcounter{section}{1}
\renewcommand{\theequation}{A.\arabic{equation}}
\renewcommand{\thesubsection}{B.\arabic{subsection}}
\setcounter{equation}{0}

\begin{center}{\bf \Large Supplementary Material to ``Asymptotic Theory of Eigenvectors for Random Matrices with Diverging Spikes"}
	
\bigskip
	
Jianqing Fan, Yingying Fan, Xiao Han and Jinchi Lv
\end{center}

\noindent This Supplementary Material contains additional technical details. In particular, we present in Section \ref{SecB} the proofs of all the lemmas and provide in Section \ref{SecC} some further technical details on under what regularity conditions the asymptotic normality can hold for the asymptotic expansion in Theorem \ref{0518-1}. Section \ref{sec:cond2ii} contains the technical details on relaxing the spike strength condition when considering scenario ii) of Condition \ref{as3} in place of scenario i), as well as the proof sketch for results in Section \ref{sec: app2}.

\section{Proofs of technical lemmas} \label{SecB}

\subsection{Proof of Lemma \ref{0524-1}} \label{SecB.1}

Let $\bx = (x_1, \cdots, x_n)^T$ and $\by = (y_1, \cdots, y_n)^T$ be two arbitrary $n$-dimensional unit vectors. Since $\bW$ is a symmetric random matrix of independent entries above the diagonal, it is easy to show that
\begin{equation} \label{neweq020}
\bbx^T\bbW\bby-\bbx^T\mathbb{E}\bbW\bby=\sum_{1 \leq i, j \leq n,\, i < j}w_{ij}(x_iy_j+x_jy_i)+\sum_{1 \le i \leq n}(w_{ii}-\mathbb{E}w_{ii})(x_iy_i)
\end{equation}
and
\begin{equation} \label{neweq021}
s^2_n\equiv\mathbb{E}(\bbx^T\bbW\bby-\bbx^T\mathbb{E}\bbW\bby)^2=\sum_{1 \leq i, j \leq n,\, i < j}\mathbb{E}w^2_{ij}(x_iy_j+x_jy_i)^2+\sum_{1 \le i \leq n}\mathbb{E}(w_{ii}-\mathbb{E}w_{ii})^2x^2_iy^2_i.
\end{equation}
Since $w_{ij}$ with $1\le i<j\le n$ and $w_{ii}-\mathbb{E}w_{ii}$ with $1\le i\le n$ are independent random variables with zero mean, by the Lyapunov condition (see, for example, Theorem 27.3 of \cite{B1995}) we can see that if
$$\frac{1}{s_n^3}\left[\sum_{1 \leq i, j \leq n,\, i < j}\mathbb{E}|w_{ij}|^3|x_iy_j+x_jy_i|^3+\sum_{1 \le i \leq n}\mathbb{E}|w_{ii}-\mathbb{E}w_{ii}|^3|x_iy_i|^3\right]\rightarrow 0,$$
then it  holds that
$$\frac{\bbx^T\bbW\bby-\bbx^T\mathbb{E}\bbW\bby}{s_n}\toD N(0,1).$$

Since by assumption $\max_{1\le i, j\le n}|w_{ij}|\le 1$ and $\|\bbx\|_\infty \|\bby\|_{\infty} \ll s_n$, we have
\begin{align} \label{neweq022}
	\frac{1}{s_n^3} & \left[\sum_{1 \leq i, j \leq n,\, i < j}\mathbb{E}|w_{ij}|^3|x_iy_j+x_jy_i|^3+\sum_{1 \le i \leq n}\mathbb{E}|w_{ii}-\mathbb{E}w_{ii}|^3|x_iy_i|^3\right]\non
	& \quad \le  \frac{2}{s_n^3}\left[\sum_{1 \leq i, j \leq n,\, i < j}\mathbb{E}|w_{ij}|^2|x_iy_j+x_jy_i|^3+\sum_{1 \le i \leq n}\mathbb{E}|w_{ii}-\mathbb{E}w_{ii}|^2|x_iy_i|^3\right] \nonumber \\
	& \quad
	\ll \frac{2s_n}{s_n^3}\left[\sum_{1 \leq i, j \leq n,\, i < j}\mathbb{E}|w_{ij}|^2|x_iy_j+x_jy_i|^2+\sum_{1 \le i \leq n}\mathbb{E}|w_{ii}-\mathbb{E}w_{ii}|^2|x_iy_i|^2\right] \leq 2,
	\end{align}
which completes the proof of Lemma \ref{0524-1}.

\subsection{Proof of Lemma \ref{0524-1h}} \label{SecB.2}

The technical arguments for the proof of Lemma \ref{0524-1h} are similar to those for the proof of Lemma \ref{0524-1} in Section \ref{SecB.1}. For the case of $\bbx^T(\bbW^2-\mathbb{E}\bbW^2)\bby$, let us first consider the term $\bbx^T\bbW^2\bby$. Such a term can be written as
\begin{align} \label{0614.5}
& \sum_{1\le k,i,l\le n} w_{ki}w_{il}x_ky_l=\sum_{1\le k,i,l\le n, \, k>l}w_{ki}w_{il}(x_ky_l+x_ly_k)+\sum_{1\le k,i\le n}w_{ki}^2x_ky_k \nonumber\\
&=\sum_{ 1\le k,i,l\le n,\, k>l,\, k<i}w_{ki}w_{il}(x_ky_l+x_ly_k)+\sum_{1\leq k,i,l \leq n,\, k>l,\, k>i}w_{ki}w_{il}(x_ky_l+x_ly_k)\nonumber\\
&\quad+\sum_{1 \leq l < k \leq n}w_{kk}w_{kl}(x_ky_l+x_ly_k)+\sum_{1\le k,i\le n}w_{ki}^2x_ky_k\nonumber\\
&=\sum_{1\le k,i,l\le n,\, k>l,\, k<i}w_{ki}w_{il}(x_ky_l+x_ly_k)+\sum_{1\le k,i,l\le n,\, i>l,\, i>k}w_{ik}w_{kl}(x_iy_l+x_ly_i)\nonumber\\
&\quad+\sum_{1 \leq l < k \leq n}w_{kk}w_{kl}(x_ky_l+x_ly_k)+\sum_{1\le k,i\le n}w_{ki}^2x_ky_k\nonumber\\
&=\sum_{1\le k<i\le n}w_{ki}\Big(x_k\sum_{1 \leq l < k \leq n}w_{il}y_l+y_k\sum_{1\le l<k\le n}w_{il}x_l+x_i\sum_{1\le l<i\le n}w_{kl}y_l+y_i\sum_{1\le l<i\le n}w_{kl}x_l\Big)\nonumber\\
&\quad+\sum_{ 1\le l<k\le n}w_{kk}w_{kl}(x_ky_l+x_ly_k)+\sum_{1\le k<i\le n}w_{ki}^2(x_ky_k+x_iy_i)+\sum_{1\leq k \leq n}w_{kk}^2x_ky_k.
\end{align}
Then it follows from (\ref{0614.5}) and the independence of entries $w_{ki}$ with $1 \leq k\le i \leq n$ that
\[
\mathbb{E}\bbx^T\bbW^2\bby=\sum_{1\le k,i\le n,\, k<i}\mathbb{E}w_{ki}^2(x_ky_k+x_iy_i)+ \sum_{1 \leq k \leq n}\mathbb{E}w_{kk}^2x_ky_k.
\]

To ease the technical presentation, let us define some new notation  $\omega_{kk}=2^{-1}w_{kk}$ and $\sigma_{kk}^2=\mathbb{E}\omega_{kk}^2$. We can further show that
\begin{align}\label{0524.2}
&\bbx^T(\bbW^2-\mathbb{E}\bbW^2)\bby =\sum_{1\le k,i\le n,\,k<i}w_{ki}\Big[x_k\sum_{1\le l<k\le n}w_{il}y_l+y_k\sum_{1\le l<k\le n}w_{il}x_l+x_i\sum_{1\leq l<i\le n}w_{kl}y_l\nonumber\\
&\quad+y_i\sum_{1\le l<i\le n}w_{kl}x_l+\mathbb{E}w_{ii}(x_iy_k+x_ky_i)\Big]+\sum_{1\le k,i\le n,\,k< i}\Big[(w_{ki}^2-\sigma_{ki}^2)(x_ky_k+x_iy_i)\nonumber\\
&\quad+2(\omega_{kk}^2-\sigma_{kk}^2)(x_ky_k+x_iy_i)\Big]+\sum_{1\le k\le n} 2 (\omega_{kk}-\mathbb{E}\omega_{kk})\Big(x_k\sum_{1\le l<k\le n}w_{kl}y_l \nonumber\\
&\quad+y_k\sum_{1\le l<k\le n}w_{kl}x_l\Big),
\end{align}
where $\sigma_{ki}^2=\mathbb{E}w^2_{ki}$ denotes the variance of entry $w_{ki}$ as defined before.

We next define a $\sigma$-algebra $\mathcal{F}_{t}=\sigma\{\mathfrak{w}_1,\cdots,\mathfrak{w}_t\}$, where $\mathfrak{w}_t=w_{kl}$ with $t=k+2^{-1}l(l-1)$ and $1\le k\le l\le n$. Clearly we have $t\le 2^{-1} n(n+1)$. In fact, there is a one to one correspondence between $t\le 2^{-1} n(n+1)$ and $(k,l)$ with $k\le l$. Suppose that such a statement is not true. Then there exist two different pairs  $(k_1,l_1)$ and $(k_2,l_2)$ with $1\le k_1\le l_1\le n$ and $1\le k_2\le l_2\le n$ such that
\begin{equation}\label{0928.17h}
k_1+\frac{l_1(l_1-1)}{2}=k_2+\frac{l_2(l_2-1)}{2}.
\end{equation}
It is easy to see that we must have $k_1\neq k_2$ and $l_1\neq l_2$. Without loss of generality, let us assume that $l_1<l_2$. Then by (\ref{0928.17h}), it holds that
\[
\frac{l_2(l_2-1)}{2}-\frac{l_1(l_1-1)}{2}=k_1-k_2\le k_1-1.
\]
On the other hand, since $l_1<l_2$ we have
\[ \frac{l_2(l_2-1)}{2}-\frac{l_1(l_1-1)}{2}\ge \frac{l_1(l_1+1)}{2}-\frac{l_1(l_1-1)}{2}\ge l_1\ge k_1, \]
which contradicts the previous inequality. Thus we have shown that there is indeed a one to one correspondence between $t\le 2^{-1} n(n+1)$ and $(k,l)$ with $k\le l$.

Assume that $t_1\le t_2$ with $t_1=k_1+2^{-1} l_1(l_1-1)$ and $t_2=k_2+2^{-1}l_2(l_2-1)$. Then using the similar arguments we can show that $l_1\le l_2$ and further $k_1\le k_2$ when $l_1=l_2$. This means that for $t=k+2^{-1}l(l-1)$ with $1\le k\le l\le n$, we have
 $\mathcal{F}_{t}=\sigma\{\mathfrak{w}_1,\cdots,\mathfrak{w}_t\}=\sigma\{w_{ij}: 1\le i\le j< l \text{ or } 1\le i\le k\le j=l\}$. With such a representation, we can see that the expression in (\ref{0524.2}) is in fact a sum of martingale differences with respect to the  $\sigma$-algebra $\mathcal{F}_{k+2^{-1}i(i-1)}$. This fact entails that for $1 \leq k\le i \leq n$,
 \[ \mathbb{E}\left[(w_{ki}-\mathbb{E}w_{ki})b_{ki}+(w^2_{ki}-\mathbb{E}w^2_{ki})c_{ki} | \mathcal{F}_{k+2^{-1}i(i-1)-1}\right]=0,
 \]
 where $b_{ki}=x_k\sum_{1\le l<k\le n}w_{il}y_l+y_k\sum_{1\le l<k\le n}w_{il}x_l+x_i\sum_{1\leq l<i\le n}w_{kl}y_l+y_i\sum_{1\le l<i\le n}w_{kl}x_l+(1 - \delta_{ki})\mathbb{E}w_{ii}(x_iy_k+x_ky_i)$ with $\delta_{ki} = 1$ when $k = i$ and 0 otherwise, and $c_{ki}=x_ky_k+x_iy_i$. The conditional variance is given by
\begin{align} \label{0524.3}
&\sum_{1\le k,i\le n,\, k<i}\mathbb{E}\left\{\left[w_{ki}b_{ki}+(w_{ki}^2-\sigma_{ki}^2)c_{ki}\right]^2|\mathcal{F}_{k+2^{-1}i(i-1)-1}\right\}\non
&\quad+\sum_{1\le k\le n}\mathbb{E}\left\{\left[(\omega_{kk}-\mathbb{E}\omega_{kk})b_{kk}+2(\omega_{kk}^2-\sigma_{kk}^2)c_{kk}\right]^2|\mathcal{F}_{2^{-1}k(k+1)-1}\right\}\non
&=\sum_{1\le k,i\le n,\, k\le i}\sigma^2_{ki}b_{ki}^2+2\sum_{1\le k,i\le n,\, k\le i}\gamma_{ki}b_{ki}c_{ki}+\sum_{1\le k,i\le n,\, k\le i}\kappa_{ki}c_{ki}^2,
\end{align}
where $\gamma_{ki}=\mathbb{E}w_{ki}^3$  
and $\kappa_{ki}=\mathbb{E}(w_{ki}^2-\sigma_{ki}^2)^2$ for $k\neq i$,  and $\gamma_{kk}=2(\mathbb{E}\omega_{kk}^3-\sigma_{kk}^2\mathbb{E}\omega_{kk})$ and  $\kappa_{kk}=4\mathbb{E}(\omega_{kk}^2-\sigma_{kk}^2)^2$.

The mean of the random variable in (\ref{0524.3}) can be calculated as
\begin{align}\label{0525.1}
s^2_{\sbx,\sby}&=\mathbb{E}(\ref{0524.3})=\sum_{1\leq k,i\leq n,\, k\le i}\Big[\kappa_{ki}(x_ky_k+x_iy_i)^2+\sigma^2_{ki}\sum_{1\leq l<k\leq n}\sigma_{il}^2(x_ky_l+y_kx_l)^2\nonumber \\
&\quad+\sigma^2_{ki}\sum_{1\leq l<i\leq n}\sigma_{kl}^2(x_iy_l+y_ix_l)^2\Big]+\sum_{1\leq k,i\leq n,\, k\le i}\sigma_{ki}^2(1 - \delta_{ki})\left[\mathbb{E}(w_{ii}+w_{kk})\right]^2\nonumber\\
&\quad \times (x_ky_i+x_iy_k)^2.
\end{align}
Moreover, the variance of the random variable in (\ref{0524.3}) is given by
\begin{align}\label{0929.2h}
\kappa_{\bx,\by} & =\text{var}(\ref{0524.3})=\sum_{1\leq k_1,i_1,k_2,i_2 \leq n,\, k_1\le i_1,\, k_2\le i_2}\mathbb{E}\Big\{\big[\sigma^2_{k_1i_1}(z_{k_1i_1}^2-\mathbb{E}z_{k_1i_1}^2)\nonumber\\
&\quad+2\gamma_{k_1i_1}(x_{k_1}y_{k_1}+x_iy_i)z_{k_1i_1}\big] \big[\sigma^2_{k_2i_2}(z_{k_2i_2}^2-\mathbb{E}z_{k_2i_2}^2)\nonumber \\
&\quad+2\gamma_{k_2i_2}(x_{k_2}y_{k_2}+x_iy_i)z_{k_2i_2}\big]\Big\},
\end{align}
where  $z_{ki} = \sum_{1\leq l<k\leq n}w_{il}(x_ky_l+y_kx_l)+\sum_{1\leq l<i\leq n}w_{kl}(x_iy_l+y_ix_l)+(1-\delta_{ki})\mathbb{E}w_{ii}(x_iy_k+x_ky_i)$.

Let us recall the classical martingale CLT; see, for example, Lemma 9.12 of \cite{BS06}. If a martingale difference sequence $(Y_t)$ with respect to a $\sigma$-algebra $\mathcal{F}_t$ satisfies the following conditions:
  \begin{itemize}\label{0929.1h}
  \item [a)] $\frac{\sum_{t=1}^T\mathbb{E}(Y^2_t|\mathcal{F}_{t-1})}{\sum_{t=1}^T\mathbb{E} Y^2_t}\toP 1$,
  \item [b)]  $\frac{\sum_{t=1}^T\mathbb{E}[Y^2_tI(|Y_t|/\sqrt{\sum_{t=1}^T\mathbb{E}Y^2_t}|\ge \ep)]}{\sum_{t=1}^T\mathbb{E} Y^2_t}\le \frac{\sum_{t=1}^T\mathbb{E}Y^4_t}{\ep^2(\sum_{t=1}^T\mathbb{E}Y^2_t)^2}\rightarrow 0$ for any $\ep>0$,
  \end{itemize}
  then we have  $\frac{\sum_{t=1}^TY_t}{\sqrt{\sum_{t=1}^T\mathbb{E} Y^2_t}} \toD N(0,1)$ as $T \rightarrow \infty$, where $I(\cdot)$ denotes the indicator function. It follows from the assumption of  $\kappa_{\sbx,\sby}^{1/4}\ll s_{\sbx,\sby}$ that
\[ \frac{(\ref{0524.3})}{\mathbb{E}(\ref{0524.3})}\toP 1, \]
which shows that condition a) above is satisfied.
Moreover, by the simple fact that for any fixed $i$, $\mathbb{E}w_{li}^2y^2_l\le 1$, and the assumptions that $s_{\sbx,\sby}\rightarrow \infty$ and $\|\bbx\|_\infty \|\bby\|_{\infty}\rightarrow 0$, we have
\begin{align*}
&\sum_{1\leq k,i\leq n,\, k<i}\mathbb{E}\Big\{w_{ki}\Big[x_k\sum_{1\leq l<k\leq n}w_{il}y_l+y_k\sum_{1\leq l<k\leq n}w_{il}x_l+x_i\sum_{1\leq l<i\leq n}w_{kl}y_l \\
&\quad +y_i\sum_{1\leq l<i\leq n}w_{kl}x_l+\mathbb{E}w_{ii}(x_iy_k+x_ky_i)\Big]\Big\}^4+\sum_{1 \leq k \leq n}\mathbb{E}\Big[2(\omega_{kk}-\mathbb{E}\omega_{kk}) \\
&\quad \times \Big(x_k\sum_{1\leq l<k\leq n}w_{kl}y_l+y_k\sum_{1\leq l<k\leq n}w_{kl}x_l\Big)\Big]^4 \\
&\quad+\sum_{1\leq k,i\leq n,\, k< i } \Big\{\mathbb{E}\big[(w_{ki}^2-\sigma_{ki}^2)(x_ky_k+x_iy_i)\big]^4+\mathbb{E}\big[(\omega_{kk}^2-\sigma_{ki}^2)(x_ky_k+x_iy_i)\big]^4 \Big\} \ll s_{\sbx,\sby}^4,
\end{align*}
which entails that condition b) above is also satisfied. Therefore, an application of the martingale CLT concludes the proof of Lemma \ref{0524-1h}.

\subsection{Further technical details on conditions of Lemma \ref{0524-1h}} \label{SecB.2.new}

Let us gain some further insights into the technical conditions in Lemma \ref{0524-1h}. Define $a_{kl}=x_ky_l+y_kx_l$ and note that $\kappa_{ij}=\mathbb{E}(w_{ij}^2-\sigma_{ij}^2)^2=\mathbb{E}w_{ij}^4-\sigma_{ij}^4$. By the assumption of $|w_{ij}|\le 1$, it is easy to see that $0\le \kappa_{ij}\le \mathbb{E}w_{ij}^4 \le \mathbb{E}w_{ij}^2=\sigma_{ij}^2$. Then we can show that the random variable in (\ref{0524.3}) subtracted by its mean $s_{\sbx,\sby}^2$ can be represented as
\begin{align}\label{0614.1}
	&(\ref{0524.3})-s_{\sbx,\sby}^2=\sum_{1\le k,i\le n,\,k\le i}\sigma_{ki}^2\Big[\sum_{1\le l<k\le n}(w_{il}^2-\sigma_{il}^2)a_{kl}^2+\sum_{1\le l<i\le n}(w_{kl}^2-\sigma_{kl}^2)a_{il}^2 \nonumber\\
	&\quad+\sum_{1\le l_1,l_2<k\leq n,\,l_1\neq l_2}w_{il_1}w_{il_2}a_{kl_1}a_{kl_2}+\sum_{1\le l_1,l_2<i\leq n,\,l_1\neq l_2}w_{kl_1}w_{kl_2}a_{il_1}a_{il_2}\Big]\nonumber\\
	&\quad+2\sum_{1\le k,i \leq n,\,k\le i,\, 1\leq l_1<k\leq n,\, 1\leq l_2<i\leq n}\sigma_{ki}^2w_{il_1}w_{kl_2}a_{kl_1}a_{il_2} \nonumber\\
	&\quad+2\sum_{1\le k,i\le n,\,k\le i}\left[\gamma_{ki}a_{kk}+\sigma^2_{ki}a_{ki}(1-\delta_{ki})\mathbb{E}w_{ii}\right]\Big(\sum_{1\le l<k\le n}w_{il}a_{kl}+\sum_{1\le l<i\le n}w_{kl}a_{il}\Big).
	\end{align}
By (\ref{0614.1}) and (\ref{0928.18h}), we have
\begin{align} \label{0614.2}
	\kappa_{\sbx,\sby}& =\mathbb{E}\left[(\ref{0524.3})-s_{\sbx,\sby}^2\right]^2\le C\Big\{\mathbb{E}\Big[\sum_{1\le k,i\le n,\,k\le i}\sigma_{ki}^2\sum_{1\le l<k\le n}(w_{il}^2-\sigma_{il}^2)a_{kl}^2\Big]^2\non
&\quad+\mathbb{E}\Big[\sum_{1\le k,i\le n,\,k\le i}\sigma_{ki}^2\sum_{1\le l<i\le n}(w_{kl}^2-\sigma_{kl}^2)a_{il}^2\Big]^2+\mathbb{E}\Big(\sum_{1\le k,i\le n,\,k\le i}\sigma_{ki}^2\sum_{1\le l_1,l_2<k\leq n,\,l_1\neq l_2}w_{il_1}w_{il_2}\nonumber\\
&\quad \times a_{kl_1}a_{kl_2}\Big)^2+\mathbb{E}\Big(\sum_{1\le k,i\le n,\,k\le i}\sigma_{ki}^2\sum_{1\le l_1,l_2<i\le n,\,l_1\neq l_2}w_{kl_1}w_{kl_2}a_{il_1}a_{il_2}\Big)^2\non
	&\quad+\mathbb{E}\Big(\sum_{1\le k,i\leq n,\, k\le i, \, 1\leq l_1<k\leq n,\,1\leq l_2<i\le n}\sigma_{ki}^2w_{il_1}w_{kl_2}a_{kl_1}a_{il_2}\Big)^2+\mathbb{E}\sum_{1\le k,i\le n,\,k\le i}\big[\gamma_{ki}a_{kk}\non
&\quad+\sigma^2_{ki}a_{ki}(1-\delta_{ki})\mathbb{E}w_{ii}\big]\Big(\sum_{1\le l<k\le n}w_{il}a_{kl}+\sum_{1\le l<i\le n}w_{kl}a_{il}\Big)\Big\}\non
& \le C\Big\{\Big(\sum_{1\le k,i\le n,\,k\le i}\sigma_{ki}^2\Big)^2\Big(\sum_{1\leq l<k\leq n}\kappa_{il}a_{kl}^4+\sum_{1\leq l<i\leq n}\kappa_{kl}a_{il}^4\Big)\nonumber\\
&\quad+\sum_{1\le k_1,k_2,l_1,l_2\leq n,\, l_1\neq l_2,\, l_1<k_1,\, l_2<k_2}\sigma_{k_1i}^2\sigma_{k_2i}^2\sigma_{il_1}^2\sigma_{il_2}^2a_{k_1l_1}a_{k_1l_2}a_{k_2l_1}a_{k_2l_2}\non
	&\quad+\sum_{1\le k,i_1,i_2,l_1,l_2 \leq n,\, l_1\neq l_2 < \min\{i_1,i_2\}}\sigma_{ki_1}^2\sigma_{ki_2}^2\sigma_{kl_1}^2\sigma_{kl_2}^2a_{i_1l_1}a_{i_1l_2}a_{i_2l_1}a_{i_2l_2}\nonumber\\
	&\quad+\sum_{1\le k,i,l_1,l_2\leq n,\, k<i,\,l_1<k,\,l_2<i}\sigma_{ki}^2\sigma_{il_1}^2\sigma_{kl_2}^2a_{kl_1}^2a_{il_2}^2\non
	&\quad+\sum_{1\le k,i\le n,\,k<i}\big\{\gamma^2_{ki}a^2_{kk}+\sigma^4_{ki}a^2_{ki}(1-\delta_{ki})[\mathbb{E}(w_{ii}+w_{kk})]^2\big\}\Big(\sum_{1\le l<k\le n}\sigma_{il}^2a_{kl}^2+\sum_{1\le l<i\le n}\sigma_{kl}^2a_{il}^2\Big)\Big\}\non
	&=O\left\{n\sigma^8_n\|\bbx\|_\infty^4\|\bby\|_{\infty}^4\right\},
	\end{align}
where $C$ is some positive constant.

Given $\|\bbx\|_\infty \|\bby\|_{\infty}\rightarrow 0$, it follows from (\ref{0525.1}) that
\begin{align} \label{0614.3}
s_{\sbx,\sby}^2&=\sum_{1\leq k,i\leq n,\,k<i}\Big[\kappa_{ki}(x_ky_k+x_iy_i)^2+\sum_{1\leq l<k\leq n}\sigma_{il}^2(x_ky_l+y_kx_l)^2 \nonumber\\
&\quad+\sum_{1\leq l<i\leq n}\sigma_{kl}^2(x_iy_l+y_ix_l)^2\Big]\nonumber\\
&\ge \sigma_{\min}^2\sum_{1\leq k,i\leq n,\,k<i}\Big[\sum_{1\leq l<k\leq n}(x_ky_l+y_kx_l)^2+\sum_{1\leq l<i\leq n}(x_iy_l+y_ix_l)^2\Big] \nonumber\\
&\ge c\sigma_{\min}^2n,
\end{align}
where $\sigma_{\min}^2$ is defined in Condition \ref{as7}. Then we can exploit the upper bound on $\kappa_{\sbx,\sby}$ in (\ref{0614.2}) and the lower bound on $s_{\sbx,\sby}^2$ in (\ref{0614.3}) to simplify the conditions of Lemma \ref{0524-1h}, which can be reduced to
\begin{equation} \label{neweq025}
\|\bbx\|_\infty \|\bby\|_{\infty}\rightarrow 0, \  \frac{ \alpha^4_n\|\bbx\|_\infty^2\|\bby\|_{\infty}^2}{n^{1/2}  \sigma_{\min}^2}\rightarrow 0, \text{ and } \sigma_{\min}^2 n\rightarrow \infty.
\end{equation}
Therefore, the conclusions of Lemma \ref{0524-1h} hold as long as condition (\ref{neweq025}) is satisfied.

\subsection{Proof of Lemma \ref{lem: define-t}} \label{SecB.3}

In view of the definition of the function $f_k(z)$ defined in (\ref{0515.3.1}), we have
\begin{align} \label{neweq026h}
f'_k(z)&=d_k\Big\{\mathcal{R}(\bbv_k,\bbv_k,z)-\mathcal{R}(\bbv_k,\bbV_{-k},z)\left[\bbD_{-k}^{-1}+\mathcal{R}(\bbV_{-k},\bbV_{-k},z)\right]^{-1} \nonumber\\
& \quad \times \mathcal{R}(\bbV_{-k},\bbv_k,z)\Big\}'.
\end{align}
For $z\in [a_k,b_k]$, it follows from Lemma \ref{1212-1h}, Condition \ref{as3}, and the definition of $\mathcal R$ in \eqref{0619.0} that
\begin{align}\label{1004.1h}
\nonumber &\left\|\mathcal{R}(\bbV_{-k},\bbV_{-k},z)+z^{-1}\bbI\right\|=\left\|-\sum_{l=2}^L z^{-(l+1)}\bbV_{-k}^T\mathbb{E}\bbW^l\bbV_{-k} \right\| \\
& \quad \leq \sum_{l=2}^Lz^{-(l+1)}\left\|\bbV_{-k}^T\mathbb{E}\bbW^l\bbV_{-k}\right\|
 =O(\alpha_n^2 |z|^{-3}).
\end{align}
 Without loss of generality, we assume that $k\neq 1$. For $l$ such that  $|d_l|>|d_k|$, by \eqref{1004.1h} the diagonal entry of  $\bbD_{-k}^{-1}+\mathcal{R}(\bbV_{-k},\bbV_{-k},z)$ corresponding to $d_l$ is given by
 \[ d_l^{-1}-z^{-1}+O(\alpha_n^2 |z|^{-3})=(z-d_l)/(zd_l)+O(\alpha_n^2 |z|^{-3}). \]
 By Condition \ref{as3}, there exists some positive constant $c$ such that $\max\{|a_k|,|b_k|\}\le (1-c)|d_l|$. It follows that $|(z-d_l)/(zd_l)|\ge c/|z|$ and thus $|(z-d_l)/(zd_l)+O(\alpha_n^2 |z|^{-3})|^{-1}=O(|z|)$.  For the remaining diagonal entry with $|d_l|<|d_k|$, there exists some positive constant $c_1$ such that $\min\{|a_k|,|b_k|\}\ge (1+c_1)|d_l|$ and similarly we have $|(z-d_l)/(zd_l)+O(\alpha_n^2 |z|^{-3})|^{-1}=O(|z|)$. Thus it follows from (\ref{1004.1h}) that the off diagonal entries of $\bbD_{-k}^{-1}+\mathcal{R}(\bbV_{-k},\bbV_{-k},z)$ are dominated by the diagonal ones, leading to
 \begin{equation} \label{1010.1h}
 \left\|\left[\bbD_{-k}^{-1}+\mathcal{R}(\bbV_{-k},\bbV_{-k},z)\right]^{-1}\right\|=O(|z|)
 \end{equation}
 for all $z\in [a_k,b_k]$.

Next an application of Lemma \ref{1212-1h} gives
\[
\mathcal{R}'(\bbv_k,\bbv_k,z)=\sum_{{l=0,\,l\neq 1}}^L\frac{l+1}{z^{l+2}}\bbv_k^T\mathbb{E}\bbW^l\bbv_k=\frac{1}{z^2}+O(\alpha_n^2 |z|^{-4}).
\]
By (\ref{neweq026h}) and Condition \ref{as3}, we have
\[
\left\{\mathcal{R}(\bbv_k,\bbV_{-k},z)\left[\bbD_{-k}^{-1}+\mathcal{R}(\bbV_{-k},\bbV_{-k},z)\right]^{-1}\mathcal{R}(\bbV_{-k},\bbv_k,z)\right\}'=O(\alpha_n^4 |z|^{-6})=o(\alpha_n^2 |z|^{-4}).
\]
Thus in view of (\ref{neweq026h}), it holds that
\begin{equation} \label{neweq026}
f'_k(z)=d_k z^{-2}\left[1+o(1)\right]
\end{equation}
for $z\in [a_k, b_k ]$. We can see from (\ref{neweq026}) that $f_k(z)$ is a  monotone function over $z\in [a_k, b_k ]$ when matrix size $n$ is large enough.

Now recall  that
\[ f_k(d_k)=1+d_k\left\{\mathcal{R}(\bbv_k,\bbv_k,d_k)-\mathcal{R}(\bbv_k,\bbV_{-k},d_k)\left[\bbD_{-k}^{-1}+\mathcal{R}(\bbV_{-k},\bbV_{-k},d_k)\right]^{-1}\mathcal{R}(\bbV_{-k},\bbv_k,d_k)\right\}.
\]
By Lemma \ref{1212-1h}, we have
\[ 1+d_k\mathcal{R}(\bbv_k,\bbv_k,d_k)=1-\sum_{{l=0,\,l\neq 1}}^L\frac{1}{d^l_k}\bbv_k^T\mathbb{E}\bbW^l\bbv_k=O( \alpha_n^2 d_k^{-2}) \]
and
\[
d_k\mathcal{R}(\bbv_k,\bbV_{-k},d_k)\left[\bbD_{-k}^{-1}+\mathcal{R}(\bbV_{-k},\bbV_{-k},d_k)\right]^{-1}\mathcal{R}(\bbV_{-k},\bbv_k,d_k)=O(\alpha_n^2 d_k^{-2}).
\]
Thus it holds that $f_k(d_k)=O( \alpha_n^2 d_k^{-2})=o(1)$. Noticing that the derivative $f_k'(z)=d_k z^{-2}\left[1+o(1)\right]\sim d_k z^{-2}\sim |d_k|^{-1}$ and by the mean value theorem, we have $f_k(a_k)\sim o(1)+ |d_k|^{-1}(a_k-d_k) $ and $f_k(b_k)\sim o(1)+ |d_k|^{-1}(b_k-d_k)$, where $\sim$ represents the asymptotic order. Therefore, we see that $f_k(a_k)f_k(b_k)<0$ and  consequently the equation $f_k(z)=0$ has a unique solution for $z\in [a_k, b_k]$, which solution satisfies that $t_k=d_k+o(d_k)$. This completes the proof of Lemma \ref{lem: define-t}.

\subsection{Proof of Lemma  \ref{0426-1}} \label{SecB.4}
The asymptotic bounds characterized in Lemma \ref{0426-1} play a key role in establishing the more general asymptotic theory in Theorems \ref{0511-1} and \ref{0518-1}. We first assume that all the diagonal entries of $\bW = (w_{ij})_{1\leq i,j \leq n}$ are zero, that is, $w_{ii}=0$. The general case of possibly $w_{ii} \neq 0$ will be dealt with later. The main idea of the proof is to calculate the moments by counting the number of nonzero terms involved in $\mathbb{E}(\bbx^T\bbW^l\bby-\mathbb{E}\bbx^T\bbW^l\bby)^2$, which is a frequently used idea in random matrix theory; see, for example, Chapter 2 of \cite{BS06}. An important difference is that bounding the order of $\mathbb{E}(\bbx^T\bbW^l\bby-\mathbb{E}\bbx^T\bbW^l\bby)^2$ by simply counting the number of nonzero terms inside is too rough for our setting since the variances of the entries of $\bbW$ can be very different from each other. Observe that the nonzero terms of the variance involve the product of $w_{ij}^m$ with $m\ge 2$. We thus collect all such terms with the same index $i$ but different  index $j$, which means that we will bound $\sum_{j=1}^n\mathbb{E}|w_{ij}|^m\le \alpha_n^2$ instead of using  $\mathbb{E}|w_{ij}|^m\le 1$. Then we can obtain a more accurate order since $\alpha_n^2$ can be much smaller than $n$ in general. Our technical arguments here provide useful refinements to the classical idea of counting the number of nonzero terms from the random matrix theory.

Let $\bx = (x_1, \cdots, x_n)^T$ and $\by = (y_1, \cdots, y_n)^T$ be two arbitrary $n$-dimensional unit vectors, and $l \geq 1$ an integer. Expanding $\mathbb{E}(\bbx^T\bbW^l\bby-\mathbb{E}\bbx^T\bbW^l\bby)^2$ yields
\begin{align} \label{1212.1h}
\nonumber & \mathbb{E} (\bbx^T\bbW^l\bby-\mathbb{E}\bbx^T\bbW^l\bby)^2\\
& =\sum_{1\le i_1,\cdots,i_{l+1},j_1,\cdots,j_{l+1}\le n,\atop
i_s\neq i_{s+1},\, j_s\neq j_{s+1},\, 1\le s\le l}\mathbb{E} \Big[\left(x_{i_1}w_{i_1i_2}w_{i_2i_3}\cdots w_{i_li_{l+1}}y_{i_{l+1}}-\mathbb{E}x_{i_1}w_{i_1i_2}w_{i_2i_3}\cdots w_{i_li_{l+1}}y_{i_{l+1}}\right) \nonumber \\
& \quad \times \left(x_{j_1}w_{j_1j_2}w_{j_2j_3}\cdots w_{j_lj_{l+1}}y_{j_{l+1}}-\mathbb{E}x_{j_1}w_{j_1j_2}w_{j_2j_3}\cdots w_{j_lj_{l+1}}y_{j_{l+1}}\right)\Big].
\end{align}
Let $\bbi=(i_1,\cdots,i_{l+1})$ and $\bbj=(j_1,\cdots,j_{l+1})$ be two vectors taking values in $\{1,\cdots,n\}^{l+1}$. For any given vector $\bbi$, we define a graph $\mathcal{G}_{\bbi}$  whose vertices represent distinct values of the components of $\bbi$. Vertices $i_s$ and $i_{s+1}$ of $\mathcal G_{\bbi}$ are connected by undirected edges for $1\le s\le l$. Similarly we can also define graph  $\mathcal{G}_{\bbj}$ corresponding to $\bbj$. It can be seen that $\mathcal{G}_{\bbi}$ is a connected graph, which means that there exists some path from $i_s$ to $i_{s'}$ for any $1\le s\neq s'\le n$. Thus for each product
\begin{align}\label{0928.3h}
&\mathbb{E} \Big[\left(x_{i_1}w_{i_1i_2}w_{i_2i_3}\cdots w_{i_li_{l+1}}y_{i_{l+1}}-\mathbb{E}x_{i_1}w_{i_1i_2}w_{i_2i_3}\cdots w_{i_li_{l+1}}y_{i_{l+1}}\right) \nonumber \\
&\quad\times \left(x_{j_1}w_{j_1j_2}w_{j_2j_3}\cdots w_{j_lj_{l+1}}y_{j_{l+1}}-\mathbb{E}x_{j_1}w_{j_1j_2}w_{j_2j_3}\cdots w_{j_lj_{l+1}}y_{j_{l+1}}\right)\Big],
 \end{align}
there exists a corresponding graph $\mathcal{G}_{\bbi}\cup\mathcal{G}_{\bbj}$. If $\mathcal{G}_{\bbi}\cup\mathcal{G}_{\bbj}$ is not a connected graph, then the corresponding expectation
\begin{align*}
& \mathbb{E} \Big[\left(x_{i_1}w_{i_1i_2}w_{i_2i_3}\cdots w_{i_li_{l+1}}y_{i_{l+1}}-\mathbb{E}x_{i_1}w_{i_1i_2}w_{i_2i_3}\cdots w_{i_li_{l+1}}y_{i_{l+1}}\right) \\
&\quad \times \left(x_{j_1}w_{j_1j_2}w_{j_2j_3}\cdots w_{j_lj_{l+1}}y_{j_{l+1}}-\mathbb{E}x_{j_1}w_{j_1j_2}w_{j_2j_3}\cdots w_{j_lj_{l+1}}y_{j_{l+1}}\right)\Big]=0.
\end{align*}
This shows that in order to calculate the order of $\mathbb{E} (\bbx^T\bbW^l\bby-\mathbb{E}\bbx^T\bbW^l\bby)^2$, it suffices to consider the scenario of connected graphs $\mathcal{G}_{\bbi}\cup\mathcal{G}_{\bbj}$.

To analyze the term in (\ref{0928.3h}), let us calculate how many distinct vertices are contained in the connected graph $\mathcal{G}_{\bbi}\cup\mathcal{G}_{\bbj}$. Since there are $2l$ edges in $\mathcal{G}_{\bbi}\cup \mathcal{G}_{\bbj}$ and $\mathbb{E}w_{ss'}=0$ for $s\neq s'$, in order to get a nonzero value of (\ref{0928.3h}) each edge in $\mathcal{G}_{\bbi}\cup \mathcal{G}_{\bbj}$ has at least one copy. Thus for each nonzero (\ref{0928.3h}), we have $l$ distinct edges in  $\mathcal{G}_{\bbi}\cup \mathcal{G}_{\bbj}$. Since graph  $\mathcal{G}_{\bbi}\cup \mathcal{G}_{\bbj}$ is connected, there are at most $l+1$ distinct vertices in  $\mathcal{G}_{\bbi}\cup \mathcal{G}_{\bbj}$. Denote by $\mathcal{S}$ the set of all such pairs $(\bbi, \bbj)$. Combining the above arguments, we can conclude that
\begin{align}\label{0928.4h}
(\ref{1212.1h})&=\sum_{(\bbi, \bbj)\in \mathcal S}\mathbb{E} \Big[\left(x_{i_1}w_{i_1i_2}w_{i_2i_3}\cdots w_{i_li_{l+1}}y_{i_{l+1}}
 -\mathbb{E}x_{i_1}w_{i_1i_2}w_{i_2i_3}\cdots w_{i_li_{l+1}}y_{i_{l+1}}\right)
\nonumber \\
& \quad \times \left(x_{j_1}w_{j_1j_2}w_{j_2j_3}\cdots w_{j_lj_{l+1}}y_{j_{l+1}}-\mathbb{E}x_{j_1}w_{j_1j_2}w_{j_2j_3}\cdots w_{j_lj_{l+1}}y_{j_{l+1}}\right)\Big].
\end{align}
For notational simplicity, we denote $j_1,\cdots,j_{l+1}$ by $i_{l+2},\cdots,i_{2l+2}$ and define $\widetilde \bbi = (i_1,\cdots,i_{l+1},\\j_1,\cdots,j_{l+1})=(i_1,\cdots,i_{2l+2}).$ We also denote $\mathcal{G}_{\bbi}\cup \mathcal{G}_{\bbj}$ by $\mathcal{F}_{\widetilde \bbi}$ which has at most $l+1$ distinct vertices and $l$ distinct edges, with each edge having at least two copies. Then it holds that
\begin{align}\label{0928.5h}
& |(\ref{0928.4h})|=\Big|\sum_{\mathcal{F}_{\widetilde \bbi},\, \widetilde\bbi \in \mathcal S}\mathbb{E} \big[(x_{i_1}w_{i_1i_2}w_{i_2i_3}\cdots w_{i_li_{l+1}}y_{i_{l+1}}-\mathbb{E}x_{i_1}w_{i_1i_2}w_{i_2i_3}\cdots w_{i_li_{l+1}}y_{i_{l+1}})\nonumber\\
& \ \ \times (x_{i_{l+2}}w_{i_{l+2}i_{l+3}}w_{i_{l+3}i_{l+4}}\cdots w_{i_{2l+1}i_{2l+2}}y_{i_{2l+2}}-\mathbb{E}x_{i_{l+2}}w_{i_{l+2}i_{l+3}}w_{i_{l+3}i_{l+4}}\cdots w_{i_{2l+1}i_{2l+2}}y_{i_{2l+2}})\big]\Big|\non
&\le\sum_{\mathcal{F}_{\widetilde \bbi},\, \widetilde\bbi \in \mathcal S}\mathbb{E}\big|x_{i_1}w_{i_1i_2}w_{i_2i_3}\cdots w_{i_li_{l+1}}y_{i_{l+1}}x_{i_{l+2}}w_{i_{l+2}i_{l+3}}w_{i_{l+3}i_{l+4}}\cdots w_{i_{2l+1}i_{2l+2}}y_{i_{2l+2}}\big|\non
&\ \ +\sum_{\mathcal{F}_{\widetilde \bbi},\, \widetilde\bbi \in \mathcal S}\mathbb{E}\big|x_{i_1}w_{i_1i_2}w_{i_2i_3}\cdots w_{i_li_{l+1}}y_{i_{l+1}}\big|\mathbb{E}\big|x_{i_{l+2}}w_{i_{l+2}i_{l+3}}w_{i_{l+3}i_{l+4}}\cdots w_{i_{2l+1}i_{2l+2}}y_{i_{2l+2}}\big|.
\end{align}

Observe that each expectation in (\ref{0928.5h}) involves the product of some independent random variables, and $x_{i_1}w_{i_1i_2}w_{i_2i_3}\cdots w_{i_li_{l+1}}y_{i_{l+1}}$ and $x_{i_{l+2}}w_{i_{l+2}i_{l+3}}w_{i_{l+3}i_{l+4}}\cdots w_{i_{2l+1}i_{2l+2}}y_{i_{2l+2}}$ may share some dependency through factors $w_{ab}^{m_1}$ and $w_{ab}^{m_2}$, respectively, for some $w_{ab}$ and nonnegative integers $m_1$ and $m_2$. Thus in light of the inequality
\[ \mathbb{E}|w_{ab}|^{m_1}\mathbb{E}|w_{ab}|^{m_2}\le \mathbb{E}|w_{ab}|^{m_1+m_2}, \]
we can further bound (\ref{0928.5h}) as
\begin{align}\label{0928.6h}
(\ref{0928.5h}) & \le 2\sum_{\mathcal{F}_{\widetilde \bbi},\, \widetilde\bbi \in \mathcal S}\mathbb{E}\big|x_{i_1}w_{i_1i_2}w_{i_2i_3}\cdots w_{i_li_{l+1}}y_{i_{l+1}}x_{i_{l+2}}w_{i_{l+2}i_{l+3}}w_{i_{l+3}i_{l+4}}\cdots \nonumber\\
&\quad \times w_{i_{2l+1}i_{2l+2}}y_{i_{2l+2}}\big|.
\end{align}

To facilitate our technical presentation, let us introduce some additional notation. Denote by $\psi(2l+2)$ the set of partitions of the edges $\{(i_1,i_2),(i_2,i_3),\cdots,(i_{2l+1},i_{2l+2})\}$ and $\psi_{\ge 2}(2l+2)$ the subset of $\psi(2l+2)$ whose blocks have size at least two. Let  $P(\widetilde \bbi)\in\psi_{\ge 2}(2l+2)$ be the partition of $\{(i_1,i_2),(i_2,i_3),\cdots,(i_{2l+1},i_{2l+2})\}$ that is associated with the equivalence relation $(i_{s_1},i_{s_1+1})\sim (i_{s_2},i_{s_2+1})$ which is defined as if and only if $(i_{s_1},i_{s_1+1})=(i_{s_2},i_{s_2+1})$ or $(i_{s_1},i_{s_1+1})=(i_{s_2+1},i_{s_2})$. Denote by $|P(\widetilde \bbi)|=m$ the number of groups in the partition $P(\widetilde \bbi)$ such that the edges are equivalent within each group. We further denote the distinct edges in the partition $P(\widetilde \bbi)$ as $(s_1,s_2), (s_3,s_4), \cdots, (s_{2m-1},s_{2m})$ and the corresponding counts in each group as  $r_1, \cdots, r_m$, and define $\widetilde \bbs=(s_1,s_2,\cdots,s_{2m})$. For the vertices, let $\phi(2m)$ be the set of partitions of $\{1,2,\cdots,2m\}$ and $Q(\widetilde \bbs)\in\phi(2m)$ the partition that is associated with the equivalence relation $a\sim b$ which is defined as if and only if $s_a=s_b$. Note that $s_{2j-1}\neq s_{2j}$ since the diagonal entries of $\bbW$ are assumed to be zero for the moment. Then we have
\begin{align}\label{0928.7h}
&\sum_{\mathcal{F}_{\widetilde \bbi},\, \widetilde\bbi \in \mathcal S}\mathbb{E}\big|x_{i_1}w_{i_1i_2}w_{i_2i_3}\cdots w_{i_li_{l+1}}y_{i_{l+1}}x_{i_{l+2}}w_{i_{l+2}i_{l+3}}w_{i_{l+3}i_{l+4}}\cdots w_{i_{2l+1}i_{2l+2}}y_{i_{2l+2}}\big|\non
&\le \sum_{1\le |P(\widetilde \bbi)|=m\le l\atop
P(\widetilde \bbi)\in \psi_{\ge 2}(2l+2)}\sum_{\widetilde \bbi \text{ with partition } P(\widetilde \bbi)\atop
r_1,\cdots,r_m\ge 2}\sum_{Q(\widetilde \bbs)\in\phi(2m)}\sum_{\widetilde \bbs \text{ with partition } Q(\widetilde \bbs)\atop 1\le s_1,\cdots,s_{2m}\le n}|x_{i_1}y_{i_{l+1}}x_{i_{l+2}}y_{i_{2l+2}}|\nonumber\\
&\quad \times \prod_{j=1}^m\mathbb{E}\big|w_{s_{2j-1}s_{2j}}|^{r_j}.
\end{align}

We denote by $\mathcal{F}_{\widetilde \bbs}$ the graph constructed by the edges of $\widetilde \bbs$. Since the edges in $\widetilde \bbs$ are the same as those of the graph $\mathcal{F}_{\widetilde \bbi}$, we see that $\mathcal{F}_{\widetilde \bbs}$ is also a connected graph. In view of  (\ref{0928.7h}), putting term $|x_{i_1}y_{i_{l+1}}x_{i_{l+2}}y_{i_{2l+2}}|$ aside we need to analyze the summation
\[ \sum_{\widetilde \bbs \text{ with partition } Q(\widetilde \bbs)\atop 1\le s_1,\cdots, s_{2m}\le n}\prod_{j=1}^m\mathbb{E}\big|w_{s_{2j-1}s_{2j}}|^{r_j}. \]
If index $s_{2k-1}$ satisfies that $s_{2k-1}\neq s$ for all $s\in \{s_1,\cdots,s_{2m}\}\setminus \{s_{2k-1}\}$, that is, index $s_{2k-1}$ appears only in one $w_{s_{2j-1}s_{2j}}$, we call $s_{2k-1}$ a single index (or single vertex). If there exists some single index $s_{2k-1}$, then we have
 \begin{align}\label{0928.10h}
&\sum_{\widetilde \bbs \text{ with partition } Q(\widetilde \bbs)\atop 1\le s_1,\cdots,s_{2m}\le n}\prod_{j=1}^m\mathbb{E}\big|w_{s_{2j-1}s_{2j}}|^{r_j}  \nonumber \\
 &\quad \le \sum_{\tiny\substack{\widetilde \bbs\setminus \{s_{2k-1}\} \text{ with partition } Q(\widetilde \bbs\setminus \{s_{2k-1}\})\\ 1\le s_1,\cdots,s_{2k-2},s_{2k+2},s_{2m}\le n\\
   s_{2k}=s_j \text{ for some } 1\le j\le 2m}} \prod_{j=1}^m\mathbb{E}\big|w_{s_{2j-1}s_{2j}}|^{r_j} \sum_{s_{2k-1}=1}^n\mathbb{E}\big|w_{s_{2k-1}s_{2k}}|^{r_k}.
   \end{align}
Note that since graph $\mathcal{F}_{\widetilde \bbs}$ is connected and  index $s_{2k-1}$ is single, there exists some $j$ such that $s_j=s_{2k}$, which means that in the summation $\sum_{s_{2k-1}=1}^n\mathbb{E}\big|w_{s_{2k-1}s_{2k}}|^{r_k}$ index $s_{2k}$ is fixed. It follows from the definition of $\alpha_n$, $|w_{ij}|\le 1$, and $r_k\ge 2$ that
\[ \sum_{s_{2k-1}=1}^n\mathbb{E}\big|w_{s_{2k-1}s_{2k}}|^{r_k}\le \alpha_n^2. \]

After taking the summation over index $s_{2k-1}$, we see that there is one less edge in $\mathcal{F}(\widetilde \bbs)$. That is, by taking the summation above we will have one additional $\alpha_n^2$ in the upper bound while removing one edge from graph $\mathcal{F}(\widetilde \bbs)$. For the single index $s_{2k}$, we also have the same bound.  If $s_{2k_1-1}$ is not a single index, without loss of generality we assume that $s_{2k_1-1}=s_{2k-1}$. Then this vertex $s_{2k-1}$ need to deal with carefully. By the assumption of $|w_{ij}|\le 1$, we have
\[ \mathbb{E}|w_{2k-1,2k}|^{r_k}|w_{2k_1-1,2k_1}|^{r_{k_1}}\le \mathbb{E}|w_{2k-1,2k}|^{r_k}+\mathbb{E}|w_{2k_1-1,2k_1}|^{r_{k_1}}.
\]
 Then it holds that
 \begin{align}\label{0928.9h}
&\sum_{\widetilde \bbs \text{ with partition } Q(\widetilde \bbs)\atop 1\le s_1,\cdots,s_{2m}\le n}\prod_{j=1}^m\mathbb{E}\big|w_{s_{2j-1}s_{2j}}|^{r_j}  \nonumber\\
 &\le \sum_{\widetilde \bbs\setminus (s_{2k-1},s_{2k_1-1}) \text{ with partition } Q(\widetilde \bbs\setminus(s_{2k-1},s_{2k_1-1}))\atop 1\le s_1,\cdots,s_{2m}\le n}\prod_{j=1,\,j\neq k}^m\mathbb{E}\big|w_{s_{2j-1}s_{2j}}|^{r_j} \nonumber\\
 &\quad+\sum_{\widetilde \bbs\setminus(s_{2k-1},s_{2k_1-1}) \text{ with partition } Q(\widetilde \bbs\setminus (s_{2k-1},s_{2k_1-1}))\atop 1\le s_1,\cdots,s_{2m}\le n}\prod_{j=1,\,j\neq k_1}^m\mathbb{E}\big|w_{s_{2j-1}s_{2j}}|^{r_j}
 \end{align}
Note that since $\mathcal{F}_{\widetilde \bbs}$ is a connected graph, if we delete either edge $(s_{2k-1},s_{2k})$ or edge $(s_{2k_1-1},s_{2k_1})$ from graph $\mathcal{F}_{\widetilde \bbs}$ the resulting graph is also connected. Then the two summations on the right hand side of (\ref{0928.9h}) can be reduced to the case in (\ref{0928.10h}) for the graph with edge $(s_{2k-1},s_{2k})$ or $(s_{2k_1-1},s_{2k_1})$ removed, since $s_{2k-1}$ or $s_{2k_1-1}$ is a single index in the subgraph. Similar to (\ref{0928.10h}), after taking the summation over index $s_{2k-1}$ or $s_{2k_1-1}$ there are two less edges in graph $\mathcal{F}_{\widetilde \bbs}$ and thus we now obtain $2\alpha_n^2$ in the upper bound.

For the general case when there are $m_1$ vertices belonging to the same group, without loss of generality we denote them by $w_{ab_1},\cdots,w_{ab_{m_1}}$. If for any $k$ graph $\mathcal{F}_{\widetilde \bbs}$ is still connected after deleting edges $(a,b_1),\cdots,(a,b_{k-1}),(a,b_{k+1}),\cdots,(a,b_{m_1})$, then we repeat the process in (\ref{0928.9h}) to obtain a new connected graph by deleting $k-1$ edges in $w_{ab_1},\cdots,w_{ab_{m_1}}$ and thus obtain $k\alpha_n^2$ in the upper bound. Motivated by the key observations above, we carry out an iterative process in calculating the upper bound as follows.
\begin{itemize}
	\item[(i)] If there exists some single index in $\widetilde \bbs$, using  (\ref{0928.10h}) we can calculate the summation over such an index and then delete the edge associated with this vertex in $\mathcal{F}_{\widetilde \bbs}$. The corresponding vertices associated with this edge are also deleted. For simplicity, we also denote the new graph as $\mathcal{F}_{\widetilde \bbs}$. In this step, we obtain $\alpha_n^2$ in the upper bound.
	
	\item[(ii)] Repeat (i) until there is no single index in graph $\mathcal{F}_{\widetilde \bbs}$.
	
	\item[(iii)] If there exists some index associated with $k$ edges such that graph $\mathcal{F}_{\widetilde \bbs}$ is still connected after deleting any $k-1$ edges. Without loss of generality, let us consider the case of $k=2$. Then we can apply (\ref{0928.10h}) to obtain $\alpha_n^2$ in the upper bound. Moreover, we delete $k$ edges associated with this vertex in $\mathcal{F}_{\widetilde \bbs}$.

  \item[(iv)] Repeat (iii) until there is no such index.

 	\item[(v)] If there still exists some single index, turn back to (i).  Otherwise stop the iteration.
\end{itemize}

Completing the graph modification process mentioned above, we can obtain a final graph $\bbQ$ that enjoys the following properties:
\begin{itemize}
	\item[i)] Each edge does not contain any single index;
	
	
	\item[ii)] Deleting any vertex makes the graph disconnected.
\end{itemize}
Let $\bbS_{\bQ}$ be the spanning tree of graph $\bbQ$, which is defined as the subgraph of $\bbQ$ with the minimum possible number of edges. Since $\bbS_{\bQ}$ is a subgraph of $\bQ$, it also satisfies property ii) above. Assume that $\bbS_{\bQ}$ contains $p$ edges. Then the number of vertices in $\bbS_{\bQ}$ is $p+1$. Denote by $q_1,\cdots,q_{p+1}$ the vertices of $\bbS_{\bQ}$ and $\text{deg}(q_i)$ the degree of vertex $q_i$. Then by the degree sum formula, we have $\sum_{i=1}^{p+1}\text{deg}(q_i)=2p$. As a result, the spanning tree has at least two vertices with degree one and thus there exists a subgraph of $\bbS_{\bQ}$ without either of the vertices that is connected. This will result in a contradiction with property ii) above unless the number of vertices in graph $\bbQ$ is exactly one. Since $l$ is a bounded constant, the numbers of partitions $P(\widetilde \bbi)$ and $Q(\widetilde \bbs)$ are also bounded. It follows that
\begin{equation}\label{0928.11h}
 (\ref{0928.7h}) \le Cd_{\sbx}^2d_{\sby}^2\sum_{\widetilde \bbs \text{ with partition } Q(\widetilde \bbs)\atop 1\le s_1,\cdots,s_{2m}\le n}\prod_{j=1}^m\mathbb{E}\big|w_{s_{2j-1}s_{2j}}|^{r_j},
\end{equation}
where $d_{\sbx}=\|\bbx\|_{\infty}$, $d_{\sby}=\|\bbx\|_{\infty}$, and $C$ is some positive constant determined by $l$. Combining these arguments above and noticing that there are at most $l$ distinct edges in graph $\mathcal{F}_{\tilde \bbs}$, we can obtain
\begin{align*}
(\ref{0928.11h}) & \le Cd_{\sbx}^2d_{\sby}^2\alpha_n^{2l-2}\sum_{1\le s_{2k_0-1}, s_{2k_0}\le n,\,  (s_{2k_0-1},s_{2k_0})=\bbQ}\mathbb{E}\big|w_{s_{2k_0-1}s_{2k_0}}|^{r_{k_0}}\\
&\le Cd_{\sbx}^2d_{\sby}^2\alpha_n^{2l}n.
\end{align*}
Therefore, we have established a simple upper bound of $O\{d_{\bx}d_{\by}\alpha_n^l n^{1/2}\}$.

In fact, we can improve the aforementioned upper bound to $O(\alpha_n^{l-1})$. Note that the process mentioned above did not utilize the condition that both $\bx$ and $\by$ are unit vectors, that is, $\|\bbx\|=\|\bby\|=1$. Since term $|x_{i_1}y_{i_{l+1}}x_{i_{l+2}}y_{i_{2l+2}}|$ is involved in (\ref{0928.7h}), we can analyze them together with random variables $w_{ij}$. There are four different cases to consider.

1). Two pairs of indices $i_1$, $i_{l+1}$, $i_{l+2}$, $i_{2l+2}$ in $\mathcal{F}_{\widetilde \bbi}$ are equal. Without loss of generality, let us assume that $i_1=i_{l+1} \neq i_{l+2}=i_{2l+2}$. Then it holds that  $|x_{i_1}y_{i_{l+1}}x_{i_{l+2}}y_{i_{2l+2}}|=|x_{i_1}y_{i_1}x_{i_{l+2}}y_{i_{l+2}}|\le 4^{-1}(x_{i_1}^2+y^2_{i_1})(x_{i_{l+2}}^2+y^2_{i_{l+2}})$.
Let us consider the bound for
\begin{equation}\label{0928.12h}
\sum_{\widetilde \bbs \text{ with partition } Q(\widetilde \bbs)\atop 1\le s_1,\cdots,s_{2m}\le n}x_{i_1}^2x_{i_{l+2}}^2\prod_{j=1}^m\mathbb{E}\big|w_{s_{2j-1}s_{2j}}|^{r_j}.
\end{equation}
We assume without loss of generality that $i_1=s_1$ and $i_{l+2}=s_2$ for this partition. Then the summation in (\ref{0928.12h}) becomes
\[ \sum_{\widetilde \bbs \text{ with partition } Q(\widetilde \bbs)\atop 1\le s_1,\cdots,s_{2m}\le n}x_{s_1}^2x_{s_2}^2\prod_{j=1}^m\mathbb{E}\big|w_{s_{2j-1}s_{2j}}|^{r_j}. \]
By repeating the iterative process (i)--(v) mentioned before, we can bound the summation for fixed $s_2$ and obtain an alternative upper bound
\[ \sum_{s_1=1}^nx_{s_1}^2\mathbb{E}\big|w_{s_1s_2}|^{r_j}\le \sum_{s_1=1}^nx_{s_1}^2=1 \]
since $\bx$ is a unit vector. Thus for this step of the iteration, we obtain $1$ instead of $\alpha_n^2$ in the upper bound.  Since the graph is always connected during the iteration process, there exists another vertex $b$ such that $w_{s_2b}$ is involved in (\ref{0928.12h}). For index $s_2$, we do not delete the edges containing $s_2$ in the graph during the iterative process (i)--(v). Then after the iteration stops, the final graph $\bQ$ satisfies properties i) and ii) defined earlier except for vertex $s_2$. Since there are at least two vertices with degree one in  $\bbS_{\bQ}$, we will also reach a contradiction unless the number of vertices in graph $\bbQ$ is exactly one. As a result, we can obtain the upper bound
\begin{align}\label{0928.13h}
 (\ref{0928.7h})\le C\alpha_n^{2l-4}\sum_{1\le s_{2},b\le n, \,  (s_2,b)=\bbQ}\mathbb{E}x_{s_2}^2\big|w_{s_2b}|^{r}\le C\alpha_n^{2l-2}
\end{align}
with $C$ some positive constant. Therefore, the improved bound of $O(\alpha_n^{l-1})$ is shown for this case.

2). Indices $i_1$, $i_{l+1}$, $i_{l+2}$, $i_{2l+2}$ in $\mathcal{F}_{\widetilde \bbi}$ are all distinct. Then by the triangle inequality, we have
$|x_{i_1}y_{i_{l+1}}x_{i_{l+2}}y_{i_{2l+2}}|\le 4^{-1}(x_{i_1}^2+x^2_{i_{l+2}})(y_{i_{l+2}}^2+y^2_{i_{2l+2}})$. Thus this case reduces to case 1 above.

3). Indices $i_1$, $i_{l+1}$, $i_{l+2}$, $i_{2l+2}$ in $\mathcal{F}_{\widetilde \bbi}$ are all equal. Then it holds that  $|x_{i_1}y_{i_{l+1}}x_{i_{l+2}}y_{i_{2l+2}}|=x_{i_1}^2y^2_{i_1}\le x_{i_1}^2$. We see that there are at most $[(2l+2-2)/2]=l$ distinct vertices in the chain $\prod_{s=1}^{2l-1}w_{i_si_{s+1}}$ and for this case there are at most $l-1$ distinct edges in $\mathcal{F}_{\widetilde \bbi}$, where $[\cdot]$ denotes the integer part of a number. Compared to case 1, the maximum number of edges in the graph becomes smaller. Therefore, for this case we have
\begin{align}\label{0928.14h}
 (\ref{0928.7h})\le C\alpha_n^{2l-4}\sum_{1\le s_{1},b\le n,\,  (s_1,b)=\bbQ}\mathbb{E}x_{s_1}^2\big|w_{s_1b}|^{r}\le C\alpha_n^{2l-2},
\end{align}
where $C$ is some positive constant and we have assumed that $i_1=s_1$ without loss of generality.

4). Three of the indices $i_1$, $i_{l+1}$, $i_{l+2}$, $i_{2l+2}$ in $\mathcal{F}_{\widetilde \bbi}$ are equal. For such a case, without loss of generality let us write $|x_{i_1}y_{i_{l+1}}x_{i_{l+2}}y_{i_{2l+2}}|=|x_{i_1}^2y_{i_1}y_{i_{2l+2}}|$.  Then there are at most $[(2l+2-1)/2]=l$ distinct vertices in the chain
$\prod_{s=1}^{2l-1}w_{i_si_{s+1}}$ and thus for this case there are at most $l-1$ distinct edges in $\mathcal{F}_{\widetilde \bbi}$. Therefore, this case reduces to case 3 above.

In addition, we can also improve the upper bound to $O(\min\{d_{\sbx}\alpha_n^l, d_{\sby}\alpha_n^l\})$. The technical arguments for this refinement are similar to those for the improvement to order $O(\alpha_n^{l-1})$ above. As an example, we can bound the components of $\bby$ by $d_{\sby} = \|\by\|_\infty$, which leads to
$|x_{i_1}y_{i_{l+1}}x_{i_{l+2}}y_{i_{2l+2}}|\le d_{\sby}^2(x^2_{i_1}+x^2_{i_{l+1}})/2$. Then the analysis becomes similar to that for case 3 above. The only difference is that the length of graph $\mathcal{F}_{\widetilde \bbi}$ is at most $l$ instead of $l-1$. Thus similar to (\ref{0928.14h}), for this case we have
\begin{align}\label{0928.15h.new}
 (\ref{0928.7h})\le Cd_{\sby}^2\alpha_n^{2l-2}\sum_{1\le s_{2},b\le n,\,  (s_2,b)=\bbQ}\mathbb{E}x_{s_1}^2\big|w_{s_1b}|^{r}\le Cd_{\sby}^2\alpha_n^{2l},
\end{align}
where $C$ is some positive constant and we have assumed that $i_1=s_1$ or $x_{i_{l+1}}=s_1$ without loss of generality. The other one can then be used to remove a factor of $\alpha_n$. Thus we can obtain the claimed upper bound $O(\min\{d_{\sbx}\alpha_n^l, d_{\sby}\alpha_n^l\})$. Therefore, combining the two aforementioned improved bounds yields the desired upper bound of $O_p(\min\{\alpha_n^{l-1}, d_{\sbx}\alpha_n^l,d_{\sby}\alpha_n^l\})$.

We finally return to the general case of possibly $w_{ii}\neq 0$. Let us rewrite $\bbW$ as $\bbW=\bbW_0+\bbW_1$ with $\bbW_1=\diag(w_{11},\cdots,w_{nn})$. Then it holds that
\[ \bbx^T\bbW^l\bby-\mathbb{E}\bbx^T\bbW^l\bby=\bbx^T(\bbW_0+\bbW_1)^l\bby-\mathbb{E}\bbx^T(\bbW_0+\bbW_1)^l\bby.
\]
Recall the classical inequality
\begin{equation}\label{0928.18h}
\mathbb{E}(X_1+\cdots+X_m)^2\le m(\mathbb{E}X_1^2+
\cdots+\mathbb{E}X_m^2),
  \end{equation}
  where $X_1,\cdots, X_m$ are $m$ random variables with finite second moments. Define a function
\begin{equation} \label{neweq.100}
f(\bbh)=\prod_{i=1}^l \bW_{h_i},
\end{equation}
where the vector $\bbh=(h_1,\cdots,h_l)$ with $h_i = 0$ or $1$. Then we have
\begin{align}\label{0928.20h-1}  \mathbb{E}&\left[\bbx^T(\bbW_0+\bbW_1)^l\bby-\mathbb{E}\bbx^T(\bbW_0+\bbW_1)^l\bby\right]^2=\mathbb{E}\Big\{\sum_{\bbh}\bbx^T[f(\bbh)-\mathbb{E}f(\bbh)]\bby\Big\}^2 \nonumber\\
&\le 2^{l}\sum_{\bbh}\mathbb{E}\Big\{\bbx^T[f(\bbh)-\mathbb{E}f(\bbh)]\bby\Big\}^2.
\end{align}
This shows that we need only to consider terms of form  $\mathbb{E}\{\bbx^T[f(\bbh)-\mathbb{E}f(\bbh)]\bby\}^2$, each of which is a polynomial of $\bbW_0$ and $\bbW_1$.

As an example, let us analyze the term  $\mathbb{E}(\bbx^T\bbW_1\bbW_0^{l-1}\bby-\mathbb{E}\bbx^T\bbW_1\bbW_0^{l-1}\bby)^2$. Similar to (\ref{1212.1h}), it can be shown that
 \begin{align} \label{0928.15h}
&\mathbb{E}(\bbx^T\bbW_1\bbW_0^{l-1}\bby-\mathbb{E}\bbx^T\bbW_1\bbW_0^{l-1}\bby)^2\non
&=\sum_{1\le i_1,\cdots,i_{l},j_1,\cdots,j_{l}\le n,\atop
i_s\neq i_{s+1},\,  j_s\neq j_{s+1},\, 1\le s\le l}\mathbb{E} \Big[\left(x_{i_1}w_{i_1i_1}w_{i_1i_2}\cdots w_{i_{l-1}i_{l}}y_{i_{l}}-\mathbb{E}x_{i_1}w_{i_1i_1}w_{i_1i_2}\cdots w_{i_{l-1}i_{l}}y_{i_{l}}\right) \nonumber \\
&\quad \times \left(x_{j_1}w_{j_1j_1}w_{j_1j_2}\cdots w_{j_{l-1}j_{l}}y_{j_{l}}-\mathbb{E}x_{j_1}w_{j_1j_1}w_{j_1j_2}\cdots w_{j_{l-1}j_{l}}y_{j_{l}}\right)\Big].
\end{align}
Repeating the arguments from (\ref{1212.1h})--(\ref{0928.6h}), we can obtain
\begin{align*}
(\ref{0928.15h}) & \le 2\sum_{\mathcal{F}_{\widetilde \bbi}}\mathbb{E}\big|x_{i_1}w_{i_1i_1}w_{i_1i_2}\cdots w_{i_{l-1}i_{l}}y_{i_{l}}x_{i_{l+1}}w_{i_{l+1}i_{l+1}}w_{i_{l+1}i_{l+2}}\cdots w_{i_{2l-1}i_{2l}}y_{i_{2l}}\big| \\
&\le 2\sum_{\mathcal{F}_{\widetilde \bbi}}\mathbb{E}\big|x_{i_1}w_{i_1i_2}\cdots w_{i_{l-1}i_{l}}y_{i_{l}}x_{i_{l+1}}w_{i_{l+1}i_{l+2}}\cdots w_{i_{2l-1}i_{2l}}y_{i_{2l}}\big|.
\end{align*}
Comparing to (\ref{0928.6h}), we can see that by replacing the diagonal entries with $1$ in the expectations, the number of edges in this graph is no more than the original one in (\ref{0928.6h}). Thus repeating all the steps before (\ref{0928.15h}), we can deduce the bound
\[
\mathbb{E}(\bbx^T\bbW_1\bbW_0^{l-1}\bby-\mathbb{E}\bbx^T\bbW_1\bbW_0^{l-1}\bby)^2=O(\min\{\alpha_n^{2(l-1)},d_{\bx}^2\alpha_n^{2l}, d_{\by}^2\alpha_n^{2l}\}).
\]
For the other expectations $\mathbb{E}\{\bbx^T[f(\bbh)-\mathbb{E}f(\bbh)]\bby\}^2$, by the same reason that $\bbW_1$ is a diagonal matrix we can obtain a similar expression as (\ref{0928.15h}) with the number of edges no larger than the original one for $\mathbb{E}(\bbx^T\bbW_0^{l}\bby-\mathbb{E}\bbx^T\bbW_0^{l}\bby)^2$. Thus all the technical arguments above can be applied to $\mathbb{E}\{\bbx^T[f(\bbh)-\mathbb{E}f(\bbh)]\bby\}^2$ so we can have the same order for the upper bound as before. This shows that all the previous arguments can indeed be extended to the general case of possibly $w_{ii}\neq 0$, which concludes the proof of Lemma  \ref{0426-1}.

\subsection{Proof  of Lemma \ref{1212-1h}} \label{SecB.6}
The main idea of the proof is similar to that for the proof of Lemma \ref{0426-1} in Section \ref{SecB.4}. We first consider the case when all the diagonal entries of $\bW = (w_{ij})_{1\leq i,j \leq n}$ are zero, that is, $w_{ii}=0$. Then we can derive a similar expression as (\ref{1212.1h})
\begin{align}\label{0928.19h}
\mathbb{E}\bbx^T\bbW^l\bby&=\sum_{1\le i_1,\cdots,i_{l+1}\le n\atop i_s\neq i_{s+1},\, 1\le s\le l}\mathbb{E} \left(x_{i_1}w_{i_1i_2}w_{i_2i_3}\cdots w_{i_li_{l+1}}y_{i_{l+1}}\right).
\end{align}
By the definition of graph $\mathcal{G}_{\bbi}$ in the proof of Lemma \ref{0426-1}, we can obtain a similar expression as (\ref{0928.6h})
\begin{align} \label{1108.2h}
|(\ref{0928.19h})| & \le \sum_{\mathcal{G}_{\bbi} \text{ with at most $[l/2]$ distinct edges and $[l/2]+1$ distinct vertices}}\mathbb{E}\big|x_{i_1}w_{i_1i_2}w_{i_2i_3}\cdots \nonumber\\
&\quad \times w_{i_li_{l+1}}y_{i_{l+1}}\big|.
\end{align}
Using similar arguments for bounding the order of the summation through the iterative process as for case 3 in the proof of Lemma \ref{0426-1} and noticing that $|x_{i_1}y_{i_{l+1}}|\le 2^{-1} (x_{i_1}^2+y_{i_{l+1}}^2)$, we can deduce the desired bound
\begin{equation}\label{0928.21h}
\mathbb{E}\bbx^T\bbW^l\bby=O(\alpha_n^{l-1}),
\end{equation}
where the diagonal entries of $\bbW$ have been assumed to be zero.

For the general case of $\bbW$ with possibly nonzero diagonal entries, we can apply the similar expansion as in the proof of Lemma \ref{0426-1} to get
\begin{equation}\label{0928.20h}
\mathbb{E}\bbx^T(\bbW_0+\bbW_1)^l\bby=\sum_{\bbh}\mathbb{E}\bbx^Tf(\bbh)\bby,
\end{equation}
where $\bbW=\bbW_0+\bbW_1$ with $\bbW_1=\diag(w_{11},\cdots,w_{nn})$, and vector $\bbh$ and function $f(\bbh)$ are as defined in (\ref{neweq.100}). Since by assumption $\bbW_1$ is a diagonal matrix with bounded entries, an application of similar arguments as in the proof of Lemma \ref{0426-1} gives
	\[
	\mathbb{E}\bbx^Tf(\bbh)\bby=O(\alpha_n^{l-1}).
	\]
	To see this, with similar arguments as below (\ref{0928.20h-1}) let us analyze the term $\mathbb{E}\bbx^T\bbW_1\bbW_0^{l-1}\bby$ as an example. Similar to (\ref{0928.19h}), it holds that
	\begin{align}\label{1108.1h}
	\mathbb{E}\bbx^T\bbW_1\bbW_0^{l-1}\bby=\sum_{1\le i_1,\cdots,i_{l}\le n\atop i_s\neq i_{s+1},\, 1\le s\le l-1}\mathbb{E} \left(x_{i_1}w_{i_1i_1}w_{i_1i_2}\cdots w_{i_{l-1}i_l}y_{i_l}\right).
	\end{align}
	By the assumption of $\max_{1\le i\le n}|w_{ii}|\le 1$, we can derive a similar bound as (\ref{1108.2h})
	\begin{align}\label{1108.3h}
	|(\ref{1108.1h})| & \le \sum_{\mathcal{G}_{\bbi} \text{ with at most $[(l-1)/2]$ distinct edges and $[(l-1)/2]+1$ distinct vertices}}\mathbb{E}\big|x_{i_1}w_{i_1i_2}w_{i_2i_3}\cdots \nonumber\\
	&\quad \times w_{i_{l-1}i_l}y_{i_l}\big|.
	\end{align}
	Since the number of edges is no more than that in (\ref{1108.2h}), we can obtain the same bound
	\[
	\mathbb{E}\bbx^T\bbW_1\bbW_0^{l-1}\bby=O(\alpha_n^{l-1}).
	\]
For the other terms in (\ref{0928.20h}), by the same reason that $\bbW_1$ is a diagonal matrix with bounded entries we can derive similar expression as (\ref{1108.3h}) with the number of edges no more than that in (\ref{1108.2h}). Therefore, since $l$ is a bounded constant we can show that $\mathbb{E}\bbx^T\bbW^l\bby=O(\alpha_n^{l-1})$ for the general case of $\bbW$ with possibly nonzero diagonal entries. This completes the proof of Lemma \ref{1212-1h}.

\subsection{Lemma \ref{0505-1} and its proof} \label{SecB.5}

\begin{lem} \label{0505-1}
	The random matrix $\bW$ given in (\ref{model}) satisfies that for any positive constant $L$, there exist some positive constants $C_L$ and $\sigma$ such that
	\begin{equation} \label{neweq028}
	\mathbb{P}\left\{\|\bbW\|\ge C_L (\log n)^{1/2}\alpha_n\right\} \le n^{-L},
	\end{equation}
	where $\|\cdot\|$ denotes the matrix spectral norm and  $\alpha_n=\|\mathbb{E}(\bbW-\mathbb{E}\bbW)^2\|^{1/2}$.
\end{lem}

\noindent \textit{Proof}. The conclusion of Lemma \ref{0505-1} follows directly from Theorem 6.2 of \cite{T12}.

\section{Further technical details on when asymptotic normality holds for Theorem \ref{0518-1} } \label{SecC}

We now consider the joint distribution of the three random variables specified in expression (\ref{0605.5}) in the proof of Theorem \ref{0518-1} in Section \ref{SecA.6}. To establish the joint asymptotic normality under some regularity conditions, it suffices to show that the random vector $(\tr[(\bbW-\mathbb{E}\bbW)\bbJ_{\sbx,\sby,k,t_k}-(\bbW^2-\mathbb{E}\bbW^2)\bbL_{\sbx,\sby,k,t_k}],\tr((\bbW-\mathbb{E}\bbW)\bbv_k\bbv_k^T), \tr((\bbW-\mathbb{E}\bbW)\bbQ_{\sbx,\sby,k,t_k}))$ tends to some multivariate normal distribution as matrix size $n$ increases, where we consider the de-meaned version of this random vector for simplicity. Consequently, we need to show that for any constants $c_1$, $c_2$, and $c_3$ such that $c_1^2+c_2^2+c_3^2=1$, the linear combination
\begin{align}\label{0625.2}
& c_1\tr[(\bbW-\mathbb{E}\bbW)\bbJ_{\sbx,\sby,k,t_k}-(\bbW^2-\mathbb{E}\bbW^2)\bbL_{\sbx,\sby,k,t_k}]+c_2\tr((\bbW-\mathbb{E}\bbW)\bbv_k\bbv_k^T) \nonumber\\
&\quad +c_3\tr((\bbW-\mathbb{E}\bbW)\bbQ_{\sbx,\sby,k,t_k})
\end{align}
converges to a normal distribution asymptotically. Define $\bbS=\bbv_k\bbv_k^T$ and let $\bbJ$, $\bbL$, and $\bbQ$ be the rescaled versions of $\bbJ_{\sbx,\sby,k,t_k}$, $\bbL_{\sbx,\sby,k,t_k}$, and $\bbQ_{\sbx,\sby,k,t_k}$, respectively, such that the asymptotic variance of each of the above three terms is equal to one. Then it remains to analyze the asymptotic behavior of the random variable
\begin{align}\label{0606.6}
&\sum_{1\leq k,i\leq n,\, k\le i}w_{ki}\Big\{c_1\Big[\sum_{1\leq l<k\leq n}w_{il}\bbL_{kl}+\sum_{1\leq l<i\leq n}w_{kl}\bbL_{il}+\bbJ_{ki}+(1 - \delta_{ki})(\bbL_{ki}+\bbL_{ik})\mathbb{E}w_{ii}\Big]\nonumber\\
&\quad + (1-\delta_{ki})(c_2\bbS_{ki}+ c_3\bbQ_{ki})\Big\}+c_1\sum_{1\leq k,i\leq n,\, k\le i}(w_{ki}^2-\sigma_{ki}^2)(\bbL_{kk}+\bbL_{ii}),
\end{align}
where $\bA_{ij}$ indicates the $(i,j)$th entry of a matrix $\bA$ and $\delta_{ki} = 1$ when $k = i$ and $0$ otherwise.

Using similar arguments as in (\ref{0524.3}), we can show that (\ref{0606.6}) is in fact a sum of martingale differences with respect to the $\sigma$-algebra $\mathcal{F}_{k+2^{-1} i(i-1)-1}$. The conditional variance of the random variable given in (\ref{0606.6})  can be calculated as
\begin{align}\label{0605.7}
&\sum_{1\le k,i\le n,\, k\le i}\sigma^2_{ki}\Big\{c_1\Big[\sum_{1\le l<k\le n}w_{il}\bbL_{kl}+\sum_{1\le l<i\le n}w_{kl}\bbL_{il}+\bbJ_{ki}+(1-\delta_{ki})(\bbL_{ki}+\bbL_{ik})\mathbb{E}w_{ii}\Big]\nonumber\\
&\quad+ (1-\delta_{ki})(c_2\bbS_{ki}+ c_3\bbQ_{ki})\Big\}^2 +c_1^2\sum_{1\le k,i\le n,\, k\le i}\kappa_{ki}(\bbL_{kk}+\bbL_{ii})^2\nonumber\\
&\quad +2c_1\sum_{1\leq k,i\leq n,\, k\le i}\gamma_{ki}(\bbL_{kk}+\bbL_{ii}) \Big\{c_1\Big[\sum_{1\le l<k\le n}w_{il}\bbL_{kl}+\sum_{1\le l<i\le n}w_{kl}\bbL_{il}+\bbJ_{ki} \nonumber\\
&\quad +(1-\delta_{ki})(\bbL_{ki}+\bbL_{ik})\mathbb{E}w_{ii}\Big] + (1-\delta_{ki})(c_2\bbS_{ki}+ c_3\bbQ_{ki})\Big\}.
\end{align}
Moreover, the expectation of the random variable given in (\ref{0605.7}) can be shown to take the form
\begin{align}\label{0605.8}	
&c_1^2\sum_{1\le k,i\le n,\,k\le i}\Big\{\sigma^2_{ki}\Big[\sum_{1\le l<k\le n}\sigma^2_{il}\bbL^2_{kl}+\sum_{1\le l<i\le n}\sigma^2_{kl}\bbL^2_{il}+\bbJ^2_{ki}+(1-\delta_{ki})(\bbL_{ki}+\bbL_{ik})^2(\mathbb{E}w_{ii})^2\Big]\non
&\quad+\kappa_{ki}(\bbL_{kk}+\bbL_{ii})^2\Big\}+ c^2_2\sum_{1\le k,i\le n,\,k\leq i}\sigma^2_{ki}\bbS^2_{ki}+ c^2_3\sum_{1\le k,i\le n,\,k\leq i}\sigma^2_{ki}\bbQ^2_{ki}\nonumber\\
&\quad+2\sum_{1\le k,i\le n,\,k\leq i}\Big[\sigma^2_{ki}(c_2\bbS_{ki}+ c_3\bbQ_{ki})(\bbL_{ki}+\bbL_{ik})\mathbb{E}w_{ii}\Big]\non
&\quad+2c_1 c_2\sum_{1\le k,i\le n,\,k\leq i}\gamma_{ki}\bbS_{ki}(\bbL_{kk}+\bbL_{ii})+2c_1 c_3\sum_{1\le k,i\le n,\,k\leq i}\gamma_{ki}\bbQ_{ki}(\bbL_{kk}+\bbL_{ii})\non
&\quad+2c_2c_3\sum_{1\le k,i\le n,\,k\leq i}\sigma^2_{ki}\bbS_{ki}\bbQ_{ki}+2c_1^2\sum_{1\le k,i\le n,\,k\leq i}\Big[\kappa_{ki}(\bbL_{kk}+\bbL_{ii})(\bbL_{ki}+\bbL_{ik})\mathbb{E}w_{ii}\Big].
\end{align}

Let us consider the following three regularity conditions.
\begin{itemize}
	\item[i)] Assume that the six individual summation terms in (\ref{0605.8}) tend to some constants asymptotically. Then (\ref{0605.8}) tends to some constant $C$ asymptotically. Without loss of generality, we assume that $C \neq 0$; otherwise (\ref{0606.6}) tends to zero in probability.
	
	\item[ii)] Assume that $\text{SD}(\ref{0605.7})\ll (\ref{0605.8})$, where  $\text{SD}$ stands for the standard deviation of a random variable.
	
	\item[iii)] Assume that
\begin{align} \label{0615.1}
& \sum_{1\leq k,i\leq n,\,k<i}\kappa_{ki}\Big\{\mathbb{E}\Big[\sum_{1\leq l<k\leq n}w_{il}\bbL_{kl}+\sum_{1\leq l<i\leq n}w_{kl}\bbL_{il}+\bbJ_{ki}+(1 - \delta_{ki})(\bbL_{ki}+\bbL_{ik})\mathbb{E}w_{ii}\Big]^4 \nonumber\\
&\quad + (1 - \delta_{ki})(\bbS^4_{ki}+ \bbQ^4_{ki})\Big\} +\sum_{1\leq k,i\leq n,\,k<i}\mathbb{E}(w_{ki}^2-\sigma_{ki}^2)^4(\bbL_{kk}+\bbL_{ii})^4\ll 1.
\end{align}
\end{itemize}
We can see that conditions i) and ii) entail condition a) in the proof of Lemma \ref{0524-1h} in Section \ref{SecB.2} below (\ref{0929.2h}), while condition iii) entails condition b). Therefore, (\ref{0606.6}) converges to a normal distribution asymptotically.

\renewcommand{\thesubsection}{D.\arabic{subsection}}
\section{Relaxing the spike strength condition and proof sketch for results in Section \ref{sec: app2}} \label{sec:cond2ii}

The main goal of this section is to show that all the results continue to hold when Condition  \ref{as3}i)  is replaced with Condition \ref{as3}ii), which is a weaker assumption on the spike strength. Thus from now on, we will assume Condition \ref{as3}ii) instead of Condition \ref{as3}i). Moreover, we provide the proof sketch for results in Section \ref{sec: app2}.

\subsection{Replacing  Condition \ref{as3}i) with Condition \ref{as3}ii)} \label{SecD.1}

Checking the proofs of our theorems, we can see that it is sufficient to show that the asymptotic expansion of $\bbx^T\bbG(z)\bby$ remains to hold under Condition \ref{as3}ii). In other words,  we need to prove \eqref{0512.1} and \eqref{0928.2h} under Condition \ref{as3}ii). To accommodate the  smaller magnitude of $d_K$ in Condition \ref{as3}ii), the key  idea is to carefully examine the asymptotic expansions \eqref{0512.1} and \eqref{0928.2h} as $L\rightarrow \infty$.
 To this end, we choose $L=\log n$ and  
 	define $c'=c/(1+2^{-1}c_0)$.
  Since $\alpha_n\le n^{1/2}$, we  have the following improved version of inequality \eqref{eq011}
\begin{equation} \label{eq011h}
\frac{\alpha_n^{L+1}(C\log n)^{(L+1)/2}}{\min\{|a_K|,|b_K|\}^{L-2}}\le\frac{\alpha_n^3(C\log n)^{(L+1)/2}}{(c'\log n)^{L-2}} \le \frac{C^{(\log n+1)/2}n^{3/2}}{(\log n)^{(\log n-5)/2}c'^{\log n-2}}\rightarrow 0
\end{equation}
for any positive constant $C$.
	
We first show that \eqref{0512.1} holds with the choice of $L=\log n$. In view of \eqref{0508.2}, it is sufficient to establish the following two equations
\begin{equation}\label{r1}
\sum_{l=L+1}^{\infty}z^{-(2l+2)}\bbx^T\bbW^l\bby=O_p(\frac{1}{|z|^4})
\end{equation}
and
\begin{equation}\label{r2}
\sum_{l=2}^{L}z^{-(2l+2)}\bbx^T(\bbW^l-\mathbb{E}\bbW^l)\bby=O_p(\frac{\alpha_n}{|z|^3})
\end{equation}
for $z\in \Omega_k$. In fact, \eqref{r1} is a direct consequence of Lemma \ref{0505-1} and  \eqref{eq011h}.  
	In light of the definitions of $a_k$ and $b_k$ below \eqref{0515.3.1}, we can conclude that for any $z\in \Omega_k$, $|z|> 4c_1\alpha_n\log n$.
Thus we see that
\begin{equation}\label{r0}
\left\{\frac{\alpha_n^{2l}(4c_1\log n)^{2l}}{|z|^{2l}}\right\} \ \text{is a decreasing  sequence when $l$ is increasing for } z \in \Omega_k.
\end{equation}
Then it follows from Lemma \ref{0426-1h} and \eqref{r0} that
\begin{eqnarray}\label{r3}
\sum_{l=2}^{9}z^{-(2l+2)}\bbx^T(\bbW^l-\mathbb{E}\bbW^l)\bby=O_p(\frac{\alpha_n}{|z|^3}),
\end{eqnarray}
\begin{eqnarray}\label{r4}
&&\left|\mathbb{E}\left[\sum_{l=\sqrt L}^{L}z^{-(2l+2)}\bbx^T(\bbW^l-\mathbb{E}\bbW^l)\bby\right]^2\right|\le L\sum_{l=\sqrt L}^{L}|z|^{-(2l+2)}\mathbb{E}\left[\bbx^T(\bbW^l-\mathbb{E}\bbW^l)\bby\right]^2\non
&&\le CL\sum_{l=\sqrt L}^{L}\frac{(4c_1l)^{2l}\alpha_n^{2l-2}}{|z|^{2l+2}}
\le CL\sum_{l=\sqrt L}^{L}\frac{(4c_1\log n)^{2l}\alpha_n^{2l-2}}{|z|^{2l+2}}\non
&&\le C(\log n)^2\frac{\alpha_n^4(4c_1\log n)^{2\sqrt{\log n}}}{|z|^8(c'\log n)^{2\sqrt{\log n}-4}}\le \frac{C(4c_1)^6\alpha_n^4(\log n)^6}{|z|^8(c'/(4c_1))^{2\sqrt{\log n}-4}}\ll \frac{\alpha_n^4}{|z|^8},
\end{eqnarray}
and
\begin{eqnarray}\label{r5}
&&\left|\mathbb{E}\left[\sum_{l=10}^{\sqrt L}z^{-(2l+2)}\bbx^T(\bbW^l-\mathbb{E}\bbW^l)\bby\right]^2\right|\le \sqrt L\sum_{l=10}^{\sqrt L}|z|^{-(2l+2)}\mathbb{E}\left[\bbx^T(\bbW^l-\mathbb{E}\bbW^l)\bby\right]^2\non
&&\le C\sqrt L\sum_{l=10}^{\sqrt L}\frac{(4c_1l)^{2l}\alpha_n^{2l-2}}{|z|^{2l+2}}\le C\sqrt L\sum_{l=10}^{\sqrt L}\frac{(4c_1\sqrt{\log n})^{2l}\alpha_n^{2l-2}}{|z|^{2l+2}}\non
&&\le C\log n\frac{\alpha_n^4((4c_1)^2\log n)^{10}}{|z|^8(c'\log n)^{20-4}}\ll \frac{\alpha_n^4}{|z|^8}.
\end{eqnarray}
Therefore, combining \eqref{r3}--\eqref{r5} yields \eqref{r2}.

To establish \eqref{0928.2h} with $L=\log n$, we need only to prove \eqref{r1} and
\begin{equation}\label{r6}
\sum_{l=3}^{L}z^{-(2l+2)}\bbx^T(\bbW^l-\mathbb{E}\bbW^l)\bby=O_p(\frac{\alpha^2_n}{|z|^4}),
\end{equation}
where the former has been shown before.  By Lemma \ref{0426-1h}, we can deduce
\begin{eqnarray}\label{r7}
\sum_{l=3}^{9}z^{-(2l+2)}\bbx^T(\bbW^l-\mathbb{E}\bbW^l)\bby=O_p(\frac{\alpha^2_n}{|z|^4}).
\end{eqnarray}
Thus (\ref{r6}) holds by combining \eqref{r4}, \eqref{r5}, and \eqref{r7}. This concludes the proofs of the desired results.

\subsection{Improvement of Lemmas \ref{0426-1} and  \ref{1212-1h} under Condition \ref{as3}ii)} \label{SecD.2}

\begin{lem}\label{0426-1h}
For any $n$-dimensional unit vectors $\bbx$ and  $\bby$, there exists some positive constant $C$ independent of $l$ such that
	\begin{equation}\label{0506.1h}
	\mathbb{E}\left[\bbx^T(\bbW^l-\mathbb{E}\bbW^l)\bby\right]^2\le C(4c_1l)^{2l}(\min\{\alpha_n^{l-1},d_{\bx}\alpha_n^l, d_{\by}\alpha_n^l\})^2
	\end{equation}
	with $l \geq 1$ some positive integer and $d_{\bbx} = \|\bbx\|_\infty$. 
\end{lem}
\begin{lem}\label{1212-1hh}
For any $n$-dimensional unit vectors $\bbx$ and  $\bby$, there exists some positive constant $C$ independent of $l$ such that
	\begin{equation}
	\left|\mathbb{E}\bbx^T\bbW^l\bby\right| \le C(2c_1l)^l\alpha_n^{l}
	\end{equation}
with $l\geq 2$ some bounded positive integer.
\end{lem}

\subsection{Proof of Lemma  \ref{0426-1h}} \label{SecD.3}
The proof of Lemma  \ref{0426-1h} is a modification of that for Lemma \ref{0426-1}. Thus we highlight only the differences of the technical arguments here. We work directly on the general case allowing for  $\mathbb{E}w_{ii}\neq 0$. In view of \eqref{1212.1h}, we have
\begin{align} \label{1212.1ht}
\nonumber & \mathbb{E} (\bbx^T\bbW^l\bby-\mathbb{E}\bbx^T\bbW^l\bby)^2\\
& =\sum_{1\le i_1,\cdots,i_{l+1},j_1,\cdots,j_{l+1}\le n}\mathbb{E} \Big[\left(x_{i_1}w_{i_1i_2}w_{i_2i_3}\cdots w_{i_li_{l+1}}y_{i_{l+1}}-\mathbb{E}x_{i_1}w_{i_1i_2}w_{i_2i_3}\cdots w_{i_li_{l+1}}y_{i_{l+1}}\right) \nonumber \\
& \quad \times \left(x_{j_1}w_{j_1j_2}w_{j_2j_3}\cdots w_{j_lj_{l+1}}y_{j_{l+1}}-\mathbb{E}x_{j_1}w_{j_1j_2}w_{j_2j_3}\cdots w_{j_lj_{l+1}}y_{j_{l+1}}\right)\Big].
\end{align}
Let $\bbi=(i_1,\cdots,i_{l+1})$ and $\bbj=(j_1,\cdots,j_{l+1})$ be two vectors taking values in $\{1,\cdots,n\}^{l+1}$. For any given vector $\bbi$, we define a graph $\mathcal{G}_{\bbi}$ whose
vertices represent the components of $\bbi$.
Vertices $i_s$ and $i_{s+1}$ of $\mathcal G_{\bbi}$ are connected by undirected edges for $1\le s\le l$. Similarly we can also define graph  $\mathcal{G}_{\bbj}$ corresponding to $\bbj$. It can be seen that $\mathcal{G}_{\bbi}$ is a connected graph, which means that there exists some path from $i_s$ to $i_{s'}$ for any $1\le s\neq s'\le n$.  One should notice that here we allow for $i_s=i_{s+1}$ or $j_s=j_{s+1}$. Such relaxation will affect only the number of pairs $(\bbi,\bbj)$, but will not affect the main arguments of the proof which are similar to the graph arguments for proving Lemma \ref{0426-1}. Thus for each product
\begin{align}\label{0928.3h}
&\mathbb{E} \Big[\left(x_{i_1}w_{i_1i_2}w_{i_2i_3}\cdots w_{i_li_{l+1}}y_{i_{l+1}}-\mathbb{E}x_{i_1}w_{i_1i_2}w_{i_2i_3}\cdots w_{i_li_{l+1}}y_{i_{l+1}}\right) \nonumber \\
&\quad\times \left(x_{j_1}w_{j_1j_2}w_{j_2j_3}\cdots w_{j_lj_{l+1}}y_{j_{l+1}}-\mathbb{E}x_{j_1}w_{j_1j_2}w_{j_2j_3}\cdots w_{j_lj_{l+1}}y_{j_{l+1}}\right)\Big],
 \end{align}
there exists a corresponding graph $\mathcal{G}_{\bbi}\cup\mathcal{G}_{\bbj}$. If $\mathcal{G}_{\bbi}\cup\mathcal{G}_{\bbj}$ is not a connected graph, then the corresponding expectation
\begin{align*}
& \mathbb{E} \Big[\left(x_{i_1}w_{i_1i_2}w_{i_2i_3}\cdots w_{i_li_{l+1}}y_{i_{l+1}}-\mathbb{E}x_{i_1}w_{i_1i_2}w_{i_2i_3}\cdots w_{i_li_{l+1}}y_{i_{l+1}}\right) \\
&\quad \times \left(x_{j_1}w_{j_1j_2}w_{j_2j_3}\cdots w_{j_lj_{l+1}}y_{j_{l+1}}-\mathbb{E}x_{j_1}w_{j_1j_2}w_{j_2j_3}\cdots w_{j_lj_{l+1}}y_{j_{l+1}}\right)\Big]=0.
\end{align*}
This shows that in order to calculate the order of $\mathbb{E} (\bbx^T\bbW^l\bby-\mathbb{E}\bbx^T\bbW^l\bby)^2$, it suffices to consider the scenario of connected graphs $\mathcal{G}_{\bbi}\cup\mathcal{G}_{\bbj}$.

To analyze the term in (\ref{0928.3h}), let us calculate how many distinct vertices are contained in the connected graph $\mathcal{G}_{\bbi}\cup\mathcal{G}_{\bbj}$. We say that  $(i_s,i_{s+1})\in \mathcal{G}_{\bbi}$ is an \textit{efficient edge} if $i_s\neq i_{s+1}$. Since there are at most $2l$ efficient  edges in $\mathcal{G}_{\bbi}\cup \mathcal{G}_{\bbj}$ and $\mathbb{E}w_{ss'}=0$ for $s\neq s'$, in order to get a nonzero value of (\ref{0928.3h}) each efficient  edge in $\mathcal{G}_{\bbi}\cup \mathcal{G}_{\bbj}$ has at least one copy. Thus for each nonzero (\ref{0928.3h}), we have at most $l$ distinct efficient  edges in  $\mathcal{G}_{\bbi}\cup \mathcal{G}_{\bbj}$. Since graph  $\mathcal{G}_{\bbi}\cup \mathcal{G}_{\bbj}$ is connected, there are at most $l+1$ distinct vertices in  $\mathcal{G}_{\bbi}\cup \mathcal{G}_{\bbj}$. Denote by $\mathcal{S}$ the set of all such pairs $(\bbi, \bbj)$. Combining the above arguments, we can conclude that
\begin{align}\label{0928.4h}
(\ref{1212.1h})&=\sum_{(\bbi, \bbj)\in \mathcal S}\mathbb{E} \Big[\left(x_{i_1}w_{i_1i_2}w_{i_2i_3}\cdots w_{i_li_{l+1}}y_{i_{l+1}}
 -\mathbb{E}x_{i_1}w_{i_1i_2}w_{i_2i_3}\cdots w_{i_li_{l+1}}y_{i_{l+1}}\right)
\nonumber \\
& \quad \times \left(x_{j_1}w_{j_1j_2}w_{j_2j_3}\cdots w_{j_lj_{l+1}}y_{j_{l+1}}-\mathbb{E}x_{j_1}w_{j_1j_2}w_{j_2j_3}\cdots w_{j_lj_{l+1}}y_{j_{l+1}}\right)\Big].
\end{align}
For notational simplicity, we denote $j_1,\cdots,j_{l+1}$ as  $i_{l+2},\cdots,i_{2l+2}$ and define $\widetilde \bbi = (i_1,\cdots,i_{l+1},\\j_1,\cdots,j_{l+1})=(i_1,\cdots,i_{2l+2}).$ We also denote $\mathcal{G}_{\bbi}\cup \mathcal{G}_{\bbj}$ as  $\mathcal{F}_{\widetilde \bbi}$ which has at most $l+1$ distinct vertices and $l$ distinct efficient  edges, with each edge having at least two copies. Then it holds that
\begin{align}\label{0928.5h}
& |(\ref{0928.4h})|=\Big|\sum_{\mathcal{F}_{\widetilde \bbi},\, \widetilde\bbi \in \mathcal S}\mathbb{E} \big[(x_{i_1}w_{i_1i_2}w_{i_2i_3}\cdots w_{i_li_{l+1}}y_{i_{l+1}}-\mathbb{E}x_{i_1}w_{i_1i_2}w_{i_2i_3}\cdots w_{i_li_{l+1}}y_{i_{l+1}})\nonumber\\
& \ \ \times (x_{i_{l+2}}w_{i_{l+2}i_{l+3}}w_{i_{l+3}i_{l+4}}\cdots w_{i_{2l+1}i_{2l+2}}y_{i_{2l+2}}-\mathbb{E}x_{i_{l+2}}w_{i_{l+2}i_{l+3}}w_{i_{l+3}i_{l+4}}\cdots w_{i_{2l+1}i_{2l+2}}y_{i_{2l+2}})\big]\Big|\non
&\le\sum_{\mathcal{F}_{\widetilde \bbi},\, \widetilde\bbi \in \mathcal S}\mathbb{E}\big|x_{i_1}w_{i_1i_2}w_{i_2i_3}\cdots w_{i_li_{l+1}}y_{i_{l+1}}x_{i_{l+2}}w_{i_{l+2}i_{l+3}}w_{i_{l+3}i_{l+4}}\cdots w_{i_{2l+1}i_{2l+2}}y_{i_{2l+2}}\big|\non
&\ \ +\sum_{\mathcal{F}_{\widetilde \bbi},\, \widetilde\bbi \in \mathcal S}\mathbb{E}\big|x_{i_1}w_{i_1i_2}w_{i_2i_3}\cdots w_{i_li_{l+1}}y_{i_{l+1}}\big|\mathbb{E}\big|x_{i_{l+2}}w_{i_{l+2}i_{l+3}}w_{i_{l+3}i_{l+4}}\cdots w_{i_{2l+1}i_{2l+2}}y_{i_{2l+2}}\big|.
\end{align}

Observe that each expectation in (\ref{0928.5h}) involves the product of some independent random variables, and $x_{i_1}w_{i_1i_2}w_{i_2i_3}\cdots w_{i_li_{l+1}}y_{i_{l+1}}$ and $x_{i_{l+2}}w_{i_{l+2}i_{l+3}}w_{i_{l+3}i_{l+4}}\cdots w_{i_{2l+1}i_{2l+2}}y_{i_{2l+2}}$ may share some dependency through factors $w_{ab}^{m_1}$ and $w_{ab}^{m_2}$, respectively, for some $w_{ab}$ and nonnegative integers $m_1$ and $m_2$. Thus with the aid of the inequality
\[ \mathbb{E}|w_{ab}|^{m_1}\mathbb{E}|w_{ab}|^{m_2}\le \mathbb{E}|w_{ab}|^{m_1+m_2}, \]
we can further bound (\ref{0928.5h}) as
\begin{align}\label{0928.6h}
(\ref{0928.5h}) & \le 2\sum_{\mathcal{F}_{\widetilde \bbi},\, \widetilde\bbi \in \mathcal S}\mathbb{E}\big|x_{i_1}w_{i_1i_2}w_{i_2i_3}\cdots w_{i_li_{l+1}}y_{i_{l+1}}x_{i_{l+2}}w_{i_{l+2}i_{l+3}}w_{i_{l+3}i_{l+4}}\cdots \nonumber\\
&\quad \times w_{i_{2l+1}i_{2l+2}}y_{i_{2l+2}}\big|.
\end{align}

To facilitate our technical presentation, let us introduce some additional notation. Denote by $\psi(2l+2)$ the set of partitions of the edges $\{(i_1,i_2),(i_2,i_3),\cdots,(i_{2l+1},i_{2l+2}), i_s\neq i_{s+1}, s=1,\cdots,2l+1\}$ and $\psi_{\ge 2}(2l+2)$ the subset of $\psi(2l+2)$ whose blocks have size at least two. Let  $P(\widetilde \bbi)\in\psi_{\ge 2}(2l+2)$ be the partition of $\{(i_1,i_2),(i_2,i_3),\cdots,(i_{2l+1},i_{2l+2}), i_s\neq i_{s+1}, s=1,\cdots,2l+1\}$ that is associated with the equivalence relation $(i_{s_1},i_{s_1+1})\sim (i_{s_2},i_{s_2+1})$ which is defined as if and only if $(i_{s_1},i_{s_1+1})=(i_{s_2},i_{s_2+1})$ or $(i_{s_1},i_{s_1+1})=(i_{s_2+1},i_{s_2})$. Denote by $|P(\widetilde \bbi)|=m$ the number of groups in the partition $P(\widetilde \bbi)$ such that the edges are equivalent within each group. We further denote the distinct edges in the partition $P(\widetilde \bbi)$ as $(s_1,s_2), (s_3,s_4), \cdots, (s_{2m-1},s_{2m})$ and the corresponding counts in each group as  $r_1, \cdots, r_m$, and define $\widetilde \bbs=(s_1,s_2,\cdots,s_{2m})$. For the vertices, let $\phi(2m)$ be the set of partitions of $\{1,2,\cdots,2m\}$ and $Q(\widetilde \bbs)\in\phi(2m)$ the partition that is associated with the equivalence relation $a\sim b$ which is defined as if and only if $s_a=s_b$. Note that $s_{2j-1}\neq s_{2j}$ since in the partition, we consider only the off-diagonal entries (efficient edges) and for diagonal entries, we use the simple inequality $|w_{ii}|\le 1$. Then it holds that
\begin{align}\label{0928.7h}
&\sum_{\mathcal{F}_{\widetilde \bbi},\, \widetilde\bbi \in \mathcal S}\mathbb{E}\big|x_{i_1}w_{i_1i_2}w_{i_2i_3}\cdots w_{i_li_{l+1}}y_{i_{l+1}}x_{i_{l+2}}w_{i_{l+2}i_{l+3}}w_{i_{l+3}i_{l+4}}\cdots w_{i_{2l+1}i_{2l+2}}y_{i_{2l+2}}\big|\non
&\le \sum_{1\le |P(\widetilde \bbi)|=m\le l\atop
P(\widetilde \bbi)\in \psi_{\ge 2}(2l+2)}\sum_{\widetilde \bbi \text{ with partition } P(\widetilde \bbi)\atop
r_1,\cdots,r_m\ge 2}\sum_{Q(\widetilde \bbs)\in\phi(2m)}\sum_{\widetilde \bbs \text{ with partition } Q(\widetilde \bbs)\atop 1\le s_1,\cdots,s_{2m}\le n}|x_{i_1}y_{i_{l+1}}x_{i_{l+2}}y_{i_{2l+2}}|\nonumber\\
&\quad \times \prod_{j=1}^m\mathbb{E}\big|w_{s_{2j-1}s_{2j}}|^{r_j}\non
&\le \sum_{1\le |P(\widetilde \bbi)|=m\le l\atop
P(\widetilde \bbi)\in \psi_{\ge 2}(2l+2)}(\frac{c^2_1\alpha_n^2}{n})^m\sum_{\widetilde \bbi \text{ with partition } P(\widetilde \bbi)\atop
r_1,\cdots,r_m\ge 2}\sum_{Q(\widetilde \bbs)\in\phi(2m)}\sum_{\widetilde \bbs \text{ with partition } Q(\widetilde \bbs)\atop 1\le s_1,\cdots,s_{2m}\le n}|x_{i_1}y_{i_{l+1}}x_{i_{l+2}}y_{i_{2l+2}}|.
\end{align}

It suffices to bound the number of graphs in the above summation. In fact, since the graph is connected there are at most $m+1$ different vertices in the graph. 
Moreover, there are $2l$ edges in the original graph with at most $l$ efficient edges and the partitions corresponding to the edges have at most $(4l)^{2l}$ cases. Thus combining these arguments together we can deduce
\begin{eqnarray}\label{0910.1hh}
&&\sum_{1\le |P(\widetilde \bbi)|=m\le l\atop
P(\widetilde \bbi)\in \psi_{\ge 2}(2l+2)}(\frac{c^2_1\alpha_n^2}{n})^m\sum_{\widetilde \bbi \text{ with partition } P(\widetilde \bbi)\atop
r_1,\cdots,r_m\ge 2}\sum_{Q(\widetilde \bbs)\in\phi(2m)}\sum_{\widetilde \bbs \text{ with partition } Q(\widetilde \bbs)\atop 1\le s_1,\cdots,s_{2m}\le n}|x_{i_1}y_{i_{l+1}}x_{i_{l+2}}y_{i_{2l+2}}|\non
&&\le d^2_{\bbx}d_{\bby}^2(\frac{c^2_1\alpha_n^2}{n})^{l}\sum_{1\le |P(\widetilde \bbi)|=m\le l\atop
P(\widetilde \bbi)\in \psi_{\ge 2}(2l+2)}\sum_{\widetilde \bbi \text{ with partition } P(\widetilde \bbi)\atop
r_1,\cdots,r_m\ge 2}\sum_{Q(\widetilde \bbs)\in\phi(2m)}\sum_{\widetilde \bbs \text{ with partition } Q(\widetilde \bbs)}1\non
&&\le d^2_{\bbx}d_{\bby}^2(\frac{c^2_1\alpha_n^2}{n})^{l} (4l)^{2l}n^{l+1}\non
&&\le (4c_1l)^{2l}n\alpha_n^{2l}d^2_{\bbx}d_{\bby}^2.
\end{eqnarray}
Therefore, we can establish the simple upper bound that
\begin{equation}
	\mathbb{E} \left[\bbx^T(\bbW^l-\mathbb{E}\bbW^l)\bby\right]^2\le C(4c_1l)^{2l}n\alpha_n^{2l}d^2_{\bbx}d_{\bby}^2.
	\end{equation}

For the other upper bounds $C(4c_1l)^{2l}d^2_{\bbx}\alpha_n^{2l}$, $C(4c_1l)^{2l}d^2_{\bby}\alpha_n^{2l}$, and $C(4c_1l)^{2l}\alpha_n^{2l-2}$, the arguments are similar to those for the proof of Lemma \ref{0426-1}. The crucial steps are considering the impact of $|x_{i_1}y_{i_{l+1}}x_{i_{l+2}}y_{i_{2l+2}}|$ from \eqref{0928.12h} to \eqref{0928.15h.new}. For our case, we can directly prove the desired bounds $C(4c_1l)^{2l}d^2_{\bbx}\alpha_n^{2l}$, $C(4c_1l)^{2l}d^2_{\bby}\alpha_n^{2l}$, and $C(4c_1l)^{2l}\alpha_n^{2l-2}$ by combining  the left hand side of \eqref{0910.1hh} with the arguments from \eqref{0928.12h} to \eqref{0928.15h.new}. This completes the proof of Lemma  \ref{0426-1h}.

\subsection{Proof  of Lemma \ref{1212-1hh}} \label{sec:prof-lm8}
Similar to the proof of Lemma \ref{1212-1h}, the proof of Lemma \ref{1212-1hh} is a direct modification of that of Lemma \ref{0426-1h}. Thus we omit it for brevity.

\subsection{Proof sketch for results in Section \ref{sec: app2}} \label{sec:prof-app2}	
	By calculating the variance of $\widehat p$, we have
	\begin{equation}\label{2rv3}
	\widehat p=p+O_p\Big(\frac{\sqrt{p(1-p)}}{n}\Big)=p\left[1+O_p\Big(\frac{\sqrt{1-p}}{n\sqrt p}\Big)\right].
	\end{equation}
	Then the mean and variance of $\bbv_1^T\bbW^2\bbv_1$ in \eqref{2rv2} can be estimated as
	\begin{equation}\label{2rv4}
	\widehat{\bbv_1^T\mathbb{E}\bbW^2\bbv_1}=n\widehat p(1-\widehat p) \ \text{ and } \  \widehat{\var(\bbv_1^T\bbW^2\bbv_1)}=\widehat p(1-\widehat p)\left[2(n-1)+\widehat{p}^3+(1-\widehat p)^3\right],
	\end{equation}
	receptively.
	By Theorem \ref{0619-1}, \eqref{2rv3}, and \eqref{2rv4}, direct calculations show that if $n^{-1}\ll p <1$, then it holds that 
	$$\lambda_1-t_1=O_p\Big(\frac{1}{\sqrt{np}}+\sqrt p\Big),$$
	$$\frac{\widehat{\bbv_1^T\mathbb{E}\bbW^2\bbv_1}}{\sqrt{\widehat{\var(\bbv_1^T\bbW^2\bbv_1)}}}=\frac{{\bbv_1^T\mathbb{E}\bbW^2\bbv_1}}{\sqrt{{\var(\bbv_1^T\bbW^2\bbv_1)}}}+o_p(1).$$
	Thus if the conditions of Corollary \ref{1109-1} hold, by \eqref{2rv1} we can obtain 
	\begin{eqnarray}\label{2rv5}
	\frac{2\lambda_1^2\left(\bbv_1^T\widehat\bbv_1-1\right)+\widehat{\bbv_1^T\mathbb{E}\bbW^2\bbv_1}}{\left[\widehat{\var(\bbv_1^T\bbW^2\bbv_1)}\right]^{1/2}}&\toD& N(0,1).
	\end{eqnarray}
	Since $\bbv_1 = n^{-1/2}\bone$ under the null hypothesis, the above results together with
	\eqref{2rv5} ensure that under the null hypothesis, statistic $T_n$ is asymptotically standard normal.
	
	Next we consider the case of alternative hypothesis. It can be derived that the leading eigenvalue and eigenvector  take the following forms
	$$d_1=\frac{1} {2}\left[np+n_1(q-p)+\left(n^2p^2+2n_1(2n_1-n)p(q-p)+n_1^2(q-p)^2\right)^{1/2}\right]$$ and $\bbv_1 =(\bbv_{1,1}^T,\bbv_{1,2}^T)^T$, where $\bbv_{1,1}$ is an $n_1$-dimensional vector with all entries being  
	$$\frac{(n-n_1)p}{\sqrt{(n-n_1)(d_1-n_1q)^2+n_1(n-n_1)^2p^2}}$$
	and $\bbv_{1,2}$ is an  ($n-n_1$)-dimensional vector with all entries being $$\frac{d_1-n_1q}{\sqrt{(n-n_1)(d_1-n_1q)^2+n_1(n-n_1)^2p^2}}.$$ 
	With some direct calculations, we can show that under the alternative hypothesis, 
	\begin{equation}\label{2rv12}
	n^{-1/2}\bone^T\bbv_1=\frac{(n-n_1)(d_1-n_1(q-p))}{\sqrt{n((n-n_1)(d_1-n_1q)^2+n_1(n-n_1)^2p^2)}}.
	\end{equation}
	Since $n_1 = o(n)$, $n^{-1}\ll p <q$, and $p \sim q$, by the Taylor expansion we can deduce 
	$$d_1=np+n_1^2(q-p)\frac{4p+5(q-p)}{4np}+O(\frac{n_1^3(q-p)^2}{n^2p})$$ 
	and 
	\begin{align*}
	&\sqrt{n((n-n_1)(d_1-n_1q)^2+n_1(n-n_1)^2p^2)}\non
	&=\sqrt{n(n-n_1)}(d_1-n_1q)+\frac{n_1\sqrt n(n-n_1)^2p^2}{2\sqrt{(n-n_1)}(d_1-n_1q)}+\frac{n_1^2p}{4}+O(\frac{n_1^3p}{n})\non
	&=\sqrt{n(n-n_1)}\left[np-n_1(q-p)-\frac{n_1p}{2}-\frac{n_1^2p}{4n}+n_1^2(q-p)\frac{4p+5(q-p)}{4np}+O(\frac{n_1^3p}{n^2})\right].
	\end{align*}
	
	Substituting the above two equations into \eqref{2rv12} yields 
	\begin{align}\label{2rv12-a}
	n^{-1/2}\bone^T\bbv_1 & -1
	=\sqrt{\frac{n-n_1}{n}} \Big[1+\frac{n_1}{2n}-\frac{n_1^2}{4n^2}-\frac{n_1^2(q-p)^2}{n^2p^2}-n_1^2(q-p)\frac{4p+5(q-p)}{4n^2p^2}\non
	&\quad+O(\frac{n_1^3}{n^3}+\frac{n_1^3(q-p)^2}{n^3p^2})\Big]-1\non
	&=-\frac{n_1^2(q-p)^2}{n^2p^2}-n_1^2(q-p)\frac{4p+5(q-p)}{4n^2p^2}+O(\frac{n_1^3}{n^3}+\frac{n_1^3(q-p)^2}{n^3p^2}).
	\end{align}
	If the conditions of Corollary \ref{1109-1} hold, by \eqref{2rv1} we have
	\begin{eqnarray}
	\frac{2\lambda_1^2\left({\bbv}_1^T\widehat\bbv_1-1\right)+{\bbv_1}^T\mathbb{E}\bbW^2{\bbv_1}}{\left[\var({\bbv_1}^T\bbW^2{\bbv_1})\right]^{1/2}}&\toD& N(0,1).
	\end{eqnarray}
	This entails that
	\begin{equation}\label{2rv13}
	{\bbv}_1^T\widehat\bbv_1-1=O_p(\frac{{\bbv_1}^T\mathbb{E}\bbW^2{\bbv_1}+\left[\var({\bbv_1}^T\bbW^2{\bbv_1})\right]^{1/2}}{t_1^2})=O_p(\frac{1}{np}),
	\end{equation}
	where the last step is obtained by directly calculating the mean and variance of $\bbv_1^T\bbW^2\bbv_1$ and noting that $t_1\sim np$.
	Since $\bbv_1,\cdots, \bbv_n$ form an orthonormal basis, it follows from \eqref{2rv13} that 
	\begin{equation}\label{2rv14}
	\sum_{j=2}^n({\bbv}_j^T\widehat\bbv_1)^2=1-(\bbv_1^T\widehat\bbv_1)^2=O_p(\frac{1}{np}).
	\end{equation}
	
	Similarly, by \eqref{2rv12} and the assumptions of  $n_1=o(n)$ and $q\sim p$, we can deduce 
	\begin{equation}\label{2rv15}
	\sum_{j=2}^n(n^{-1/2}\bone^T{\bbv}_j)^2=O(\frac{n_1^3}{n^3}+\frac{n_1^3(q-p)^2}{n^3p^2}).
	\end{equation}
	Then it follows from \eqref{2rv12}, \eqref{2rv13}, and \eqref{2rv15} that 
	\begin{align}
	&n^{-1/2}\bone^T\widehat\bbv_1-1=n^{-1/2}\bone^T\bbv_1\bbv^T_1\widehat\bbv_1-1+n^{-1/2}\bone^T\sum_{j=2}^n\bbv_j\bbv^T_j\widehat\bbv_1\non
	&\quad=-\Big[\frac{n_1^2(q-p)^2}{n^2p^2}+n_1^2(q-p)\frac{4p+5(q-p)}{4n^2p^2}\Big]+O_p\Big[\frac{n_1^3}{n^3}+\frac{n_1^3(q-p)^2}{n^3p^2}+\frac{1}{np}\Big].
	\end{align}
	Under the alternative hypothesis, it can be shown that the estimators in \eqref{2rv4-a} are of orders
	$n\widehat p(1-\widehat p)=O_p(np)$ and $\widehat p(1-\widehat p)\left[2(n-1)+\widehat{p}^3+(1-\widehat p)^3\right]=O_p(np)$, respectively, and in addition, $t_1\sim np$. Therefore, if the conditions of Corollary \ref{1109-1} holds and $\frac{n_1^2(q-p)^2}{np}+\frac{n_1^2(q-p)}{n}\gg 1$, with probability tending to one we have
	$$T_n\rightarrow -\infty,$$
	which means that the power can tend to one asymptotically. This concludes the proof sketch for the results in Section \ref{sec: app2}.
	
\end{document}